\documentclass{amsart}

\usepackage{amsmath,amssymb,amsfonts,amscd}
\usepackage[all,cmtip]{xy}
\usepackage{enumerate}
\usepackage{ amssymb, latexsym, amsmath}

\usepackage{lingmacros}
\usepackage{tree-dvips}
\usepackage{amssymb,amsfonts,amsthm,amsmath}

\usepackage{mathrsfs}
\usepackage[english]{babel}
\usepackage[all,cmtip]{xy}
\usepackage{tikz}
\usepackage{tikz-cd}
\usetikzlibrary{cd,quotes}
\usetikzlibrary{calc,babel}
\usepackage{mathtools}

\usepackage{tensor}
\usepackage{csquotes}
\usepackage[breaklinks]{hyperref}
\usetikzlibrary{arrows,chains,matrix,positioning,scopes}

\usepackage{fullpage}
\usepackage{url}
\usepackage{hyperref}
\usepackage{ marvosym }
\usepackage{amssymb}

\numberwithin{equation}{section}

\theoremstyle{plain}
\newtheorem{theorem}{Theorem}[section]

\newtheorem{prop}[theorem]{Proposition}
\newtheorem{lemma}[theorem]{Lemma}

\newtheorem{cor}[theorem]{Corollary}

\theoremstyle{definition}
\newtheorem{defn}[theorem]{Definition}

\newcommand{\CP}{\mathbb{CP}}
\newcommand{\R}{\mathbb{R}}
\newcommand{\N}{\mathbb{N}}
\newcommand{\C}{\mathbb{C}}
\newcommand{\Z}{\mathbb{Z}}
\newcommand{\RP}{\mathbb{RP}}

\newcommand{\T}{\mathcal{T}}
\newcommand{\F}{\mathcal{F}}
\newcommand{\E}{\mathcal{E}}

\theoremstyle{remark}
\newtheorem{remark}[theorem]{Remark}

\makeatletter

\newcommand{\Rmnum}[1]{\expandafter\@slowromancap\romannumeral #1@}
\makeatother

\newcommand{\Id}{\operatorname{Id}}

\newcommand{\n}{\mathfrak{n}}

\newcommand{\p}{\mathfrak{p}}

\newcommand{\m}{\mathfrak{m}}

\newcommand{\M}{\overline{\mathcal{M}}}

\begin{document}
    
    \begin{abstract}
    Let $L_0$ and $L_1$ be a pair of transversely intersecting closed Lagrangians in a closed symplectic manifold, equipped with a relative-spin structure with respect to a common background class. We construct a curved, filtered, $A_\infty$-bimodule structure over the integral Fukaya algebras of $L_0$ and $L_1$ as given in \cite{rabah2024fukaya}, \cite{rabah2025fukaya}; on the free module generated by the intersection points of $L_0$ and $L_1$, over the integral Novikov ring. When both Lagrangians are unobstructed, we define a Lagrangian Floer homology of the pair over the integral Novikov ring, extending the work of Fukaya-Oh-Ohta-Ono in \cite{fukaya2013lagrangian} to general closed symplectic manifolds. We use this framework to bound the number of intersection points between a Lagrangian and its transverse Hamiltonian image. Such bound incorporates the torsion subgroup in the integral homology of the Lagrangian. In return, it improves the characteristic-zero Betti number bound.
\end{abstract}

		\title
		{A Lagrangian Floer Theory over the Integers and Finite Fields}
		\author{Mohamad Rabah}
		\maketitle
		\tableofcontents

\section{Introduction}

Let $(X,\omega)$ be a closed symplectic manifold and $L\subset X$ a closed Lagrangian. For $f:L\rightarrow\R$ a Morse function on $L$ and $\epsilon>0$ small enough, $L$ and the graph of $\epsilon df$ intersect transversely in a Weinstein neighborhood of $L$ in $X$. The integral Morse inequalities then provide us with a lower bound on the number of intersection points between $L$ and the graph of $\epsilon df$, involving the Betti numbers of $L$ and the torsion in $H_*(L;\Z)$. In this paper, we study the $\textit{Lagrangian intersection problem}$ of $L$ with its Hamiltonian image by any arbitrary degenerate Hamiltonian in a general closed symplectic manifold.

In their books \cite{Fukaya2009}, \cite{fukaya2010lagrangian}, Fukaya-Oh-Ohta-Ono developed a filtered $A_\infty$-approach to this intersection problem and in \cite{fukaya2013lagrangian}, the authors constructed an integral version in spherically-positive closed symplectic manifolds. In this paper, we extend the author's thesis \cite{rabah2025fukaya} on the integral Fukaya algebra for a single Lagrangian to a pair of transversely intersecting Lagrangians. In particular and following Floer's original approach \cite{floer1988morse}, for such a pair of transversely intersecting Lagrangians, we construct a chain complex generated by their intersection points, where its homology is what we refer to as $\textit{Lagrangian Floer homology}$. A new feature of our construction is that it is done over Novikov rings with coefficients in $\Z$ or finite fields $\Z/p$ for $p$ prime, without assuming spherical-positivity on the ambient symplectic manifold as in \cite{fukaya2013lagrangian}.

To be more precise, let $(X,\omega)$ be a closed symplectic manifold and $L\subset X$ be a closed, connected relatively-spin Lagrangian. We say $L$ is unobstructed over $\Lambda_0$ the Novikov ring with integer coefficients (c.f. Definition \ref{UnobstructedDefn}), where $\Lambda_0$ is given by 
\[
\Lambda_0:=\{\sum_ia_iT^{\lambda_i}e^{n_i}:a_i,n_i\in\Z, 0\leq\lambda_1\leq\lambda_2\leq\dots,\lim_{i\rightarrow\infty}\lambda_i=\infty,\#\{i:a_i\neq0, \lambda_i\leq E\}<\infty,\forall E>0\},
\]
where $T,e$ are formal variables of degrees $0,2$ respectively; if we can find a degree-$1$ cochain $b$ of positive-valuation solving the Maurer-Cartan equation
\[
\sum_{k\geq0}\m_k(b,\dots, b)=0.
\]
Denote by $\iota:L\hookrightarrow X$ the canonical inclusion and let $t_L$ be the minimum number generators of the torsion subgroup of $\bigoplus_j H_j(L;\Z)$. Under the assumption that $L$ is unobstructed over the Novikov ring with integer coefficients and that $\iota_*:H_*(L;\Z)\rightarrow H_*(X;\Z)$ admits a left-inverse, we show (c.f. Theorem \ref{IntergralArnoldConjI}) that for any given Hamiltonian diffeomorphism $\phi:X\rightarrow X$ of $X$ such that $L$ and $\phi(L)$ intersect transversely, their number of intersection points is bounded from below by
\[
\#(L\cap \phi(L))\geq \sum_j \mathrm{rank}_\Z H_j(L;\Z)+2t_L.
\]
The condition on $\iota_*$ as stated above, is a bit restrictive. Nevertheless they are realized in the case of the diagonal embedding of a symplectic manifold to its product symplectic manifold. On the other hand, to give bounds on the number of intersection points $\#(L\cap\phi(L))$ without assuming that $\iota_*:H_*(L;\Z)\rightarrow H_*(X;\Z)$ admits a left-inverse, we develop a Lagrangian Floer theory over finite fields.

For $p$ any prime number and if we denote by $\mathbb F_p$ the field of $p$-elements, we give $H^*(L;\mathbb F_p)$ the structure of an $A_\infty$-algebra (c.f. Theorem \ref{A_inftyoverF}) as in \cite{Fukaya2009}, \cite{fukaya2013lagrangian}, \cite{fukaya2013lagrangian}, \cite{rabah2024fukaya},\cite{rabah2025fukaya}. As a consequence, if $L$ is unobstructed over the Novikov ring with coefficients in $\mathbb F_p$ (c.f. Definition \ref{unobstructedoverFdefn}) and we denote by $A_p:=\sum_j\dim_{\mathbb{F}_p}H_j(L;\mathbb F_p)$ and $B_p$ to be the total dimension of the kernel of $\iota_*:H_*(L;\mathbb F_p)\rightarrow H_*(X;\mathbb F_p)$ then (c.f. Theorem \ref{LagIntFThm}), 
\[
\#(L\cap \phi(L))\geq A_p-2B_p.
\]
This result is Theorem $I$ of \cite{Fukaya2009} and Theorem $15.6$ of \cite{fukaya2013lagrangian}, yet we prove it in full generality with no restrictions on the ambient symplectic manifold. In fact, in our framework, we show (c.f. Proposition \ref{UnobstroverZimpliesF}) that if $L$ is unobstructed over the Novikov ring with integer coefficients then, it is unobstructed over the Novikov ring with coefficients in $\mathbb F_p$ for any prime $p$. As a consequence (c.f. Corollary \ref{corLagIntoverF}), we can give the following stronger lower bound 
\[
\#(L\cap \phi(L))\geq \max_p A_p-2B_p.
\]

Another feature of our framework that we would like to highlight, is that we realize Remark $14.3$ \cite{fukaya2013lagrangian} and prove convergence of the associated energy filtration spectral sequence without using Assumption $15$ of \cite{fukaya2013lagrangian}.

\begin{remark}
    In principle, the methods developed in this paper could be extended to prove integral or even a finite field version of Theorem $J$ of \cite{Fukaya2009}, as was done in \cite{fukaya2013lagrangian} without the need of spherical-positivity of the ambient symplectic manifold. This would require using the Hamiltonian-perturbed moduli spaces as considered in \cite{fukaya2013displacement}, in the proof of Theorem \ref{HamiltonianInvariance}, at the cost of adapting the $\textit{framing construction}$ of \cite{rezchikov2022integral} or \cite{BaiXu2026Foundations}, using a Hamiltonian perturbed $2$-form. Nevertheless, we did not pursue this here as we were not able to find any examples in which such integral or finite field version of Theorem $J$ of \cite{Fukaya2009} in a general closed symplectic manifold, would yield a stronger lower bounds on the displacement energy of some Lagrangian. On the other hand, the coefficient dependent $\textit{phenomena}$ showed in \cite{KawamotoShelukhin2024} on Spectral invariants, suggests that an integral Lagrangian Floer theory may contain additional quantitative information.
\end{remark}

\subsection{Organization}
As mentioned above, this work is a continuation of \cite{rabah2024fukaya}, \cite{rabah2025fukaya}. Thus, in order to reduce redundancy and repetition, we will use the arguments and results of \cite{rabah2024fukaya}, \cite{rabah2025fukaya} freely and refer the reader to the relevant Sections therein.

In $\textbf{Section 2}$, we construct normally complex Global Kuranishi charts in the sense of \cite{bai2022arnold}, for the compactified moduli space of marked pseudo-holomorphic strips. In $\textbf{Section 3}$, we construct a filtered open-closed map from the Hochschild chains of the Fukaya algebra of a single Lagrangian, to the Morse complex of the ambient symplectic manifold. In $\textbf{Section 4}$, we use two-paged metric ribbon trees, as in \cite{ma2010geometric}, \cite{ma2015quilted}, \cite{charest2015floer} and \cite{venugopalan2020fukaya}, to associate for a pair of transversely intersecting Lagrangians, the structure of an $A_\infty$-bimodule on the algebra generated by their intersection points.

In $\textbf{Section 5}$, we use bounding cochains to define Floer homology over the integral Novikov ring for a single Lagrangian and a pair of transversely intersecting Lagrangians. Using the constructions of Section $3$, we prove the convergence of the energy spectral sequence. Finally, in $\textbf{Section 6}$, we prove the lower bounds on $\#(L\cap \phi(L))$ as advertised above, and finish this paper by giving some examples.

\subsection*{Acknowledgment}
 We carried out most of this work during my doctoral studies at Stony Brook
University. We thank the Department of Mathematics and the Simons Center
for Geometry and Physics for their hospitality during those years. We are deeply indebted to my advisors, Professors Kenji Fukaya and Mark McLean. Most, if not all, of the ideas and results in this paper grew
out of their guidance and sustained support. We are also grateful to Professor Ivan Smith for bringing to our attention the papers used to formulate the examples given in Section \ref{AppSection}, to Professor Guangbo Xu for always taking the time to answer my questions and to Professor Zhengyi Zhou for his interest in this work. Finally, I thank my girlfriend, Tracy Zhang Ke, for her constant encouragement and support.\\
AI/LLM were $\textbf{not}$ used in the preparation of this manuscript.

\section{Preliminary Constructions}
Let $(X, \omega)$ be a closed symplectic manifolds and let $L_0, L_1$ be a pair of transversely intersecting closed connected Lagrangians. 

\begin{defn}[Relatively-spin as a pair]
    We say the pair $L_0, L_1$ is relatively-spin if $L_i$ is oriented and there exists an $st\in H^2(X; \Z/2)$ such that $st|_{L_i}=w_1(TL_i)$ for $i=0,1$.
\end{defn}

\begin{remark}
    Note that, the notion of relatively-spin as a pair is stronger than requiring both Lagrangians to be relatively-spin. That is, we require $L_i$ to be relatively-spin with respect to the same class $st\in H^2(X; \Z/2)$.
\end{remark}

 For $p, q\in L_0\cap L_1$ two intersection points and following \cite{Fukaya2009}, we introduce the following notation:

\begin{defn}
    Let $\pi_2(p, q)\equiv \pi_2(L_0, L_1;p,q)$ be the set of homotopy classes of smooth maps $u:[0,1]_s\times [0,1]_t\rightarrow X$ such that:
    \begin{enumerate}
        \item $u(0,t)\equiv p$ and $u(1, t)\equiv q$ for all $t\in[0,1]$,
        \item $u(s, 0)\in L_0$ and $u(s, 1)\in L_1$ for all $s\in[0,1]$.
    \end{enumerate}
\end{defn}

We will use the moduli space of pseudo-holomorphic strips bounding $L_0$ and $L_1$ in our constructions. 

\begin{defn}[Marked Strips]\label{stripsdefn}
            A marked strip is a tuple $((\Sigma, \partial\Sigma), \underline{z}, \underline{w})$ representing a marked genus-zero bordered Riemann surface with $(k+2)$-cyclically ordered boundary marked points $\underline{z}$ and $l$-interior marked points $\underline{w}$, such that all marked points are non-nodal points of $(\Sigma, \partial\Sigma)$. We require the existence of $2$-distinguished boundary marked points $z_-, z_+\in\underline{z}$, respecting the cyclic order on $\underline{z}$. 
        \end{defn}

      We recall that the group of automorphisms of $(\Sigma, \partial\Sigma, \underline{z}, \underline{w})$, which we denote by $\operatorname{Aut}(\Sigma, \partial\Sigma, \underline{z}, \underline{w})$ is the group of homeomprhisms $\psi:(\Sigma, \partial\Sigma, \underline{z}, \underline{w})\rightarrow(\Sigma, \partial\Sigma, \underline{z}, \underline{w})$ fixing all marked points such that on each irreducible component of $(\Sigma, \partial\Sigma)$, $\psi$ restricts to a biholomorphism. We denote by $[\Sigma, \partial\Sigma, \underline{z}, \underline{w}]$ the set of equivalence classes given by $\operatorname{Aut}(\Sigma, \partial\Sigma, \underline{z}, \underline{w})$. Similarly, if $u:(\Sigma,\partial\Sigma, \underline{z}, \underline{w})\rightarrow X$ is a continuous map, we denote by $\operatorname{Aut}((\Sigma,\partial\Sigma, \underline{z}, \underline{w}),u)\equiv \operatorname{Aut}(u)$ the subgroup of $\operatorname{Aut}(\Sigma, \partial\Sigma, \underline{z}, \underline{w})$ such that $u\equiv u\circ\psi$ and denote by $[(\Sigma, \partial\Sigma, \underline{z}, \underline{w}),u]$ its set of equivalence classes.

      Now suppose that $(\Sigma, \partial\Sigma, \underline{z}, \underline{w})$ is marked strip, as in Definition \ref{stripsdefn}, and note that by the definition of a nodal bordered Riemann surface, it follows that $\partial\Sigma$ is connected and homeomorphic to a wedge of circles. Let $T_D$ be the tree associated to $(\Sigma, \partial\Sigma)$ where each vertex of $T_D$ corresponds uniquely to a disk-component of $(\Sigma, \partial\Sigma)$ and each edge connecting two vertices to $T_D$ corresponds to a nodal point $(\Sigma, \partial\Sigma)$ joining such disk-components. As $\partial\Sigma$ is connected, it follows that $T_D$ is a connected tree and hence the unique minimal path joining the vertex containing the marked point $z_-$ to the vertex containing the marked point $z_+$ is a $\textit{chain}$ of disk-components which we denote by $D_1, \dots, D_m$, where $m\geq1$. We also need to distinguish the nodal points of $(\Sigma, \partial\Sigma)$ corresponding to edges to $T_D$. First, we set $p_1^-\equiv z_-$ and $p_m^+=z_+$. If $i>1$, we set $p_i^-$ to be the preimage under normalization of the nodal point corresponding to the edge in $T_D$ joining $D_{i-1}$ to $D_i$. If $i<m$, we set $p_i^+$ to be the preimage under normalization of the nodal point corresponding to the edge in $T_D$ joining $D_i$ to $D_{i+1}$.  

\begin{defn}[Strip-Components]\label{stripComponentDefn}
    The disk component of $(\Sigma, \partial\Sigma)$ corresponding to some $D_i$ as above, is called a strip-component. 
\end{defn}

For a given $\omega$-tame almost complex structure $J$ on $(X, \omega)$, we introduce the following Definition:

\begin{defn}[$J$-strips]\label{J-stripsDefn}
    A $J$-strip is a tuple  $((\Sigma,\partial\Sigma, \underline{z},\underline{w} ),u)$ such that:
    \begin{enumerate}
        \item $(\Sigma, \partial\Sigma, \underline{z}, \underline{w})$ is a marked strip as in Definition \ref{stripsdefn},
        \item $\operatorname{Aut}((\Sigma, \partial\Sigma, \underline{z}, \underline{w}), u)$ is of finite order,
        \item $u$ is a $J$-holomorphic curve, that is $\Bar{\partial}_J u=0$, with boundary conditions on $L_0, L_1$ given as follows:\\
        On each strip-component $\Sigma_i$ of $\Sigma$ corresponding to $D_i$, $\partial\Sigma_i\setminus\{p_i^-, p_i^+\}$ is of the form of two arcs each mapped by $u$ to a distinct choice of $L_0$ or $L_1$, respecting the cylic order to boundary marked points and the continuity of $u$. We also require that all disk-components of $(\Sigma, \partial\Sigma)$ that are not strip-components have boundary conditions on exactly one of the $L_0$ and $L_1$.
    \end{enumerate}
\end{defn}

\begin{defn}[Moduli Space of $J$-strips]\label{ModuliJ-stripsDefn}
    We denote by $\M_{k+2, l}(L_0, L_1, p,q; J,\alpha)$ the set of equivalence classes $[(\Sigma, \partial\Sigma, \underline{z}, \underline{w}),u]$ of $J$-strips where:
    \begin{enumerate}
        \item $u$ represents the class $\alpha$, that is $[u]=\alpha$, and
        \item we require that $u(z_-)\equiv p$ and $u(z_+)\equiv q$.
    \end{enumerate}
\end{defn}

We give $\M_{k+2, l}(L_0, L_1,p,q; J, \alpha)$ the $\textit{Gromov Convergence}$ topology, as in \cite{fukaya1996arnold} or \cite{mcduff2012j}, making it a compact Hausdorff space.

        The goal of this section is to give the moduli space of $J$-strips, as in Definition \ref{ModuliJ-stripsDefn}, admits a normally complex Global Kuranishi chart representation. Before we state the main Theorem of this section, we start by recalling the relevant definitions following \cite{abouzaid2021complex}, \cite{bai2022arnold}, \cite{rezchikov2022integral}, \cite{hirschi2024global}, \cite{abouzaid2024gromov} and \cite{rabah2024fukaya}.

	\begin{defn}[Global Kuranishi Chart]\label{GKCdefn}
		Let $\mathcal X$ be a metrizable space. A Global Kuranishi chart of $\mathcal X$ or a Global Kuranishi chart representation of $\mathcal X$ is a tuple $(G, \T, \E, s)$ such that:
		\begin{enumerate}
			\item $G$ is a compact Lie group, thought of as the $\textit{symmetry group}$.
			\item $\T$ and $\E$ are smooth $G$-manifolds, possibly with boundaries and corners, such that $G$ acts locally-linear and with finite stabilizers on $\T$ and $\E$. We call $\T$ the $\textit{thickening}$. 
			\item $\E\rightarrow\T$ is a $G$-vectorbundle, thought of as the $\textit{obstruction bundle}$.
			\item $s:\T\rightarrow\E$ is a $G$-equivariant contionuous section, which we call the $\textit{Kuranishi}$ section or map.
		\end{enumerate}
        Where we have an indentification $s^{-1}(0)/G\cong \mathcal X$ as orbispaces.
        \end{defn}
Observe that, from item $2$ of the above Definition, it follows that, after quotienting by the $G$-action, $\mathcal X$ is identified, at the level of orbispaces, with the zero locus of the derived orbifold chart $(U, E,S)$. With this in mind and referring to Definition $2.1$ of \cite{bai2022arnold} or Definition $1.9$ of \cite{rabah2024fukaya}, we say $(G, \T, \E, s)$ is a normally complex Global Kuranishi chart of $\mathcal X$, if the induced derived orbifold chart $(U,E,S)$ is normally complex.

We also recall that for a given metrizable space $\mathcal X$, we may have more than one Global Kuranishi chart representations of $\mathcal X$. Following the terminology as in \cite{abouzaid2021complex}, we say such Global Kuranishi charts are $\textit{equivalent}$ and in fact, are related to each others by the following set of operations:\\
Suppose that $(G, \T, \E, s)$ is a Global Kuranishi chart of $\mathcal{X}$,
\begin{enumerate}
    \item (Germ Equivalence): If $s^{-1}(0)\subseteq U\subseteq\T$ is a $G$-invariant open set of $s^{-1}(0)$ then, $(G, U, \E|_U, s|_U)$ is a Global Kuranishi chart of $\mathcal X$ equivalent to $(G, \T, \E, s)$.
    \item (Group Enlargement): If $G'$ is a compact Lie group and $\pi:P\rightarrow\T$ is a $G$-equivariant principal $G'$-bundle then, $(G'\times G, P, \pi^*\E, \pi^*s)$ is a Global Kuranishi chart of $\mathcal X$ equivalent to $(G, \T, \E, s)$.
    \item (Stabilization): If $\pi: V\rightarrow\T$ is a $G$-equivariant vectorbundle then, $(G, V, \pi^*\E\oplus\pi^*V, \pi^*s\oplus\tau_V)$ is a Global Kuranishi chart of $\mathcal X$ equivalent to $(G, \T, \E, s)$, where $\tau_V:V\rightarrow \pi^*V$ is the tautological section.
\end{enumerate}

        \begin{theorem}\label{GlobalChartThmI}
             Let $k,l\geq0$ be integers, $p,q\in L_0\cap L_1$ be two distinct intersection points and $\beta\in\pi_2(p,q)$. The moduli space $\overline{\mathcal M}_{k+2, l}(L_0, L_1, p,q; J, \beta)$ as given in Definition \ref{ModuliJ-stripsDefn}, admits a normally complex Global Kuranishi chart representation $(G, \T, \E, s)$.
        \end{theorem}

\subsection{Framing}
To avoid notational clutter, we denote by $\mathbf{L}\equiv L_0\cup L_1$.
		\begin{lemma}\label{Existinceof2-form}
			There exists a precompact neighborhood $U$ of $\mathbf L$ and a closed $2$-form $\Omega\in\Omega^2(X)$ supported on the complement of $U$ such that:
			\begin{enumerate}
				\item $\Omega$ is non-degenerate and tame with $J$ on its support and
				\item there exists a positive integer $N>0$ so that $N[\Omega]\in H^2(X, \mathbf L; \mathbb Z)$ is an integral lift of $[\omega]\in H^2(X, \mathbf L; \mathbb R)$. 
			\end{enumerate} 
		\end{lemma}
		
		\begin{proof}
			The proof is essentially given in Lemma $3.5$ of \cite{rabah2024fukaya}, after replacing Weinstein tubular neighborhood of a single Lagrangian by that of two Lagrangians $L_1$ and $L_2$, as in our setting, and smoothing the parts where the tubular neighborhoods intersect. To be more preciese and for the sake of simplicity, we assume that $L_1$ and $L_2$ intersect transversally. Let $g:=\omega(., J.)$ be the Riemannian metric on $X$ and fix an $\epsilon>0$. Denote by $\mathbb D_\epsilon^*(L_i)$ a Weinstein tubular neighborhood around $L_i$, for $i=1,2$. At a point $x$ of transversal intersection of $L_1\cap L_2$, we consider a $\textit{plumbing}$ neighborhood of $T^*_xL_1\cup T^*_x L_2$, namely a smoothing of $(U_1\times\mathbb D^*_{\frac{\epsilon}{3}, x}L_1)\cup (U_2\times\mathbb D^*_{\frac{\epsilon}{3}, x}L_2)$, where $U_i$ is a small enough trivializing neighborhood of $T^*L_i$ around $x\in L_i$, for $i=1,2$, and hence the exitestence of a $U$ as claimed. Now using the same bump-function argument as given in the proof of Lemma $3.5$ of \cite{rabah2024fukaya}, we get the desired result. 
		\end{proof}
		
		\begin{lemma}\label{LineBundleLemma}
			There exits a Hermitian line bundle $(H, \nabla)\rightarrow X$ of curvature form $F_\nabla=-2\pi i\Omega$, together with a trivialization $\tau:H|_U\cong U\times \mathbb C$ for some open neighborhood $U$ of $\mathbf L$.
		\end{lemma}
		
		\begin{proof}
			This follows from the proof of Lemma $3.6$ of \cite{rabah2024fukaya}, after replacing $N$ by $X\setminus\mathbf L$ and using $U$ as in Lemma \ref{Existinceof2-form} above.
         \end{proof}   

         Now consider $(\CP^d, \omega_{FS}, J_{std})$ with its Fubini-Study symplectic form and its standard integrable almost complex structure and let $\alpha\in\pi_2(\CP^d, \RP^d)$. Following the same notation as in \cite{rabah2024fukaya}, we recall the following definition.
		\begin{defn}[Space of Framing for Strips]\label{spaceFramings}
			Denote by $\F_{k+2, l}(\alpha):=\{\phi:(\Sigma, \partial\Sigma; \underline{z}, \underline{w})\rightarrow(\CP^d, \RP^d)\}$ satisfying the following conditions:
			\begin{enumerate}
				\item $(\Sigma, \partial\Sigma,\underline{z}, \underline{w})$ is a marked strip, as in Definition \ref{stripsdefn}.
				\item $\phi$ is a smooth map satisfying $\Bar{\partial}_{J_{std}}\phi=0$ on each irreducible component of $\Sigma$, representing the class $\alpha$.
				\item $\operatorname{Aut}((\Sigma, \partial\Sigma, \underline{z}, \underline{w}), \phi)$ is trivial.
				\item $H^1(\Sigma, \partial\Sigma; \phi^*\mathcal{O}(1), \phi^*_{\partial\Sigma}\mathcal{O}(1)_{\mathbb{R}})=\{0\}$.
			\end{enumerate}
		\end{defn}
		
		\begin{prop}[Proposition 3.13 of \cite{rabah2024fukaya}]\label{SpaceofFramingSmooth}
			The space of framings $\F_{k+2, l}(\alpha)$ is a smooth $O(d+1)$-manifold.
		\end{prop}

        Let $\Sigma$ be a bordered genus-zero Riemann surface and $u:(\Sigma, \partial\Sigma)\rightarrow(X, \mathbf{L})$ be a smooth map having boundary condition on $\mathbf{L}$. We denote by $H_u:=u^*H\rightarrow \Sigma$ and using the same notation as in Lemma \ref{LineBundleLemma} and after pulling-back $\tau$ by $u$, we get a totally-real line bundle $H_u^\mathbb R\rightarrow\partial\Sigma$ of $H_u|_{\partial\Sigma}$. After taking a large enough tensor power of $H\rightarrow X$ and abuse of notation, it follows by Lemma $3.10$ of \cite{rabah2024fukaya} (c.f. Section $3.4$ of \cite{katz2006enumerative}), that we may assume, without loss of generality, that $H^1(\Sigma, \partial\Sigma; H_u, H_u^\mathbb R)=\{0\}$.

        Now suppose that $[u]=\beta\in\pi_2(X, \mathbf{L})$, where $\beta$ is a class that can be represented by a pseudo-holomorphic curve, and set $d:=\Omega(\beta)\geq0$. By Lemma $3.14$ of \cite{rabah2024fukaya}, it follows that $\dim_\mathbb R H^0(\Sigma, \partial\Sigma; H_u, H_u^\mathbb R)=d+1>0$. We also recall from Section $3$ of \cite{rabah2024fukaya}, that $H^0(\Sigma, \partial\Sigma; H_u, H_u^\mathbb R)$ has a canonical pairing 
		\[
		\langle.,.\rangle: H^0(\Sigma, \partial\Sigma; H_u, H_u^\mathbb R)\times H^0(\Sigma, \partial\Sigma; H_u, H_u^\mathbb R)\rightarrow\mathbb R,
		\]
		induced from the complex structure on $\Sigma$ (c.f. Definition 3.15 of \cite{rabah2024fukaya}). 

        Let $F=\{f_0,\dots, f_d\}$ be an $\mathbb R$-basis of the real vectorspace $H^0(\Sigma,\partial\Sigma; H_u, H_u^\mathbb R)$ and consider the holomorphic embedding of $\Sigma$ into $\CP^d$ given by
	\[
	\begin{aligned}
		\phi_F: (\Sigma, \partial\Sigma) &\rightarrow (\CP^d, \RP^d) \\
		z &\mapsto [f_0(z):\dots:f_d(z)].
	\end{aligned}
	\]

\begin{defn}[Framed Marked Strips]\label{FramedStrips}
		A framed marked strip, representing the class $\beta\in\pi_2(X, \mathbf L)$, is a tuple $((\Sigma, \partial\Sigma, \underline{z}, \underline{w}), u, F)$ where:
		\begin{enumerate}
			\item $(\Sigma, \partial\Sigma; \underline{z}, \underline{w})$ is a marked strip, as in Definition \ref{stripsdefn}. 
			\item $u:(\Sigma, \partial\Sigma)\rightarrow (X, \mathbf L)$ smooth map, representing the class $\beta$ and $\operatorname{Aut}(u)$ is of finite order, with boundary conditions on $\mathbf L$ given as in item $3$ of Definition \ref{J-stripsDefn}.
 			\item $F$ is an $\R$-basis of $H^0(\Sigma, \partial\Sigma; H_u, H_u^\R)$ such that the matrix $A_F:=(\int_\Sigma\langle f_i, f_j\rangle u^*\Omega)_{i, j}$ is a positive-definite matrix.
		\end{enumerate}
	\end{defn}

\begin{remark}
    It follows from the continuity of $u$ that $u$ maps nodal points between strip-components to intersection points of $L_0\cap L_1$. 
\end{remark}
    
	\begin{prop}[Proposition 3.16 of \cite{rabah2024fukaya}]\label{FramingMap}
		Suppose that $((\Sigma, \partial\Sigma, \underline{z}, \underline{w}), u, F))$ is a framed marked strip representing the class $\beta$ and write $F:=\{f_0,\dots,f_d\}$. Then the induced holomorphic map $\phi_F$ is an element of $\F_{k+2, l}(\alpha)$, where $\alpha\equiv d[\CP^1, \RP^1]$. On the other hand, given $\phi\in\F_{k+2, l}(\alpha)$ then there exists an $\R$-basis of $H^0(\Sigma,\partial\Sigma; \phi^*\mathcal{O}(1), \phi^*_{\partial\Sigma}\mathcal{O}(1)_\R)$, which is  unique up to conjugation by $O(d+1)$, inducing $\phi$. 
	\end{prop}

\subsection{Construction of Global Kuranishi Charts}
We will follow the construction of \cite{abouzaid2024gromov} and \cite{rabah2024fukaya}. In return, we start by recalling the relevant notions:

\begin{defn}[Finite-Dimensional Approximation Scheme]\label{finite dimensional}
    Let $G$ be a compact Lie group and $\pi:V\rightarrow B$ be a smooth $G$-equivariant vectorbundle. A finite dimensional approximation scheme $(V_{m}, \lambda_m)_{m\in\mathbb{N}}$ for $C^\infty_c(V)$, the space of compactly-supported smooth sections of $V$, is a sequence of finite-dimensional $G$-representations $(V_{m})$ together with a sequence of $G$-equivariant linear maps $\lambda_m:V_{m}\rightarrow C^\infty_c(V)$, for all $m\in\mathbb{N}$ satisfying the following conditions:
    \begin{enumerate}
        \item $V_m\leq V_{m+1}$ is a sub-representation, for every $m\in\mathbb{N}$.
        \item $\lambda_{m}|V_{m-1}=\lambda_{m-1}$, for every $m\in\mathbb{N}$.
        \item $\bigcup_{m}\lambda_{m}(V_{m})\subset C^{\infty}(V)$ is dense with respect to the $C^\infty_{loc}$-topology.
    \end{enumerate}
\end{defn}

\begin{lemma}\cite{abouzaid2024gromov}\label{finiteAppScheme}
    Finite dimensional approximation schemes exists.
\end{lemma}

To simplify our notation, we denote by $\F:=\F_{k+2, l}(\alpha)$. As $\F$ is a smooth manifold (c.f. Proposition \ref{SpaceofFramingSmooth}) we consider its universal curve $\mathcal C\rightarrow \F$. Note that, $\mathcal C$ can be identified with $\F_{k+2, l+1}$ and under such identification, $\mathcal C\rightarrow \F$ is the forgetful map. In particular, the fibre of $\mathcal C$ over a point $[(\Sigma, \partial\Sigma, \underline{z}, \underline{w}),\phi]\in \F$ comes with an identification $\iota:(\Sigma, \partial\Sigma)\cong\mathcal C_\phi$. Moreover, as per Proposition \ref{SpaceofFramingSmooth}, we have that $\mathcal C$ is a smooth $O(d+1)$-manifold.

In order to establish a pairing on $H^0(\Sigma, \partial\Sigma; H_u, H_u^\R)$ as given above, that varies continuously with $\phi\in\F$, we recall the following definitions:

    \begin{defn}[Fibre-wise Metric]
    A fibre-wise metric on $\mathcal{C}$ is a continuous map $\mu:\mathcal{C}\times_\F\mathcal{C}\rightarrow\R_{\geq 0}$ such that the restriction of $\mu$ to each irreducible component of the fibre over a point in $\F$, is a metric induced from a smooth Riemannian metric on such irreducible component.
\end{defn}
\begin{defn}[Consistent Domain Metric]\label{consistentMetricDefn}
    A consistent domain metric on $\F$ is an $O(d+1)$-invariant fibre-wise metric on $\mathcal{C}$.
\end{defn}

\begin{lemma}[Lemma 3.34 of \cite{rabah2024fukaya}]\label{metricExist}
    $\F$ can be equipped with a consistent domain metric.
\end{lemma}

Notice that $\mathcal C^*\subset\mathcal C$, the complement of nodes and marked points is a smooth $O(d+1)$-manifold. To simplify our notations further, we set $G:=O(d+1)$ and consider the following definition:
\begin{defn}\label{cokernelVectorspace}
    Let $V:=\bigwedge^{0,1}_{\mathcal C^*/\F}\otimes_\C TX\rightarrow \mathcal C^*\times X$ the $G$-equivariant vectorbundle whose fibre over a point $(c,x)\in\mathcal C^*\times X$ is the space of all complex anti-linear maps $T_c(\mathcal C^*|_{\pi(c)})\rightarrow T_xX$, where $\pi:\mathcal C\rightarrow \F$ is the forgetful map.
\end{defn}
Using Lemma \ref{finiteAppScheme}, let $(V_m,\lambda_m)_{m\in\mathbb N}$ be a finite dimensional approximation scheme of $V$ as given above and fix two points $p,q\in L_0\cap L_1$.

\begin{defn}
    We define the pre-thickened moduli space $\T^{pre}\equiv\T^{pre}(p,q,\beta, (V_m,\lambda_m)_{m\in\mathbb N})$ to be the set of tuples $((\Sigma, \partial\Sigma, \underline{z}, \underline{w}), u, F, e)$ where $((\Sigma, \partial\Sigma, \underline{z}, \underline{w}), u, F)$ is framed marked strip, as in Definition \ref{FramedStrips}, representing the class $\beta$ and $e\in V_m$ for some $m\in\mathbb N$. We require $u$ to satisfy:
    \begin{equation}\label{perturbedDbar}
    \Bar{\partial}_J u|_{\mathcal C^*_{\phi_F}}+\lambda_m(e)\circ\Gamma_u =0, u(z_-)\equiv p\text{ and }u(z_+)\equiv q,
\end{equation}
where $\phi_F:(\Sigma, \partial\Sigma)\rightarrow (\CP^d, \RP^d)$ is the holomorphic map induced from $F$, as given above, and $\Gamma_u:\mathcal C_{\phi_F}\underset{\iota}{\cong}(\Sigma, \partial\Sigma)\rightarrow \mathcal C\times X$ is the graph map given by $z\mapsto(z, u(\iota(z)))$.
\end{defn}

Let $((\Sigma, \partial\Sigma, \underline{z}, \underline{w}), u, F,e)\in\T^{pre}$ and for $p>2$, let
\[
\mathcal X_u:=\{\xi\in W^{1, p}(\Sigma,\partial\Sigma,u^*TX, u|_{\partial\Sigma}^*(TL_0\oplus TL_1)):\xi(z)\equiv 0,\forall z\in\partial\Sigma\text{ where }u(z)\in L_0\cap L_1\}.
\]
The linearization of Equation \ref{perturbedDbar} with respect to $u$, having a fixed domain (equivalently for a fixed $\phi_F$) is a linear operator
\[
D:\mathcal X_u\oplus V_m\rightarrow L^p(\Sigma, \wedge^{0,1}\Sigma\otimes_\C u^*TX),
\]
    of the form
\[
D\equiv D_u\Bar{\partial}_J+\lambda_m(e)\circ\Gamma_u,
\]
where $D_u\Bar{\partial}_J:\mathcal{X}_u\rightarrow L^p(\Sigma, \wedge^{0,1}\Sigma\otimes_\C u^*TX)$ is the linearization of $\Bar{\partial}_J$. As $D_u\Bar{\partial}_J$ is a Fredholm operator (c.f. \cite{mcduff2012j} or Section $58$ of \cite{fukayalagrangian}) and as $\lambda_m\circ \Gamma_u$ is a compact operator, it follows that $D$ is a Fredholm operator. Using the same proof of Lemma $3.38$ of \cite{rabah2024fukaya}, we get the following Lemma:

        \begin{lemma}\label{truncatedFiniteDimApp}
            There exists an $m\in\mathbb N$ such that the operator $D$ as given above, is surjective at every point of $\T^{pre}$.
        \end{lemma}

 \begin{defn}[Topological Thickening]\label{topThichening} We define $\T^{top}$ to be the set of $((\Sigma, \partial\Sigma, \underline{z}, \underline{w}), u, F,e)\in\T^{pre}$ such that $e\in V_m$ where $m\in\mathbb N$ is as given in Lemma \ref{truncatedFiniteDimApp}. We give $\T^{top}$ the Gromov-Hausdorff topology on the graphs induced from the closure of the image of $\Gamma_u$, together with the vectorspace topology on $V_m$. 
        \end{defn}

        \begin{lemma}[Lemma 3.25 of \cite{rabah2024fukaya}]\label{SameTopLemma}
            The topology on $\T^{top}$ agrees with the Gromov's convergence topology.
        \end{lemma}

Note that, we have a canonical projection map $\T^{top}\rightarrow\F$ given by $((\Sigma, \partial\Sigma, \underline{z}, \underline{w}), u, F, e)\mapsto[(\Sigma, \partial\Sigma, \underline{z}, \underline{w}),\phi_F]$. With this in mind and using Lemma \ref{truncatedFiniteDimApp}, we get the following:

        \begin{lemma}[Corollary 3.28 of \cite{rabah2024fukaya}] \label{SmoothingLemma}
            $\T^{top}$ is a topological $O(d+1)$-manifold with a fibre-wise smooth structure on $\T^{top}\rightarrow\F$. Moreover, the stabilizer group at each point $((\Sigma, \partial\Sigma,\underline{z}, \underline{w}),u, F,e)\in \T^{top}$ is a finite group isomorphic to $\operatorname{Aut}((\Sigma, \partial\Sigma, \underline{z}, \underline{w}), u)$, the automorphism group of the underlying pseudo-holomorphic curve.
        \end{lemma}
        \begin{proof}
            We first note that, the $O(d+1)$-action on $\T^{top}$ is given by:
            \begin{enumerate}
                \item the $O(d+1)$-action on $V_m$,
                \item change of basis on $F$ and
                \item precomposition of $u$ with the domain automorphism induced from the $O(d+1)$-action on $\mathcal{C}|_{\phi_F}$ and $\iota$. 
            \end{enumerate}
            With this in mind, the first statement follows from Lemma \ref{truncatedFiniteDimApp} and the fact that surjectivity of $D$, is an open condition. As for the second claim, we argue as in Lemma $6.4$ of \cite{abouzaid2021complex} or Theorem $4.30$ of \cite{abouzaid2024gromov} after adhering to the structure Theorems of \cite{kwon2000structure} or \cite{lazzarini2011relative}.
        \end{proof}

             \begin{prop}[Proposition 3.29 of \cite{rabah2024fukaya}]\label{SmoothingProp}
            There exists a real $O(d+1)$-representation $R$ so that $\T^{top}\times R$ can be given an $O(d+1)$-equivariant smooth structure. Moreover, such smooth structure is unique up to stable $O(d+1)$-isotopy of smooth structures. 
        \end{prop}

        \begin{proof}
            The proof follows the same train of thought as in \cite{abouzaid2021complex}, \cite{bai2022arnold}, \cite{abouzaid2024gromov} and \cite{rabah2024fukaya}. Namely, after showing that $\pi:\T^{top}\rightarrow\F$ is an $O(d+1)$-equivariant $C^1_{loc}$-topological submersion, our claim follows from $\textit{Relative-Lashoff's Theorem}$ as in Theorem $6.9, 6.10$ of \cite{bai2022arnold} and the lifting argument of the vertical tangent microbundle as found in Section $4.1$ of \cite{abouzaid2021complex}. To avoid notational clutter, we will denote by $\mathcal M\equiv\overline{\mathcal M}_{k+1, l+l'}$ the Deligne-Mumford moduli space of genus-zero bordered Riemann surfaces with $(k+1)$-boundary marked points and $(l+l')$-interior points, where $l'\geq\dim_\R\F+3$ a fixed integer. Moreover, for a point $((\Sigma, \partial\Sigma, \underline{z}, \underline{w}), u,F,e)\in\T^{top}$ we will abuse notation and denote its image by $\pi$ by $a$.

            Let $\phi_F\in\F$ representing $a$ and choose $D_1, \dots, D_{l'}$ linear subspaces of $\CP^d$ in general position, intersecting $\phi_F$ transversally. Denote by $\F(D)\subset\F$ the open set of all elements of $\F$ such that their underlying map intersects each of the $D_i$'s transversally. Such construction provides us a local diffeomorphism $\F(D)\rightarrow\overline{\mathcal M}$ given by sending the domain of $\phi\in\F(D)$ together with the intersection points with $D_i$ to $\overline{\mathcal M}$ (c.f. Proposition $6.5$ of \cite{abouzaid2021complex}). Let $\Pi:\pi^{-1}(\F(D))\rightarrow \overline{\mathcal M}$ be the composition of the above map with $\pi$. For a fixed point $q\in\T^{top}$ let $U$ be an open set of $\overline{\mathcal M}$  so that $\F(D)\rightarrow\overline{\mathcal M}$ is a diffeomorphism onto $U$ and let $K\subseteq\Pi^{-1}(U)$ be an open set.

            We claim that there exists a continuous map $G:U\times K\rightarrow\T^{top}$ such that:
            \begin{enumerate}
                \item $G$ is a homeomorphism onto its image,
                \item for any $a\in U$, $G_a(k)\equiv G(a, k)$ is a diffeomorphism onto its image and
                \item $\pi\circ G=pr_U$ where $pr_U:U\times K\rightarrow U$ denotes the canonical projection.
            \end{enumerate}
            Indeed, fix $q\in\T^{top}$ and denote by $a:=\pi(q)$. By Lemma \ref{truncatedFiniteDimApp}, it follows that the linearization $D_q$ at $q$ is surjective. Pick a complement of $\ker(D_q)$. Theorem $3.1$ of \cite{schmaschke2016floer}, allows us to use the $\textit{Newton Iteration}$ as per \cite{fukaya2024exponential} at the domains that do not have mixed boundary conditions, while on such domains we adhere to Lemma $7.6$ and Corollary $7.7$ of \cite{schmaschke2016floer}, in such complement of $\ker (D_q)$. Thus using Theorem $8.16, 8.17$ and $8.25$ of \cite{fukaya2024exponential}, we set $G_a(k)$ to be the pre-glued curve as in Equation $5.4$ of \cite{fukaya2024exponential} and Lemma $7.9$ of \cite{schmaschke2016floer} followed by sufficiently many Newton-iterations.

            As $G_a$ is given by Newton-iterations, it follows that $G_a$ is smooth. Moreover, by Theorem $8.16$ \cite{fukaya2024exponential} and Section $7.6$ of \cite{schmaschke2016floer}, it follows that $G_a$ is surjective, while injectivity follows as we are performing Newton-iteration in the complement of $\ker(D_q)$. In particular, possibly after shrinking $U$ and $K$ to precompact neighborhoods, it follows that $G$ is a continuous bijection into a Hausdorff space and thus a homeomorphism.
            
            Now after applying the parametric Implicit Function Theorem on G, implies that $\pi:\T^{top}\rightarrow\F$ is a $C^1_{loc}$-topological submersion.
        \end{proof}

\begin{defn}[Thickening]\label{ThickeningDefnI}
            We define $\T:=\T^{top}\times R$ with its $O(d+1)$-equivariant smooth structure, where $R$ is as in Proposition \ref{SmoothingProp}.
        \end{defn}

        \begin{defn}
            Fix a consistent domain metric on $\F$, as in Definition \ref{consistentMetricDefn}, on $\F$ and let $((\Sigma, \partial\Sigma, \underline{z}, \underline{w}), u, F, e; v)\in\T$. We define $A_F:=(\int_{\mathcal C|_{\phi_F}}\langle f_i,f_j\rangle u^*\Omega)_{i,j}$, where $\langle.,.\rangle$ is the $O(d+1)$-equivariant metric on $\mathcal C|_{\phi_F}$ induced from Lemma \ref{metricExist}.
        \end{defn}

        \begin{defn}[Obstruction Bundle]\label{ObsBundleDefnI}
            We define the obstruction bundle to be the trivial bundle $\E\rightarrow\T$ of fibre $V_m\oplus\mathcal H_d\oplus R$, where $\mathcal H_d$ is the space of $(d+1)\times(d+1)$-symmetric real matrices.
        \end{defn}

        \begin{defn}[Kuranishi Section]\label{KuranishiSecDefnI}
            We define the Kuranishi section $s:\T\rightarrow\E$ by $((\Sigma, \partial\Sigma, \underline{z}, \underline{w}), u, F, e; v)\mapsto(e, \exp^{-1}(A_F); v)$, where $\exp:\mathcal{H}_d\rightarrow\mathcal H_g^+$ is the matrix exponential map into the convex cone of all symmetric positive-definite $(d+1)\times(d+1)$-matrices.
        \end{defn}

        Using the same arguments as in \cite{rabah2024fukaya}, we can now give a proof of Theorem \ref{GlobalChartThmI}.

        \begin{proof}[Proof of Theorem \ref{GlobalChartThmI}]
            Set $G\equiv O(d+1)$, $\T$ be as in Definition \ref{ThickeningDefnI}, $\E$ be as in Definition \ref{ObsBundleDefnI} and $s$ as in Definition \ref{KuranishiSecDefnI}. We study the zero-locus of $s:\T\rightarrow\E$. Namely, we define $s^{-1}(0)/G\rightarrow \M_{k+2, l}( L_0, L_1, p, q; J, \beta)$ by $[(\Sigma, \partial\Sigma, \underline{z}, \underline{w}), u, I_d, 0; 0]\mapsto[(\Sigma, \partial\Sigma, \underline{z}, \underline{w}), u]$. By Lemma \ref{SameTopLemma}, it follows that such map is a continuous bijection between compact and Hausdorff spaces and hence a homeomorphism. Using the first paragraph of the proof of Lemma \ref{SmoothingLemma}, it follows that such homeomorphism is in fact an equivalence of orbispaces. Therefore, $(G, \T,\E, s)$ serves as a Global Kuranishi chart representation of $\overline{\mathcal M}_{k+2, l}(L_0, L_1, p, q; J, \beta)$. To show that, possibly after further stabilization, we can represent $\overline{\mathcal M}_{k+2, l}(L_0, L_1, p, q; J, \beta)$ by a normally complex Global Kuranishi chart, we argue as in the second paragraph of the proof of Theorem $3.1$ of \cite{rabah2024fukaya}.
        \end{proof}
        
     \subsection{Compatibility with Codimension $1$ Corners}  
     As for now, for a given $p,q\in L_0\cap L_1$ (which may not be distinct), $\beta\in\pi_2(p,q)$ and $k,l\geq0$ integers, Theorem \ref{GlobalChartThmI} provides us with a Global Kuranishi chart $(G,\T,\E,s)$ representing $\M_{k+2,l}(L_0,L_1,p,q;J,\beta)$. On the other hand, by Floer-Gromov compactness (c.f. Section $3.7.4$ of \cite{Fukaya2009}) the codimension-$1$ corner of $\M_{k+2,l}(L_0,L_1,p,q;J,\beta)$ is of the form:
     \begin{align}
         &\bigsqcup_{r\in L_0\cap L_1}\M_{k_1+2,l_1}(L_0,L_1,p,r;J,\beta_1)\times \M_{k_2+2,l_2}(L_0,L_1,r,q;J,\beta_2)\\
         &\sqcup\bigsqcup\M_{k_1+1,l_1}(L_0;J,\beta_1){}_{ev}\times_{ev}\M_{k_2+2,l_2}(L_0,L_1,p,q;J,\beta_2)\\
         &\sqcup\bigsqcup\M_{k_1+2,l_1}(L_0,L_1,p,q;J,\beta_1){}_{ev}\times_{ev}\M_{k_2+1,l_2}(L_1;J,\beta_2),
     \end{align}
     where the union is taken over $k_1+k_2=k$, $l_1+l_2=l$ and $\beta_1+\beta_2=\beta$. More importantly, to be able to use our constructions for geometric applications as in Section \ref{AppSection}, we have to construct our Global Kuranishi charts in a $\textit{compatible}$ way.

     To be more precise, if for $i=1,2$; $k_i,l_i\geq0$ integers and $\beta_i\in\pi_2(X,\mathbf{L})$ such that:
     \begin{enumerate}
         \item $\beta_1+\beta_2=\beta$.
         \item $k_1+k_2=k$.
         \item $l_1+l_2=l$.
     \end{enumerate}
     In addition, $(G_i,\T_i,\E_i,s_i)$ is a Global Kuranishi chart representation for some moduli space appearing as a factor of the above stratification and if we denote by $\F_i$ the corresponding framing space as in Definition \ref{spaceFramings} or Definition $3.11$ of \cite{rabah2024fukaya}, we require:

\begin{enumerate} 
    \item \label{GroupEmbedding }$G_1\times G_2\hookrightarrow G$.
    \item \[\begin{tikzcd}
	{\mathcal{T}_1{}_{ev}\times_{ev}\mathcal{T}_2} && {\mathcal{T}} \\
	\\
	{\mathcal{F}_{1}{}_{ev}\times_{ev}\mathcal{F}_2} && {\mathcal{F}}
	\arrow[hook, from=1-1, to=1-3]
	\arrow[from=1-1, to=3-1]
	\arrow[from=1-3, to=3-3]
	\arrow[hook, from=3-1, to=3-3]
\end{tikzcd}\]
to be a commutative diagram of smooth manifolds.
\item \[\begin{tikzcd}
	{\iota^*(\mathcal{E}_1\oplus\mathcal{E}_2)} && {\mathcal{E}} \\
	\\
	{\mathcal{T}_1{}_{ev}\times_{ev}\mathcal{T}_2} && {\mathcal{T}}
	\arrow[hook, from=1-1, to=1-3]
	\arrow[from=1-1, to=3-1]
	\arrow[from=1-3, to=3-3]
	\arrow["{\iota^*(s_1\oplus s_2)}", bend left = 30pt, from=3-1, to=1-1]
	\arrow[hook, from=3-1, to=3-3]
	\arrow["s"', bend right = 30pt, from=3-3, to=1-3]
\end{tikzcd}\]
to be a commutative diagram, where $\iota:\T_1\times_{ev}\T_2\hookrightarrow\T_1\times\T_2$ is the canonical inclusion.
\item The fibre-product Global Kuranishi, $(G_1\times G_2, \T_1\times_{ev}\T_2, \E_1\oplus\E_2,s_1\oplus s_2)$ is equivalent to\\ $(G, \T|_{\pi^{-1}(\iota(\F_1\times_{ev}\T_2))}, \iota^*\E,\iota^*s)$, the restriction of $(G, \T, \E, s)$ to the boundary corresponding to degeneration of the form $\F_1\times_{ev}\F_2$, where $\pi:\T\rightarrow \F$ is the forgetful map.
\item $\iota^*ev=(ev_1, ev_2)$, where $ev$ is the evaluation map on $\T$ and $ev_i$ is the evaluation map on $\T_i$ for $i=1,2$.
\item If $(G_i,\T_i,\E_i,s_i)$ is a Global Kuranishi representation of either $\M_{k_i+1,l_i}(L_j;J,\beta_i)$ where $j=0,1$ then, it is the same chart used in the proof of Theorem $5.3$ of \cite{rabah2024fukaya}.
\end{enumerate}

To be consistent with \cite{rabah2024fukaya}, \cite{rabah2025fukaya}, we will refer to the above conditions as the $\textbf{compatibilty conditions}$.

Following \cite{rabah2024fukaya}, \cite{rabah2025fukaya}, we give the following definition:

\begin{defn}[Framed Marked Strip with Constraint]\label{framedStripsConstraint}
    A framed marked strip with constraint, representing the class $\beta\in\pi_2(X,\mathbf{L})$ is a framed marked strip $((\Sigma,\partial\Sigma,\underline{z},\underline{w}), u, F)$ as in Definition \ref{FramedStrips} such that, if we write $F=\{f_0,\dots,f_d\}$ then we require $[f_i(z_0)]_{i=0,\dots,d}=[1:0:\dots:0]$.
\end{defn}

Using the same proof of Lemma $3.48$ in \cite{rabah2024fukaya}, it follows that $\F':=\{\phi\in\F:\phi(z_0)=[1:0:\dots:0]\}$ where $\F$ is as in Definition \ref{spaceFramings}, is a smooth $O(d)$-manifold where $d:=\Omega(\beta)$ and $\Omega$ is as in Lemma \ref{Existinceof2-form}.

\begin{lemma}\label{refinedChart}
    The moduli space $\M_{k+2,l}(L_0,L_1,p,q;J,\beta)$ can be given a Global Kuranishi chart representation of the form $(O(d), \T,\E,s)$ together with a smooth sumbermsion $ev_i:\T\rightarrow\mathbf{L}$ extending the evaluation map at the $i^\text{th}$-boundary marked point, where we view $\mathbf{L}$ as an immersed submanifold of $X$.
\end{lemma}

\begin{proof}
    We follow the same argument of the proof of Corollary $3.49$ in \cite{rabah2024fukaya}. To be more precise, as $\F'$ is a smooth $O(d)$-manifold, the universal curve $\mathcal C'\rightarrow\F'$ over it also retains a smooth $O(d)$-manifold structure making the projection map $O(d)$-equivariant. Let $V\rightarrow \mathcal C'\times X$ be as in Definition \ref{cokernelVectorspace}. Using Lemma \ref{truncatedFiniteDimApp}, it follows that we can find $V_m$ a finite-dimensional $O(d)$-representation such that  $\T':=\{((\Sigma,\partial\Sigma,\underline{z},\underline{w}),u, F, e)\}$ such that $((\Sigma,\partial\Sigma,\underline{z},\underline{w}),u, F)$ is as in Definition \ref{framedStripsConstraint} and $e\in V_m$, is a topological $O(d)$-manifold together with a fibre-wise smooth structure on the projection map $\T'\rightarrow\F'$. Now after stabilizing $ev_i^*TL_0\oplus ev_i^* TL_1$ and using the same arguments as in the proof of Lemma $4.5$ of \cite{abouzaid2021complex} and the arguements of Proposition \ref{SmoothingProp}, it follows from $\textit{Lashoff's Smoothing Theorem}$ as in \cite{Lashof1979} that there exists a finite-dimensional real $O(d)$-representation such that $\T:=\T'\times R$ is a smooth $O(d)$-manifold. To finish the proof, we set $\E\rightarrow\T$ to be the trivial vectorbundle of fibre $V_m\oplus\mathcal H_d\oplus R$ and $s$ is given analogously to Definition \ref{KuranishiSecDefnI}, and argue as in the proof of Theorem \ref{GlobalChartThmI} above.
\end{proof}

The above Lemma settles the first requirement in the compatibility conditions listed above. Nevertheless, we still need modify our Global Kuranishi chart representations to get items $(2)$ and $(3)$ in the compatibility conditions. To do so, we will modify our construction inductively using the Lexographic partial order $\precsim$ on the monoid $\mathfrak{G}:=\{(\Omega(\beta),k,l)\}$. For $i=1,2$, we denote by:
\begin{enumerate}
    \item $d_i:=\Omega(\beta_i)$.
    \item $G_i=O(d_i)$ and $G=O(d)$.
    \item $\F^\text{strip}_i$ is the space of marked framed strips with constraints as in Definition \ref{framedStripsConstraint}, and $\F_i^\text{disk}$ is the space framed marked disks with constraints as in Definition $3.47$ of \cite{rabah2024fukaya}.
    \item $\mathcal C^\text{disk}_i\rightarrow\F^\text{disk}_i$ and $\mathcal{C}_i^\text{strip}\rightarrow\F_i^\text{disk}$ the corresponding universal curves.
\end{enumerate}

We start by giving an explicit definition of the maps $\F_1{}_{ev}\times_{ev}\F_2\hookrightarrow\F^\text{strip}$. It is important to note that for $\alpha\in\pi_2(\CP^d,\RP^d)$ and $k,l\geq0$ integers, $\F^\text{strip}_{k+2,l}(\alpha)$ is nothing more than $\F^\text{disk}_{k+2,l}$. Nevertheless, in order to prove item $(2)$ of the above compatibility conditions, we must keep track of the marked points we are concatenating that framings at. With this said and having the Equations given in Section $3.5$ of \cite{rabah2024fukaya} in mind, we will only give the Equations defining $\F_1^\text{strip}{}_{ev}\times_{ev}\F_2^\text{strip}\hookrightarrow\F^\text{strip}$, as all the other ones are a matter of notational change. Suppose that $\phi_i\in\F_i^\text{strip}$ and note that, as in the proof of Proposition $3.16$ of \cite{rabah2024fukaya}, we can write $\phi_i(z)=[f_j^i(z)]_{j=0,\dots,d_i}$. Suppose that $\Sigma=\Sigma_1\vee\Sigma_2:=\Sigma_1\cup_{z^1_+\sim z_-^2}\Sigma_2$ corresponding to a nodal point $u_1(z_+^1)=u_2^2(z^2_-)=r\in L_0\cap L_1$. Let $a_j:=f^1_j(z^1_+)$, $b:=f^2_0(z_-^2)>0$ and $c:=(\sum_{0\leq j\leq d_1}a_j^2)^\frac{1}{2}$. For $j=0,\dots,d_1$, set 
\begin{equation*}
\tilde f^1_j(z):=
    \begin{cases}
        &f^1_j(z), \text{ for }z\in\Sigma_1\\
        &\frac{a_j}{b}f_0^2(z), \text{ for }z\in\Sigma_2.
    \end{cases}
\end{equation*}
For $j=0,\dots, d_2$, set
\begin{equation*}
\tilde f^2_j(z):=
    \begin{cases}
        &0,\text{ for }z\in\Sigma_1\\
        &\frac{c}{b}f^2_j(z),\text{ for }z\in\Sigma_2.
    \end{cases}
\end{equation*}
Finally, we define $\phi_1\vee\phi_2:\Sigma\rightarrow\CP^{d_1+d_2}$ by
\[
\phi_1\vee\phi_2(z):=[\tilde f^1_0(z):\dots:\tilde f^1_{d_1}(z):\tilde f^2_1(z):\dots:\tilde f^2_{d_2}(z)].
\]
Such construction gives the desired embedding $\F_1{}_{ev}\times_{ev}\F_2\hookrightarrow\F$ making it compatible with the diagonal embedding $G_1\times G_2\hookrightarrow G$.

\begin{lemma}\label{compatibleObsBundle}
    There exists a finite-dimensional approximation scheme, as in Lemma \ref{truncatedFiniteDimApp}, for the moduli \\space $\M_{k+2,l}(L_0,L_1,p,q;J,\beta)$ such that if
    \begin{enumerate}
        \item $V_{1}^\text{strip}$ is a finite-dimensional approximation scheme for $\M_{k_1+2,l_1}(L_0,L_1,p,r;J,\beta_1)$,
        \item $V_2^\text{strip}$ is a finite-dimensional approximation scheme for $\M_{k_2+2,l_2}(L_0,L_1,r,q;J,\beta_2)$ and
        \item $V^\text{disk}_i$ is a finite-dimensional approximation scheme for $\M_{k_i+1,l_i}(X,L_{i-1};J,\beta_i)$ as given in Lemma $3.50$ of \cite{rabah2024fukaya}
    \end{enumerate}
       then,
    \begin{enumerate}
        \item $V_1^\text{strip}\oplus V_2^\text{strip}$,
        \item $V_1^\text{disk}\oplus V^\text{strip}_2$ and 
        \item $V^\text{strip}_1\oplus V_2^\text{disk}$
    \end{enumerate}
    is a $G$-subrepresentation of $V$.
\end{lemma}

\begin{proof}
    We argue by induction over the discrete monoid $\mathfrak{G}$ as given above. The base case follows from Lemma \ref{truncatedFiniteDimApp}. As for the inductive step, we give $\M_{k+2,l}(L_0,L_1,p,q;J,\beta)$ a finite-dimensional approximation scheme, say $V'$, using Lemma \ref{truncatedFiniteDimApp}. Then, we consider one of the possible degenerations of $\F$. That is, denote by $\iota:\F_1{}_{ev}\times_{ev}\F_2\hookrightarrow\F$ and by $\mathcal C\rightarrow\F$ the universal curve. Then, by the definition of differential forms on nodal Riemann surfaces, we will have the following commutative diagram:
    \[\begin{tikzcd}
	{\bigwedge^{0,1}_{\iota^*\mathcal{C}/\mathcal{F}_1{}_{ev}\times_{ev}\mathcal{F}_2}\otimes_{\mathbb{C}}TX} && {\bigwedge^{0,1}_{\mathcal{C}/\mathcal{F}}\otimes_{\mathbb{C}}TX} \\
	\\
	{\iota^*\mathcal{C}\times X} && {\mathcal{C}\times X}
	\arrow[hook, from=1-1, to=1-3]
	\arrow[from=1-1, to=3-1]
	\arrow[from=1-3, to=3-3]
	\arrow[hook, from=3-1, to=3-3]
\end{tikzcd}\]
Now the argument of the proof of Lemma $3.50$ of \cite{rabah2024fukaya} (c.f. Lemma $2.6$ of \cite{rezchikov2022integral}) implies that we can enlarge $V'$ to $V$, that is there exists a finite-dimensional $G$-representation $H$ such that $V:=H\oplus V'$ is a finite-dimensional approximate for $\M_{k+2,l}(L_0,L_1,p,q;J,\beta)$ having one of the items in the results of the above statement, as a $G$-subrepresentation. Now repeating this argument finitely many times, corresponding to all possible ways we have $\beta_1+\beta_2=\beta$, $k_1+k_2=k$ and $l_1+l_2=l$, we get the desired result. 
\end{proof}

\begin{theorem}\label{GKCStrips}
    There exists a Global Kuranishi chart $(G,\T,\E,s)$ representing the moduli space $\M_{k+2,l}(L_0,L_1,p,q;J,\beta)$ such for $i=1,2$ and any $\beta_i,k_i,l_i$ as above, each factor of the codimension-$1$ corner of $\M_{k+2,l}(L_0,L_1,p,q;J,\beta)$ as given above, has a Global Kuranishi chart representation $(G_i,\T_i,\E_i,s_i)$ satisfying all the compatibility conditions.
\end{theorem}

\begin{proof}
        We argue by induction over $\mathfrak{G}$ using Theorem \ref{GlobalChartThmI} and Theorem $3.51$ of \cite{rabah2024fukaya}. The base case follows from Theorem \ref{GlobalChartThmI}. As for the inductive step we only show how to extend our construction in the case of $\textit{disk bubbling}$ as the other case is easier, after noting that non-empty components of a fibre-product over a finite set is just direct product, and to ensure that our construction is compatible with that of \cite{rabah2024fukaya}, \cite{rabah2025fukaya}.

        Indeed, suppose that $(G_1,\T_1,\E_1,s_1)$ is a Global Kuranishi chart representation for $\M_{k_1+1,l_1}(X,L;J,\beta_1)$ as given in Theorem $3.51$ and $(G_2,\T_2,\E_2,s_2)$ is a given Global Kuranishi chart representation of\\$\M_{k_2+2,l_2}(L_0,L_1,p,q;J,\beta_2)$. By construction $\E_i$ is a trivial vectorbundle over $\T_i$ of fibre $V_i\oplus\mathcal H_{d_i}\oplus R_i$ where $V_i$ is a finite-dimensional approximation and $R_i$ is given by $\textit{Lashoff's Theorem}$ as per \cite{Lashof1979} or Theorems $6.9,6.10$ of \cite{bai2022arnold}. Let $\T':=\{((\Sigma,\partial\Sigma,\underline{z},\underline{w}),u,F,e, (r_1,r_2))\}$ such that $((\Sigma,\partial\Sigma,\underline{z},\underline{w}),u,e,F)$ where $((\Sigma,\partial\Sigma,\underline{z},\underline{w}),u,F,e)$ is as in Definition \ref{topThichening} with the understanding that $e\in V$ such that $V$ is given by Lemma \ref{compatibleObsBundle}, having $V_1\oplus V_2\leq V$ as a $G$-subrepresentations. If such $\textit{disk bubbling}$ is covering $\F_1{}_{ev}\times_{ev}\hookrightarrow\F$, it follows that we have the following commutative diagram:
        \[\begin{tikzcd}
	{\mathcal T_1{}_{ev}\times_{ev}\mathcal T_2} && {\mathcal T'} \\
	\\
	{\mathcal F_1{}_{ev}\times_{ev}\mathcal F_2} && {\mathcal F}
	\arrow[hook, from=1-1, to=1-3]
	\arrow[from=1-1, to=3-1]
	\arrow[from=1-3, to=3-3]
	\arrow[hook, from=3-1, to=3-3]
\end{tikzcd}\]
We note that, eventhough $\T_1,\T_2$ are smooth $G_1,G_2$-manifolds respectively, yet the upper horizontal arrow is a topological embedding. On the other hand, by Proposition \ref{SmoothingProp}, it follows that $\T'\rightarrow\F$ has a fibre-wise smooth structure. Now we use the same $\textit{trick}$ as in \cite{rabah2024fukaya}, \cite{rabah2025fukaya}. Namely, we use an auxiliary norm on $V_1, V_2$ and the fact that transversality is an open condition, to find a $G$-invariant open set in $\T'$ containing $\T_1{}_{ev}\times_{ev}\T_2$ (c.f. the last paragraph of the proof of Proposition $3.27$ of \cite{rabah2024fukaya}) to apply Theorem $6.9, 6.10$ of \cite{bai2022arnold}, which enable us to find a finite-dimensional real $G$-representation $R$ such that $\T:=\T'\times R$ admits a smooth $G$-manifold structure. Now the second item of the $\textbf{compatibility conditions}$ above, follows after viewing $R$ as $G_2$-representation and stabilizing $(G_2,\T_2,\E_2,s_2)$ by the trivial vectorbundle $R\times\T_2\rightarrow\T_2$.
Now by Lemma $3.52$ of \cite{rabah2024fukaya} and using the same arguments of the proof of Theorem $3.51$ of \cite{rabah2024fukaya}, the result follows.
\end{proof}

\begin{remark}
    We state the $\textit{trick}$ we used in \cite{rabah2024fukaya}, \cite{rabah2025fukaya} that we mentioned above more precisely. Since in all cases, we are constructing thickening for the moduli spaces of fixed degree, it follows that after fixing an auxiliary norm on the finite-dimensional approximation space, solutions of the Gromov's graph trick as in Equation \ref{perturbedDbar} will have a uniform $L^2$-norm. Then, using the argument as in the last paragraph of the proof of Proposition $3.27$ of \cite{rabah2024fukaya} (c.f \cite{rabah2025fukaya}) we can find an invariant open neighborhood of $\T_1{}_{ev}\times_{ev}\T_2\hookrightarrow\T$ before apply Relative Lashoff's Theorem as in \cite{bai2022arnold}. Such system of neighborhoods serve as $\textit{collars}$ as per \cite{bai2022arnold} or \cite{BaiXu2026Foundations} yet we do not adhere to any analytical results from Chapter $17$ of \cite{FOOO2020}.
\end{remark}

        \section{Open-Closed Map}

Following \cite{rabah2024fukaya}, let $(f, g)$ be a Morse-Smale pair on $L$ and consider the induced $A_\infty$-algebra $(\mathrm{CM}(f; \Lambda_0), \{\m_k\}_{k\geq0})$ on the Morse complex on $L$ given by Theorem $1.2$ of \cite{rabah2024fukaya}. The goal of this section is construct a filtered chain map from the chains of the Hochschild homology of $(\mathrm{CM}(f; \Lambda_0), \{\m_k\}_{k\geq0})$ to some Morse chains on $X$. To be more precise, let $g_X$ be a metric on $X$ and $F:X\rightarrow \R$ be a function such that the pair $(F,g)$ is a Morse-Smale pair. Moreover, we impose the $\textit{generic}$ condition on $(F, g)$ and require that the stable/unstable manifolds of $F$ are transverse to the unstable/stable manifolds of $f$ (c.f. \cite{venugopalan2020fukaya}). To simplify our notations, denote by $V:=\mathrm{CM}^{*+1}(f; \Lambda_0)$ and hence, the deformed $A_\infty$-operations defines a linear map $\m_k^b:V^{\otimes k}\rightarrow V$. We consider the Hochschild complex $CC^*(V):=\widehat{\bigoplus}_{K\geq1}V^{\otimes K}$ where $\delta_H:CC^*(V)\rightarrow CC^*(V)$ is defined on pure tensors by

    \begin{align*}
        \delta_H(v_0,\dots, v_K):=&\sum_{i=1}^K\sum_{j=1}^{K-j+1}(-1)^{\epsilon_i}v_0\otimes\dots\otimes v_{i-1}\otimes \m_{j}(v_i,\dots, v_{i+j-1})\otimes v_{i+j}\otimes\dots\otimes v_K\\
        &+\sum_{l,s\geq0, l+s\leq K}(-1)^{\epsilon_l}\m_{l+s+1}(v_{K-l+1},\dots, v_K, v_0,\dots, v_s)\otimes v_{s+1}\otimes\dots\otimes v_{K-l}
    \end{align*}

where $\epsilon_i:=\sum_{r=0}^{i-1}|v_r|$ and $\epsilon_l:=(\sum_{r=0}^{K-l}|v_r|)(\sum_{r=K-l+1}^K|v_r|)$.

\begin{theorem}\label{OCmapTheorem}
    There exists a degree-zero, filtered map $\p:CC^*(V)\rightarrow \mathrm{CM}^*(F; \Lambda_0)$ such that $\p\circ\delta_H=d^{\text{Morse}}_X\circ\p$.
\end{theorem}

\begin{remark}
    In fact, we will use a variant of Theorem \ref{OCmapTheorem} in \cite{RabahIntegralFukaya} to give a $\textit{weak Calabi-Yau}$  structure as per \cite{Sheridan2016} on the $\textit{integral Fukaya category}$ as per \cite{RabahIntegralFukaya}.
\end{remark}

\subsection{Geometric Preliminaries}
        Recall and following the notation of \cite{Fukaya2009} and \cite{rabah2024fukaya}, the $A_\infty$-algebra we have is $\mathfrak{G}$-gapped where $\mathfrak{G}:=\{(\omega(\beta), \mu(\beta)):\beta\in\Pi(L), \M(\beta)\neq\emptyset\}$. Moreover, using \cite{Fukaya2009} or Section $5$ of \cite{rabah2024fukaya}, we can give a total preorder $\precsim$ on $\Pi(L)\times \N$. The goal of this section is to prove the following:

        \begin{theorem}\label{TruncatedOCmap}
            For every $r\geq 1$ integer there exists a filtered, degree-zero $\Lambda_0$-linear map $\p_r:V^{\otimes r}\rightarrow \mathrm{CM}^*(F; \Lambda_0)$ of the form $\p_r\equiv\sum_{\beta\in\Pi(L)}T^{\omega(\beta)}e^{\frac{\mu(\beta)}{2}}\p_{r, \beta}$ such that for every $(\beta, j)\precsim(n, K)$ and $x_1,\dots, x_k\in\mathrm{Crit}(f)$ we have:
            \begin{align*}
                d^\text{Morse}_X\circ\p_{k, \beta}(x_1,\dots,x_k)&+\sum_{\beta_1+\beta_2=\beta}\sum_{0\leq j\leq k}\sum_{1\leq i\leq k-j}(-1)^{\epsilon_i}\p_{k-j, \beta_2}(x_1,\dots,x_{i-1},\m_{j, \beta_1}(x_i,\dots,x_{i+j-1}),x_{i+j},\dots, x_k)\\
                &+\sum_{\beta_1+\beta_2=\beta}\sum_{l,s\geq0, l+s\leq k}(-1)^{\epsilon_l}\p_{k-l+s, \beta_2}(\m_{l+s,\beta_1}(x_{k-l+1},\dots,x_k,x_1,\dots,x_s),x_{s+1},\dots,x_{k-l})\\
                &=0.
            \end{align*}
        \end{theorem}

We will follow closely the notation and constructions of Section $4$ of \cite{rabah2024fukaya}. For instance, we will denote by $\mathrm{Gr}_{k+1}$ the moduli space moduli space of isomorphism classes of metirc, rooted, ribbon trees with $(k+1)$-external vertices. As per Theorem $14.1$ of \cite{fukaya1997zero}, $\mathrm{Gr}_{k+1}$ admits the structure of a smooth manifold with corners of $\dim_\R\mathrm{Gr}_{k+1}=k-2$ and a cellular stratification $\mathrm{Gr}_{k+1}=\bigcup_{t\in G_{k+1}}\mathrm{Gr}_{k+1}(t)$, where $G_{k+1}$ is the set of all isomorphism classes of rooted, ribbon trees and $\mathrm{Gr}_{k+1}(t)$ denotes the open cell of all metric rooted ribbon trees having the same combinatorial type as $t$, that is the same number of internal vertices and internal edges as $t$. Moreover, we will denote by $\widetilde{\mathrm{Gr}}_{k+1}\rightarrow \mathrm{Gr}_{k+1}$ the universal curve over $\mathrm{Gr}_{k+1}$. In what follows, it is important to note that a metric rooted ribbon tree $(t, l)\in\mathrm{Gr}_{k+1}$ comes with a special point $rt$ the root and we will denote the edge adjacent to root by $e_0$.

\begin{defn}[Open-Closed Fixed Combinatorial Type Moduli Space]
    Let $x_1,\dots,x_k\in\mathrm{Crit}(f)$, $y\in\mathrm{Crit}(F)$ and $\beta_v\in\pi_2(X, L)$ for each internal vertex $v$ of $t$. Fix a rooted, ribbon tree $t$ and domain-dependent Morse function $\Tilde{f}\in C^\infty(\widetilde{\mathrm{Gr}}_{k+1}\times L)$ corresponding to $f$, as in \cite{charest2015floer}, \cite{venugopalan2020fukaya} or Definition $4.19$ of \cite{rabah2024fukaya}. We denote by $\mathcal{M}(t, \underline{x},y;(\beta_v ),\Tilde{f}_t):=\{(l,u): l\in\mathrm{Gr}_{k+1}(t), u:t\rightarrow X\}$ where $u$ is a continuous maps such that:
    \begin{enumerate}
        \item For each internal vertex $v$ we associate a $J$-holomorphic map $u:(\Sigma_v,\partial\Sigma_v)\rightarrow(X, L)$ with boundary condition on $L$ representing the class $\beta_v$, where $(\Sigma_v,\partial\Sigma_v)$ boundary marked corresponding to the valency of $v$.
        \item For each internal edge $e$, we associate a Morse trajectory $u_e:[0, l(e)]_s\rightarrow L$ such that
        \begin{equation}\label{MorseEqu}
            \frac{d}{ds}u_e(s)+\nabla^g\Tilde{f}_t(s, u_e(s))=0,
        \end{equation}
        where $\Tilde{f}_t:\widetilde{\mathrm{Gr}}_{k+1}(t)\rightarrow C^\infty(L)$ is given by $\title{f}_t\equiv\title{f}$.
        \item For each external edge $e$ other than $e_0$, we associate a Morse trajectory $u_e:(-\infty,0]_s\rightarrow L$ solving Equation \ref{MorseEqu} and $\lim_{s\rightarrow-\infty}u_e(s)=x_i$.
        \item For the edge $e_0$ adjacent to the root of $t$ we associate a Morse trajectory $u_{e_0}:[0, \infty)_s\rightarrow X$ such that
        \begin{equation}
            \frac{d}{ds}u_{e_0}(s)+\nabla^{g_X}F(u_{e_0}(s))=0,
        \end{equation}
        where $\lim_{s\rightarrow\infty}u_{e_0}(s)=y$.
        \item If $v_0$ is the internal vertex connected to the root by $e_0$ then, we associate for $v_0$ a $J$-holomorphic map $u_{v_0}:(\Sigma_{v_0},\partial\Sigma_{v_0})\rightarrow(X, L)$ representing the class $\beta_{v_0}$, where $\Sigma_{v_0}$ has an interior marked point denoted by $w\in\Sigma_{v_0}$. We impose the condition that, $u_{v_0}(w)=u_{e_0}(0)$. 
    \end{enumerate}
    We give $\mathcal M(t,\underline{x}, y;(\beta_v))\equiv\mathcal M(t,\underline{x},y;(\beta_v),\Tilde{f}_t)$ the $C^\infty_{loc}$-topology.
\end{defn}

From now on and to avoid repetition, we will denote by $x_1,\dots, x_k\in\mathrm{Crit}(f)$ and $y\in\mathrm{Crit}(F)$.

\begin{lemma}\label{OC ModuliLemmaI}
    Let $t\in G_{k+1}$ and $(\beta_v)$ be a sequence of relative-spherical classes indexed by the internal vertices of $t$. For a generic choice of domain-dependent Morse function $\Tilde{f}\in C^\infty(\widetilde{\mathrm{Gr}}_{k+1}\times L)$ associated to $f$, the moduli space $\mathcal{M}(t,\underline{x},y;(\beta_v),\Tilde{f}_t)$ admits a Global Kuranishi chart representation.
\end{lemma}

\begin{proof}
    This follows from the proof of Proposition $4.37$ of \cite{rabah2024fukaya} (c.f. \cite{rabah2025fukaya}) after noting that in the case of $\beta_{v_0}=0$ then, transversality is achieved from the fact that the stable/unstable manifolds of $(F,g_X)$ are transverse to the unstable/stable manifolds of $(f, g)$.
\end{proof}

\begin{defn}
    For $t\in G_{k+1}$, $\beta\in\pi_2(X,L)$ and a domain dependent Morse function $f\in C^\infty(\widetilde{\mathrm{Gr}}_{k+1}\times L)$ associated to $f$, we denote by $\mathcal{M}(\underline{x}, y; \beta, \Tilde{f}_t):=\bigsqcup_{\sum\beta_v=\beta}\mathcal{M}(t, \underline{x},y;(\beta_v ),\Tilde{f}_t)$.
\end{defn}

\begin{cor}\label{OCmoduliCor}
    For a generic choice of domain-dependent Morse function $\title{f}\in C^\infty(\widetilde{\mathrm{Gr}}_{k+1}\times L)$, $\mathcal{M}(t,\underline{x},y;\beta)$ admits a Global Kuranishi chart representation.
\end{cor}

\begin{proof}
    Noting that, by Gromov's compactness the set of $\beta_i\in\pi_2(X,L)$ such that $\omega(\beta_i)\leq\omega(\beta)$ is finite, the result follows from Proposition \ref{OC ModuliLemmaI}.
\end{proof}

Following \cite{rabah2024fukaya} and \cite{rabah2025fukaya}, we introduce the following Definition:

\begin{defn}
    For $\delta>0$, $t\in G_{k+1}$, $\beta\in\pi_2(X,L)$ and a domain-dependent Morse function $\Tilde{f}\in C^\infty(\widetilde{\mathrm{Gr}}_{k+1}\times L)$ associated to $f$, we denote by $\mathcal{M}^\delta(t,\underline{x},y;\beta,\Tilde{f}_t)\subset\mathcal{M}(t, \underline{x},y;\beta)$ such that, there exists an internal edge $e$ of $t$ with $l(e)<\delta$.
\end{defn}

We recall some of the notations from \cite{rabah2024fukaya} and \cite{rabah2025fukaya}. For $t,t'\in G_{k+1}$, we write $t>t'$ if $t'$ is formed from $t$ after shrinking a single internal edge of $t$. Moreover, we note that, by Theorem $14.1$ of \cite{fukaya1997zero} (c.f. \cite{charest2015floer} or \cite{venugopalan2020fukaya}) $\mathrm{Gr}_{k+1}$ admits the structure of a smooth manifold with corners with a smooth stratification given by $\mathrm{Gr}_{k+1}=\bigcup_{t\in G_{k+1}}\mathrm{Gr}_{k+1}(t)$, where $\mathrm{Gr}_{k+1}(t)$ denotes the open cell containing $t$. For $\epsilon>0$ and after fixing a metric on $\mathrm{Gr}_{k+1}$, we consider $\mathbb D_\epsilon(t')$ an $\epsilon$-neighborhood of $\mathrm{Gr}_{k+1}(t')$ inside $\mathrm{Gr}_{k+1}(t)$ given by the normal bundle of $\mathrm{Gr}_{k+1}(t')\subset\mathrm{Gr}_{k+1}(t)$, with the understanding that $t>t'$. In fact, the fibre of $\mathbb D_\epsilon(t')$ is diffeomorphic to $[0,\epsilon)$.

Now let $\Tilde{f}$ be a generic domain-dependent Morse function, as in Lemma \ref{OC ModuliLemmaI}, and denote the Thickening of $\mathcal{M}(t',\underline{x}, y; \beta,\Tilde{f}_{t'})$ given by Lemma \ref{OC ModuliLemmaI}, by $\T(t',\underline{x}, y;\beta, \Tilde{f}_{t'})$.

\begin{prop}\label{OCmoduliProp}
    There exists $\epsilon,\delta>0$ real numbers where $\delta$ depends on $\epsilon$ and the injectivity radius of the metric on $\mathrm{Gr}_{k+1}$ such that $\mathcal M^\delta(t,\underline{x},y;\beta, \Tilde{f}_t)$ admits a Thickening $\T^\delta(t,\underline{x}, y;\beta,\tilde f_t)$ with the understanding that:
    \begin{enumerate}
        \item $\tilde f_{t}|_{\widetilde{\mathrm{Gr}}_{k+1}(t')}\equiv\tilde f_{t'}$.
        \item $\tilde f_t(s, .)\equiv f$ on the shrinked edge of length $l(e)<\delta$.
    \end{enumerate}
    Fuirthermore, we have a diffeomorphism $\T^\delta(t,\underline{x}, y;\beta,\tilde f_t)\cong\T(t',\underline{x}, y;\beta, \tilde f_{t'})\times \mathbb D_\epsilon(t')$.
\end{prop}

\begin{proof}
    This follows from the proof of Proposition $4.42$ of \cite{rabah2024fukaya}.
\end{proof}

\begin{remark}\label{CompactificationRemark}
    It is important to note that the proof of Proposition $4.42$ of \cite{rabah2024fukaya} works in our case as we are only allowing a single edge to associate to a Morse trajectory of $(F, g_X)$. In other words, we will not allow the $\textit{bigger}$ compactification of $\mathcal M(t,\underline{x},y;\beta)$ to admit configurations where some internal edge will correspond to a Morse trajectory of $(F, g_X)$, equivalently, where $v_0$ corresponds to sphere components. In fact, if we allow such configurations, we will not be able to prove a result of the form as Theorem \ref{OCmapTheorem}. For instance, consider a configuration corresponding to a tree with $(k+1)$-external edges, two internal vertices where the non-root vertex is adjacent to $k$-external edges. Associate to the non-root internal vertex a single disk components of degree $\beta$ with $k$-inputs corresponding to Morse trajectories associated to $x_1,\dots, x_k$. Moreover, suppose that such disk component has an interior marked connecting to a finite length Morse trajectory of $(F,g)$ which in return is connected to a $J$-sphere of degree $\alpha$, and such $J$-sphere connected to a Morse trajectory corresponding to $y\in\mathrm{Crit}(F)$. Taking the $C^\infty_{loc}$-limit of such configurations will result in terms of the form $\mathrm{GW}_{0,2}^\alpha(\p_{k,\beta}(x_1,\dots,x_k))$ in the expansion of $d^\text{Morse}_X\circ\p_{k,\alpha+\beta}(x_1,\dots,x_k)$. Indeed, a similar consideration has been considered in \cite{venugopalan2020fukaya}.
\end{remark}

\begin{defn}[Open-Closed Moduli Space]\label{OCmoduliDefn}
    For $\beta\in\pi_2(X,L)$ and $\tilde f\in C^\infty(\widetilde{\mathrm{Gr}}_{k+1}\times L)$ domain-dependent Morse function associated to $f$, we denote by $\mathcal{M}(\underline{x},y;\beta,\tilde f):=\{((t, l), u)\}$ such that:
    \begin{enumerate}
        \item $(t, l)\in\mathrm{Gr}_{k+1}$.
        \item $u\in\mathcal{M}(t,\underline{x},y;\beta,\tilde f_t)$ with the understanding that $\tilde f_t\equiv\tilde f|_{\widetilde{\mathrm{Gr}}_{k+1}}$.
    \end{enumerate}
\end{defn}

\begin{prop}\label{OCPropII}
    For a generic choice of domain-dependent Morse function $\tilde f\in C^\infty(\widetilde{\mathrm{Gr}}_{k+1}\times L)$ associated to $f$, $\mathcal{M}(\underline{x},y;\beta,\tilde f):=\{((t, l), u)\}$ admits a Thickening $\T(\underline{x},y;\beta,\tilde f)$ together with a smooth submersion $\T(\underline{x},y;\beta,\tilde f)\rightarrow \mathrm{Gr}_{k+1}$.
\end{prop}

\begin{proof}
    This follows from the same argument of the proof of Theorem $4.44$ of \cite{rabah2024fukaya} (c.f. \cite{rabah2025fukaya}). Namely, we set $\T(\underline{x}, y;\beta,\tilde f):=\bigsqcup_{t\in G_{k+1}}\T(t, \underline{x}, y;\beta, \tilde f_t)$ where $\T(t, \underline{x}, y;\beta, \tilde f_t)$ is given by Corollary \ref{OCmoduliCor}. Indeed, we inductively define $\tilde f$ over the stratification of $\mathrm{Gr}_{k+1}=\bigcup_{t\in G_{k+1}}\mathrm{Gr}_{k+1}(t)$ and using the partial order $t>t'$ as described above. The base case follows from the fact that the stable/unstable manifolds of $(F,g_X)$ are transverse to the unstable/stable manifolds of $(f,g)$ and the submersion property of the evaluation map in Theorem $3.51$ of \cite{rabah2024fukaya}. The inductive step follows from Proposition \ref{OCmoduliProp} and the fact that the Banach space of all functions $\tilde f \in C^N(\widetilde{\mathrm{Gr}}_{k+1}\times L)$ such that $\tilde f|_{\widetilde{\mathrm{Gr}}_{k+1}(t)}\equiv\tilde f_t$ a pre-describe domain-dependent Morse function associated to $f$, is non-empty and closed. With this in mind, the smooth structure on $\T(\underline{x},y;\beta,\tilde f)$ is a consequence of the implicit function theorem. Indeed, a smooth coordinate chart around $((t,l),u)\in\T(\underline{x}, y;\beta,\tilde f)$ would be of the form of $V_\epsilon\times U$ where $V_\epsilon$ is neighborhood of zero in the normal fibre of $\mathrm{Gr}_{k+1}(t)\subset\mathrm{Gr}_{k+1}$ at $(t, l)$ and $u\in U\subseteq\T(t,\underline{x},y;\beta,\tilde f_t)$ is a coordinate chart around $u$ given by Corollary \ref{OCmoduliCor}. 

    As for the second claim, it follows from the second claim of Proposition \ref{OCmoduliProp}.
\end{proof}

\subsubsection{Global Kuranishi Charts}

The goal of this subsection is to give the $\textit{compactified}$ Open-Closed map moduli space a normally complex Global Kuranishi chart representation. We first give such construction for the moduli spaces in Definition \ref{OCmoduliDefn} and then $\textit{extend}$ our construction to the case where we allow broken Morse trajectories. Once again, we follow closely the constructions in \cite{rabah2024fukaya} and \cite{rabah2025fukaya}.

To avoid repetition, we fix $x_1,\dots,x_k\in\mathrm{Crit}(f)$, $y\in\mathrm{Crit}(F)$ and $\tilde f\in C^\infty(\widetilde{\mathrm{Gr}}_{k+1}\times L)$ a domain-dependent Morse function associated to $f$ as in Proposition \ref{OCPropII}. We also denote by $\T(\underline{x},y;\beta)\equiv\T(\underline{x}, y;\beta,\tilde f)$.

It is important to note that, the reason why we perturb the Morse function $f$ before the construction of the Global Kuranishi charts, as in Proposition \ref{OCPropII} for instance, is that to not include the Morse $\textit{perturbation}$ in our desired obstruction bundle. In fact, Theorem $3.51$ of \cite{rabah2024fukaya} (c.f. \cite{rabah2025fukaya}) inductively constructs Global Kuranishi chart representations for the moduli space of $J$-disks such that, if $\beta$ can be decomposed into a sum $\beta=\beta_1+\dots+\beta_j$ where each $\beta_i$ represents a $J$-disk, possibly constant $J$-disk, and $k+1=k_1+\dots+k_j$ then, the corresponding obstruction bundles given by Theorem $3.51$ of \cite{rabah2024fukaya} satisfies:
\[
\E_{k_1+1}(\beta_1)\oplus\dots\oplus\E_{k_j+1}(\beta_j)\hookrightarrow\E_{k+1}(\beta).
\]

\begin{defn}[Obstruction Bundle for the Open-Closed Moduli Space]\label{ObsBundOC}
    We define $\E(\underline{x},y;\beta)$ to be the pullback, via the projection map $\T(\underline{x},y;\beta)\rightarrow\mathrm{Gr}_{k+1}$ as given in Proposition \ref{OCPropII}, of the trivial bundle $\mathrm{Gr}_{k+1}\times\E_{k+1}(\beta)\rightarrow\mathrm{Gr}_{k+1}$. 
\end{defn}

\begin{defn}[Kuranishi Section for the Open-Closed Moduli Space]\label{KurSecOC}
    We define $s_{\underline{x},y;\beta}:\T(\underline{x},y;\beta)\rightarrow \E_{k+1}(\beta)$ by $s_{\underline{x},y;\beta}((t,l),u):=((t,l),s_{t,\underline{x},y;\beta}(u),0,\dots,0)$ where $s_{t,\underline{x},y;\beta}$ is the Kuranishi section of the Global Kuranishi chart given by Corollary \ref{OCmoduliCor} and the number of zeros is equal to the rank of the normal bundle $\mathrm{Gr}_{k+1}(t)\subseteq\mathrm{Gr}_{k+1}$ plus the rank of the normal bundle of $\E_{k_1+1}(\beta_1)\oplus\dots\oplus\E_{k_j+1}(\beta_j)\hookrightarrow\E_{k+1}(\beta)$ as described above.
\end{defn}

By construction, it follows that $\E(\underline{x},y;\beta)\rightarrow\T(\underline{x},y;\beta)$ is a smooth $O(d)$-equivariant vectorbundle, where $d:=\Omega(\beta)$ such that $\Omega$ is given by Lemma $3.5$ of \cite{rabah2024fukaya} or Lemma \ref{Existinceof2-form}. Moreover, by the same argument as in the proof of Lemma $4.47$ of \cite{rabah2024fukaya}, it follows that $s_{\underline{x},y;\beta}:\T(\underline{x},y;\beta)\rightarrow\E(\underline{x},y;\beta)$ is a well-defined continuous $O(d)$-equivariant section. Therefore, the tuple $(O(d),\T(\underline{x},y;\beta),\E(\underline{x},y;\beta),s_{\underline{x},y;\beta})$ serves as a Global Kuranishi chart for the moduli space $\mathcal{M}(\underline{x},y;\beta)$.

\begin{prop}\label{uncompSingleChart}
    The moduli space $\mathcal M(\underline{x},y;\beta,\tilde f)$ as given in Definition \ref{OCmoduliDefn}, admits a stably normally complex Global Kuranishi chart representation.
\end{prop}

\begin{proof}
    A candidate of such a Global Kuranishi chart is $(O(d),\T(\underline{x},y;\beta),\E(\underline{x},y;\beta),s_{\underline{x},y;\beta})$ as given above. The arguments presented in the proof of Theorem $4.48$ of \cite{rabah2024fukaya} shows that, possibly after stabilization and group enlargement, such Global Kuranishi chart admits a normal complex structure. Therefore, $(O(d),\T(\underline{x},y;\beta),\E(\underline{x},y;\beta,s_{\underline{x},y;\beta}))$ is a stably normally complex Global Kuranishi chart representation of $\mathcal{M}(\underline{x},y;\beta,\tilde f)$, as desired.
\end{proof}

To give a proof of Theorem \ref{TruncatedOCmap} and in view of Remark \ref{CompactificationRemark}, we will need to specify the $\textit{compactification}$ of $\mathcal M(\underline{x}, y;\beta,\tilde f )$ and give it a Global Kuranishi chart representation with some $\textit{compatibility}$ that would allow us to use the normally complex polynomial perturbations as in \cite{fukaya2001floer}, \cite{parker2013integral}, \cite{bai2022integral} and \cite{bai2022integral}.

We start by recalling some of the notions from Section $4$ of \cite{rabah2024fukaya}. We denote by $\overline{\mathrm{Gr}}_{k+1}$ the space of all broken, rooted, metric, ribbon trees as in Definition $4.11$ of \cite{rabah2024fukaya}. Following Definition $4.10$ of \cite{rabah2024fukaya}, each edge $e$ of a broken tree is decorated by a non-negative integer corresponding to the number of $\textit{breakings}$, which we will denote by $b(e)$. We also recall that a broken tree is constructed inductively out of $\textit{gluing}$ unbroken ones (c.f.  the discussion after Definition $4.11$ of \cite{rabah2024fukaya}). Moreover, we note that by Section $1.6$ of \cite{boardman2006homotopy} or Theorem $2.8$ of \cite{ma2010geometric}, $\overline{\mathrm{Gr}}_{k+1}$ has the structure of compact topological space homeomorphic to Stasheff's associahedra as in \cite{stasheff1963homotopy}, where its boundary is given by the image of the gluing map $\overline{\mathrm{Gr}}_{k_1+1}\times\overline{\mathrm{Gr}}_{k_2+1}\rightarrow \overline{\mathrm{Gr}}_{k+1}$ such that $k_1+k_2=k+1$. Now suppose that $(t, l)$ is a broken, rooted, metric, ribbon tree, formed $(t_1,l_1),\dots,(t_i, l_i)$ unbroken trees. Suppose that, for each of the $(t_i,l_t)$ we have chosen domain-dependent Morse function $\tilde f_i\in C^\infty(\widetilde{\mathrm{Gr}}_{k_i+1}\times L)$ associated to $f$ as in Corollary \ref{OCmoduliCor}. For each external edge $e$ of $t$ other than the one connected to the root $rt$, we will assign $I_e:=\R\sqcup\dots\sqcup\R\sqcup(-\infty, 0]$ where the number of $\R$-factors is equal to $b(e)$. Moreover, if $s$ denotes the coordinate on $\R\sqcup\dots\sqcup\R\sqcup(-\infty, 0]$ then, we will consider continuous maps $u_e:I_e\rightarrow L$ such that on each $\R$-factor $u_e$ solves the Morse trajectory equation of $f$
\[
\frac{d}{ds}u_e(s)+\nabla^gf(u_e(s))=0,
\]
with the condition that $\lim_{s\rightarrow-\infty}u_e(s)=x_i$ on the first $\R$-factor. While $u_e|_{(-\infty,0]}$ solves the domain-dependent Morse trajectory equation
\[
\frac{d}{ds}u_e(s)+\nabla^g\tilde f_i(s,u_e(s))=0.
\]
Similary, for the external edge connected to the root, we will assign $[0, \infty)\sqcup\R\sqcup\dots\sqcup\R$ where the number of $\R$-factors is equal to $b(e)$ and a continuous map $u_e:I_e\rightarrow X$ such that on each factor of $I_e$, $u_e$ solves the Morse trajectory equation of $(F, g_X)$. Namely,
\[
\frac{d}{ds}u_e+\nabla^{g_X}F(u_e(s))=0,
\]
with the condition that, on the last factor, we have $\lim_{s\rightarrow\infty}u_e(s)=y$. As for an internal edge $e$, we will assign $I_e=[0,\infty)\sqcup(-\infty,0]$ if $b(e)=1$ and if $b(e)\geq 2$, we will assign $I_e=[0,\infty)\sqcup\R\sqcup\dots\sqcup\R\sqcup(-\infty,0]$ where the number of $\R$ factors is $b(e)-1$. Similarly, we will also assign a continuous map $u_e:I_e\rightarrow L$ such that $u_e$ solves the Morse trajectory equation for $(f,g)$ on each $\R$-factor, while it solves the corresponding domain-dependent Morse trajectory equation on each of $[0,\infty)$ and $(-\infty,0]$ factors. Finally, we assign to each vertex $v$ of $t$ a $J$-disk $u_v:(\Sigma_v,\partial\Sigma_v)\rightarrow (X, L)$ representing the class $\beta_v$ with the incident conditions so that the resulting map $u:T\rightarrow X$ is continuous. We denote the moduli space of such maps by $\M(t, \underline{x},y;(\beta_v))$ and by $\M(t, \underline{x},y;\beta):=\bigsqcup_{\sum\beta_v=\beta}\M(t, \underline{x},y;(\beta_v))$.

\begin{defn}[Compactified Closed-Open Moduli Space]\label{CompOCmoduliDefn}
    We denote by $\M(\underline{x},y;\beta)$ the moduli space, for the $C^\infty_{loc}$-topology, of all tuples $((t, l),u)$ such that:
    \begin{enumerate}
        \item $(t,l)\in\overline{\mathrm{Gr}}_{k+1}$ is a broken, rooted, metric, ribbon tree.
        \item $u\in\M(t,\underline{x},y;\beta)$.
    \end{enumerate}
\end{defn}

Notice that, from the above description of $\M(\underline{x},y;\beta)$ and using standard gluing and compactness results of Morse trajectories (c.f. \cite{schwarz1993morse}), it follows that the codimension-$1$ $\textit{boundary}$ of $\M(\underline{x},y;\beta)$ is of the form
\begin{equation}\label{boundaryOCmoduli}
    \partial\M(\underline{x},y;\beta)=\bigcup\M(\underline{x}';\beta_1)\times\M(\underline{x}'',y;\beta_2),
\end{equation}

where $\underline{x}'\cup\underline{x}''=\underline{x}\cup\{z\}$ for some $z\in\mathrm{Crit}(f)$, $\beta_1+\beta_2=\beta$, and $\M(\underline{x}';\beta_1)$ denotes the compactified Floer-Morse trajectories passing through $\underline{x}'$ of degree $\beta_1$ as in Definition $4.49$ of \cite{rabah2024fukaya}.

\begin{prop}\label{ChartOCModuliCompact}
    Let $x_1,\dots,x_k\in\mathrm{Crit}(f)$, $y\in\mathrm{Crit}(F)$ and $\beta\in\pi_2(X,L)$. Then, for a generic choice of $\tilde f:\overline{\mathrm{Gr}}_{k+1}\rightarrow C^\infty(L)$ of domain-dependent Morse function associated to $f$, we have the following:
    \begin{enumerate}
        \item The moduli space $\M(\underline{x},y;\beta)$ as in Definition \ref{CompOCmoduliDefn}, admits an oriented normally complex Global Kuranishi chart representation $(G_\beta, \T(\underline{x},y;\beta),\E(\underline{x},y;\beta),s_{\underline{x},y;\beta})$ together with a smooth submersion $\T(\underline{x},y;\beta)\rightarrow \overline{\mathrm{Gr}}_{k+1}$.
        \item Each of the factors of $\partial\M(\underline{x},y;\beta)$ as given in \ref{boundaryOCmoduli}, admits an oriented normally complex Global Kuranishi chart representation $(G_{\beta_i},\T_i,\E_i,s_i)$ together with a smooth submersion $\T_i\rightarrow \overline{\mathrm{Gr}}_{k_i+1}$ such that:
        \begin{enumerate}
            \item $G_{\beta_1}\times G_{\beta_2}\hookrightarrow G_{\beta}$.
            \item If $(G_{\beta_i}, \T_i,\E_i,s_i)$ is a Global Kuranishi chart representation of $\M(\underline{x}_1;\beta_1)$ then, it is the one given by Theorem $4.50$ of \cite{rabah2024fukaya}.
            \item The following diagram commutes:
             \[\begin{tikzcd}
	{\mathcal{E}(\underline x_1; \beta_1)\boxplus\mathcal{E}(\underline x_2,y; \beta_2)} && {\mathcal{E}(\underline x,y;\beta)} \\
	\\
	{\mathcal{T}(\underline x_1; \beta_1)\times\mathcal{T}(\underline x_2, y; \beta_2)} && {\mathcal{T}(\underline x,y; \beta)} \\
	\\
	{\overline{Gr}_{k_1+1}\times\overline{Gr}_{k_2+1}} && {\overline{Gr}_{k+1}}
	\arrow[hook, from=1-1, to=1-3]
	\arrow[from=1-1, to=3-1]
	\arrow[from=1-3, to=3-3]
	\arrow["{s_{\underline x_1; \beta_1}\boxplus s_{\underline x_2, y; \beta_2}}", bend left = 30pt, from=3-1, to=1-1]
	\arrow[hook, from=3-1, to=3-3]
	\arrow[from=3-1, to=5-1]
	\arrow["{s_{\underline x, y; \beta}}"', bend right = 30pt, from=3-3, to=1-3]
	\arrow[from=3-3, to=5-3]
	\arrow["glue"', from=5-1, to=5-3]
\end{tikzcd}\]
and compatible with the normal complex structure and orientations, as in Remark $4.51$ of \cite{rabah2024fukaya}.
        \end{enumerate}
    \end{enumerate}
\end{prop}

\begin{proof}
    We note that for purposes, it suffice to just consider the case when the expected dimension of $\M(\underline{x},y;\beta)$ is $1$. Namely, when $\mu(\underline{x},y,\beta):=|y|-\sum_i|x_i|+\mu_L(\beta)+k-n=1$. Nevertheless, we still have to consider the moduli spaces arising in our compactification even in the case their expected dimension is negative. 
    
    We first start under the assumption that the number of $\textit{breaking}$ is $1$ and argue by an upward two-layered induction over the energy $\omega(\beta)$ and the number of leaves minus $1$, that is $k$.\\
    The base case follows from Proposition \ref{uncompSingleChart} and the fact that the stable/unstable manifolds of $(F,g_X)$ are transversal to the the unstable/stable manifolds of $(f,g)$.\\
    As for the inductive step, suppose that $\beta_1+\beta_2=\beta$, $k_1+k_2=k+1$ and for $i=1,2$ we have $\tilde f_i\in C^\infty(\widetilde{\mathrm{Gr}}_{k_i+1}\times L)$ a domain-dependent Morse function associated to $f$ such that $\M(\underline{x}_1;\beta_1)$ has a Global Kuranishi chart representation as given in Theorem $4.50$ of \cite{rabah2024fukaya}, while $\M(\underline{x}_2,y;\beta_2)$ has a Global Kuranishi chart representation as given in Proposition \ref{uncompSingleChart}. As each of $\tilde f_i$ is a domain-dependent Morse function associated to $f$ and in particular, $\tilde f_i\equiv f$ on a small enough neighborhood of $\mathrm{Crit}(f)$, it follows by Whitney's Extension Theorem \cite{whitney1934analytic} that there exists an $\tilde f\in C^\infty(\widetilde{\mathrm{Gr}}_{k+1}\times L)$ such that the pullback of $\tilde f$ by the gluing map $\mathrm{Gr}_{k_1+1}\times \mathrm{Gr}_{k_2+1}\rightarrow\overline{\mathrm{Gr}}_{k+1}$ is $\tilde f_i$ for each $i=1,2$. In particular, it follows that the Banach space all elements of $C^N(\overline{\mathrm{Gr}}_{k+1}\times L)$ that pullback to $\tilde f_i$ correspondingly under the above gluing map, is closed and non-empty. Therefore, the arguments of the proof Proposition \ref{uncompSingleChart} go through and we can construct a Global Kuranishi chart representation $(G_\beta, \T(\underline{x},y;\beta), \E(\underline{x},y;\beta),s_{\underline{x},y;\beta})$ of $\M(\underline{x},y;\beta)$ using such $\tilde f $ such that $G_{\beta_1}\times G_{\beta_2}\hookrightarrow G_{\beta}$ corresponding to the diagonal embedding of $O(d_1)\times O(d_2)\hookrightarrow O(d)$ where $d_1+d_2=d$. Now as $\tilde f$ is a domain-dependent Morse function associated to $f$ and since it pulls back by the gluing map to $\tilde f_i$ by construction, it follows using the same Morse trajectory gluing arguments presented in the proof of Theorem $4.50$ of \cite{rabah2024fukaya} that for such choice of $(G_\beta, \T(\underline{x},y;\beta), \E(\underline{x},y;\beta),s_{\underline{x},y;\beta})$ we have the following commutative diagram:
     \[\begin{tikzcd}
	{\mathcal{T}(\underline{x}_1;\beta_1)\times\mathcal{T}(\underline x_2,y; \beta_2)} && {\mathcal{T}(\underline x, y;\beta)} \\
	\\
	{Gr_{k_1+1}\times Gr_{k_2+1}} && {\overline{Gr}_{k+1}}
	\arrow[hook, from=1-1, to=1-3]
	\arrow[from=1-1, to=3-1]
	\arrow[from=1-3, to=3-3]
	\arrow[from=3-1, to=3-3]
\end{tikzcd}\]
We note that, as for now $\T_i\rightarrow \mathrm{Gr}_{k_i+1}$ is a smooth submersion, while $\T(\underline{x}_1;\beta_1)\times \T(\underline{x}_2,y;\beta_2)\hookrightarrow\T(\underline{x},y;\beta)$ is a topological embedding. On the other hand, using the fact that $G_{\beta_1}\times G_{\beta_2}\hookrightarrow G_\beta$ and adhering Theorem $6.10$ or Appendix B Theorem $B.3$ of \cite{bai2022arnold} with the understanding that $G_{\beta_i}, G_\beta$ action is trivial on $\mathrm{Gr}_{k_i+1}, \overline{\mathrm{Gr}}_{k+1}$ respectively, it follows that, possibly after further stabilizing $(G, \T(\underline{x},y;\beta), \E(\underline{x},y;\beta), s_{\underline{x},y;\beta})$ by a real $O(d)$-representation and some abuse of notation, we can assume that the above diagram commutes as smooth manifolds. Moreover and as discussed above, due to the fact that the relevant obstruction bundles have no Morse perturbations, it follows by construction, namely Definitions \ref{ObsBundOC} and \ref{KurSecOC} that we have the following commutative diagram:
\[\begin{tikzcd}
	{\mathcal E(\underline{x}_1;\beta_1)\oplus\mathcal E(\underline{x}_2,y;\beta_2)} && {\mathcal E(\underline{x};\beta)} \\
	\\
	{\mathcal T(\underline{x}_1;\beta_1)\times\mathcal T(\underline{x}_2,y;\beta_2)} && {\mathcal T(\underline{x},y;\beta)}
	\arrow[hook, from=1-1, to=1-3]
	\arrow[from=1-1, to=3-1]
	\arrow[from=1-3, to=3-3]
	\arrow["{s_{\underline{x}_1;\beta_1}\oplus s_{\underline{x}_2,y;\beta_2}}", bend left=30pt, from=3-1, to=1-1]
	\arrow[hook, from=3-1, to=3-3]
	\arrow["{s_{\underline{x},y;\beta}}"', bend right=30pt, from=3-3, to=1-3]
\end{tikzcd}\]
As for the compatibility of orientation, this follows by using the same orientation construction as in Proposition $8.1.14$ of \cite{fukaya2010lagrangian} or Remark $4.34$ of \cite{charest2015floer} (c.f. \cite{fukaya2013lagrangian}). While the compatibility of the normal complex structures, follows from the same inductive argument as in the proof of Theorem $4.50$ of \cite{rabah2024fukaya} (c.f \cite{rabah2025fukaya}).

To finish the proof, we argue by induction on the number of $\textit{breakings}$. The base case follows from Proposition \ref{uncompSingleChart}, while the inductive step follows from the above arguments.
\end{proof}

\begin{proof}[Proof of Theorem \ref{TruncatedOCmap}]
    Having Proposition \ref{ChartOCModuliCompact} in mind, for $x_1,\dots, x_k\in\mathrm{Crit}(f)$ and $\beta\in\pi_2(X,L)$ we set
    \[
    \p_{k,\beta}(x_1,\dots, x_k):=\sum_{y\in\mathrm{Crit}(F)}(-1)^*|s_\epsilon^{-1}(0)\cap\mathcal U_{free}|y,
    \]
    where the sum is taken over all $y\in\mathrm{Crit}(F)$ such that $|y|-\sum_i|x_i|+\mu_L(\beta)+k-n=0$, $(\mathcal{U}, E, \mathcal S)$ is the derived orbifold chart associated to the Global Kuranishi chart representation of $\M(\underline{x},y;\beta)$ as given in Proposition \ref{ChartOCModuliCompact} and $s_\epsilon$ is an FOP-perturbation of $\mathcal S$ such that $||\mathcal{S}-s_\epsilon||_{C^0}<\epsilon$ as in \cite{bai2022integral}, \cite{bai2022arnold}. Now set
    \[
    \p_k:=\sum_{\beta\in\Pi(L)}T^{\omega(\beta)}e^\frac{\mu(\beta)}{2}\p_{k,\beta}
    \]
    and extend $\Lambda_0$-linearily.
    It is important to note that the construction of $s_\epsilon$ is done inductively as follows:

    We construct $s_\epsilon$ by upward induction on the energy $\omega(\beta)$ and $k$ the number of leaves minus using the partial order $\precsim$ on $\mathfrak{G}$ as given in \cite{rabah2024fukaya}, \cite{Fukaya2009} using the same arguments in the proof of Theorem $5.3$ of \cite{rabah2024fukaya} or \cite{bai2022arnold}, with the requirement that FOP-perturbation chosen on $\M(\underline{x}_1;\beta_1)$ in the boundary of $\M(\underline{x},y;\beta)$ is the one given by Theorem $5.3$ of \cite{rabah2024fukaya}. 
    
    As for the $\textit{sign check}$, we $\textit{translate}$ the results of Section $15.4$ of \cite{fukaya2013lagrangian} into our setting, as we are following their orientation convention. To be more precise, let $x_1,\dots,x_k\in\mathrm{Crit}(f)$, $y\in\mathrm{Crit}(F)$ and $\beta\in\pi_2(X, L)$ such that $\mu(\underline{x},y,\beta)=1$. We have the following cases:
    \begin{enumerate}
        \item Ambient Morse trajectory breaking. This corresponds to the boundary component as in Proposition \ref{ChartOCModuliCompact} of the form $\T(\underline{x}_1;\beta_1)\times\T(\underline{x}_2,y;\beta_2)$ where $\beta_2=0$, $\underline{x}_1=\underline{x}$. Such degeneration contributes by $(+1)y$ in the expression of $d_X\circ\p_{k,\beta}(x_1,\dots,x_k)$.
        \item Degenerations of the form 
        \[
        \bigsqcup_{z\in\mathrm{Crit}(f)}\T(x_i,\dots,x_{i+j-1},z;\beta_1)\times \T(x_1,\dots,x_{i-1},z,x_{i+j},\dots,x_k,y;\beta_2)\]
        which will contribute by $(-1)^{\epsilon_i}y$ in the expansion of 
        \[
        \p_{k-j,\beta_2}(x_1,\dots,x_{i-1},\m_{j,\beta_1}(x_i,\dots,x_{i+j-1}),x_{i+j},\dots,x_k).
        \]
        \item Lastly, using the ribbon structure on our underlying trees and considering configurations starting at the internal broken edge, we will consider configurations of the form
        \[
        \bigsqcup_{z\in\mathrm{Crit}(f)}\T(x_{k-l+1},\dots,x_k,x_1,\dots,x_s,z;\beta_1)\times\T(z,x_{s+1},\dots,x_{k-l},y;\beta_2),
        \]
        which will contribute by $(-1)^{\epsilon_l}y$ in the expansion of 
        \[
        \p_{k-l+s,\beta_2}(\m_{l+s,\beta_1}(x_{k-l+1},\dots,x_k,x_1,\dots,x_s),x_{s+1},\dots,x_{k-l}).
        \]
    \end{enumerate}
    Now the result follows after using Lemma $15.9$ and Lemma $15.10$ of \cite{fukaya2013lagrangian}.
\end{proof}

\subsection{Proof of Theorem \ref{OCmapTheorem}}
The proof that we will present here, follows the general recipe as in \cite{Fukaya2009} and \cite{fukaya2010lagrangian}. Namely, we will use a $\textit{homotopy limit argument}$ on the operators given by Theorem \ref{TruncatedOCmap}, to get our desired result. In fact, we will adapt the results of Section $15.1.7$ of \cite{fukaya2017kuranishi} to our setting and not use the full generality of their work.

We start by recalling some notations from Section $5$ of \cite{rabah2024fukaya} which are inline with the conventions in \cite{Fukaya2009} and \cite{fukaya2010lagrangian} (c.f. \cite{rabah2025fukaya}). Set $\bar C:=\mathrm{CM}^*(f;\Z)$, the $\Z$-reduction of our $A_\infty$-algebra of this Section and let $\bar B:=\mathrm{CM}^*(F;\Z)[n]$ where $n=\dim_\R L$. We set $\delta_X:=(-1)^nd^\text{Morse}_X$ and denote by $B:=\bar B\widehat{\otimes}_\Z\Lambda_0$ its degree-wise completion and we linearly extended $\delta_X$ to a differential on $B$. Moreover, recall that we set $B_k\bar C[1]\equiv C[1]^{\otimes k}$. Note that, our $A_\infty$-algebra is $\mathfrak{G}$-gapped, where 
\[
\mathfrak{G}:=\{(\omega(\beta),\mu(\beta)):\beta\in\Pi(L), \mathcal M(\beta)\neq\emptyset\}
\]
and we have assigned a pre-order $\precsim$ on $\mathfrak G\times\N$ defining a partial order on its set of equivalence classes $\mathfrak G\times\N/\sim$. Let $a\geq1$ be an indexing integer and note that by Theorem $5.3$ of \cite{rabah2024fukaya}, we have that for every $(\beta,j)\precsim(n_a,K_a)$ we have linear maps
\[
\m_{j,\beta}^{(a)}:B_j\bar C[1]\rightarrow \bar C[1]
\]
of degree $1-\mu(\beta)$. Moreover, by Theorem \ref{TruncatedOCmap}, we have that for every $(\beta, k)\precsim(n_a, K_a)$ such that $k\geq1$, we have linear maps
\[
\p_{k,\beta}^{(a)}:B_k\bar C[1]\rightarrow \bar B
\]
of degree $-\mu(\beta)$.

Theorem $3.2$ above can be rewritten as 
\[
\delta_X\circ\p^{(a)}_{k,\beta}+\p_{k,\beta}^{(a)}\circ\delta_H^{\m^{(a)}}=0,
\]
for every $(\beta, k)\precsim(n, K)$ such that $k\geq 1$, with the understanding that $\delta_H^{\m^{(a)}}$ is the Hochschild differential given by the coderivation induced from $\m^{(a)}$. It is important to note that the terms in the expansion of $\p_{k,\beta}\circ\delta_X$ are linear maps of the above form for $(k',\beta')\precsim(n, K)$. Indeed, using Section $5$ of \cite{rabah2024fukaya} and the definition of $||\beta||$ there (c.f. \cite{fukaya2010lagrangian}), it is easy to see that if $\beta=\beta_1+\beta_2$ such that $\beta_1\neq 0$ then, $||\beta||\geq ||\beta_2||+1$. In particular, since the curvature term $\m_0(1)$ has positive-valuation, it follows that the operators of Theorem $3.2$ above are filtration-preserving. 

Now following Section $15.4$ of \cite{fukaya2013lagrangian}, we set
\[
\widehat{CC}(C[1]):=\varprojlim_{E\to\infty}
\left[
\left(\bigoplus_{k\ge1}B_k\bar C[1]\right)
\otimes_{\mathbb Z}
\left(\Lambda_0/T^E\Lambda_0\right)
\right].
\]

\begin{lemma}(c.f. Lemma 7.2.69 of \cite{fukaya2010lagrangian})\label{OChomotopyLimit}
    Suppose that, for any $b\geq a$ we have:
    \begin{enumerate}
        \item $\m^{(b)}_{j,\beta}\equiv\m^{(a)}_{j, \beta}$ for all $(\beta, j)\precsim(n_a, K_a)$.
        \item $\p^{(b)}_{k,\beta}\equiv\p^{(a)}_{k,\beta}$ for all $(\beta, k)\precsim(n_a, K_a)$.
    \end{enumerate}
    Then, $\p_k:=\sum_{\beta\in\mathfrak{G}}T^{\omega(\beta)}e^{\frac{\mu(\beta)}{2}}\p_{k,\beta}$ defines a degree-zero, filtered chain map
    \[
    \p:(\widehat{CC}(C[1]), \delta_H)\rightarrow (B, \delta_X)
    \]
    such that $\delta_X\circ\p+\p\circ\delta_H=0$.
\end{lemma}

\begin{proof}
    For any $a\geq 1$ indexing integer and by definition of $\p_k$ and using items $(1)$ and $(2)$ of the statement, we have that $\m_{j,\beta}=\m_{j,\beta}^{(a)}$ for every $(\beta,j)\precsim(n_a, K_a)$ and $\p_{k,\beta}=\p_{k,\beta}$ for every $(\beta, k)\precsim(n_a, K_a)$ such that $k\geq 1$. Now set $\p:=\sum_{k\geq1}\p_k$ and argue as follows:
    
    Notice that, since the operators $\p_{k,\beta}$ given by Theorem $3.2$ are filteration-preserving, it follows that $\p$ is filteration-preserving. In particular, as $B$ is complete and seperated with respect to the $T$-adic topology on $\Lambda_0$, it follows that $\p$ is well-defined and continuous with respect to the $T$-adic topology. 
    
    To show the identity $\delta_X\circ\p+\p\circ\delta_H=0$ it is important to note that we are following the differential convention for $\delta_H$ as given in Section $15.4$ of \cite{fukaya2013lagrangian}. Now fix an energy bound $E\geq 0$ and a word length of the elements of $\widehat{CC}(C[1])$. Fix a pair $(\beta, k)$ and an indexing integer $a\geq1$ such that $(\beta,k)\precsim(n_a, K_a)$, where $K_a$ is less than or equal to the fixed word length and $n_a\leq E$. As the terms in the expansion of $\delta_X\circ\p_{k,\beta}$ have operators indexed by $(\beta',k')$ such that $(\beta', k')\precsim(n_a, K_a)$, it follows by Theorem $3.2$ that $\delta_X\circ\p^{(a)}_{k,\beta}+\p_{k,\beta}^{(a)}\circ\delta_H^{\m^{(a)}}=0$. By continuity of the operators $\delta_X, \p, \delta_H$ with respect to the $T$-adic topology and seperatedness of both $B$ and $\widehat{CC}(C[1])$, it follows that $\delta_X\circ\p+\p\circ\delta_H=0$ as desired. 
\end{proof}

\begin{proof}[Proof of Theorem \ref{OCmapTheorem}]
    Having the construction of the $A_\infty$-algebra $(\mathrm{CM}^*(f; \Z),\{\m_k\}_{k\geq1})$ as given in Section $5$ of \cite{rabah2024fukaya} in mind, the result now follows from Theorem \ref{TruncatedOCmap} and Lemma \ref{OChomotopyLimit}.
\end{proof}

	\section{$A_\infty$-bimodules}

        \subsection{Algebraic Preliminaries}
    Let $(C_0, \m^0)$ and $(C_1, \m^1)$ be two filtered, possibly curved, $A_\infty$-algebras over $\Lambda_0$. Let $\bigoplus_{m\in\Z}D^m$ be a graded free filtered $\Lambda_0$-module. We require the following conditions to be satisfied our filtration $F$ (c.f. 3.2.2 of \cite{Fukaya2009}):
    \begin{enumerate}
        \item $\lambda>\lambda'\implies F^\lambda D^m\subset F^{\lambda'}D^m$,
        \item $T^{\lambda_0}. F^\lambda D^m\subset F^{\lambda+\lambda_0}D^m$,
        \item $e^nD^m\subset D^{m+2n}$,
        \item $D^m$ is complete with respect to the $F$-filtration and 
        \item $D^m$ has a basis $d_i$ such that $d_i\in F^0D^m\setminus\bigcup_{\lambda>0}F^\lambda D^m$.
    \end{enumerate}

\begin{remark}
    In our case, $D$ will be finitely-generated and hence items $(4)$ and $(5)$ will be automatically satisfied.
\end{remark}

Let $D:=\widehat{\bigoplus}_{m\in\Z} D^m$ be the completion with respect to $F$. For $k_0, k_1\geq0$ integers we consider the following $\Lambda_0$-linear map:
\begin{equation}
    n_{k_1, k_0}:B_{k_1}(C_1[1])\otimes_{\Lambda_0}D[1]\otimes_{\Lambda_0}B_{k_0}(C_0[1])\rightarrow D[1],
\end{equation}
        where $\deg n_{k_1,k_0}=1$ and $n_{k_1, k_0}$ is $F$-filtration preserving. That is,
        \[
        n_{k_1,k_0}(F^{\lambda_1}(C_1[1])^{m_1}\otimes\dots\otimes F^{\lambda_{k_1}}(C_1[1])^{m_{k_1}}\otimes F^{\lambda_0}D[1]\otimes F^{\lambda_{m'_1}}(C_0[1])^{m'_1}\otimes\dots\otimes F^{\lambda_{m'_{k_0}}}(C_0[1])^{m'_{k_0}})
        \]
\[
\subset F^{\sum\lambda_i+\lambda_0+\sum\lambda_i'}(D[1])^{\sum m_i+m_0+\sum m'_i}.
\]
With this in mind, we extend $n_{k_1, k_0}$ to the completion
\[
\n:\widehat{B}(C_1[1])\hat{\otimes}_{\Lambda_0}D[1]\hat\otimes_{\Lambda_0}\widehat{B}(C_0[1])\rightarrow D[1]
\]
of $\bigoplus_{k_0, k_1}B_{k_1}(C_1[1])\otimes_{\Lambda_0}D[1]\otimes_{\Lambda_0}\bigoplus_{k_0, k_1}B_{k_0}(C_0[1])$. In this case, we say $(D[1], \n)$ is a left $(C_1[1], \m^1)$ filtered $\Lambda_0$-module and a right $(C_0[1], \m^0)$ filtered $\Lambda_0$-module.

To give $(D[1], \n)$ the structure of an $A_\infty$-bimodule and following Section $3.7$ of \cite{Fukaya2009}, we first consider the (formal) adjoint of $\n$. Namely, let
\[
\hat{d}:\widehat{B}(C_1[1])\hat{\otimes}_{\Lambda_0}D[1]\hat{\otimes}_{\Lambda_0}\widehat{B}(C_0[1])\rightarrow\widehat{B}(C_1[1])\hat{\otimes}_{\Lambda_0}D[1]\hat{\otimes}_{\Lambda_0}\widehat{B}(C_0[1])
\]
be given as follows:\\
For $k_0, k_1\geq0$ integers and $x_{1,1}\otimes\dots\otimes x_{1, k_1}\in \widehat{B}(C_1[1])$, $y\in D[1]$ and $x_{0,1}\otimes\dots\otimes x_{0, k_0}\in \widehat{B}(C_0[1])$, we set
\[
\hat{d}(x_{1,1}\otimes\dots\otimes x_{1, k_1}\otimes y\otimes x_{0,1}\otimes\dots\otimes x_{0, k_0})
\]
\[
=\sum_{k_1'\leq k_1, k_0'\leq k_0}(-1)^{*_1}x_{1,1}\otimes\dots\otimes x_{1, k_1-k_1'}\otimes n_{k'_1, k_0'}(x_{1, k_1'-k_1+1},\dots,y, x_{0,1}, \dots, x_{0, k_0'})\otimes x_{0, k_0'+1}\otimes\dots\otimes x_{0, k_0}
\]
\[
+\hat{d}^1(x_{1, 1}\otimes\dots\otimes x_{1, k_1})\otimes y\otimes x_{0,1}\otimes\dots\otimes x_{0, k_0}+(-1)^{*_2}x_{1,0}\otimes\dots x_{1, k_1}\otimes y\otimes\hat d^0(x_{0,1}\otimes\dots\otimes x_{0, k_0}),
\]
where $*_1=\deg x_{1,1}+\dots+\deg x_{1, k_1-k_1'}+(k_1-k_1')$ and $*_2=\deg x_{1,1}+\dots+\deg x_{1, k_1}+\deg y+ k_1+1$, and $\hat{d}^i$ is given in Equation $(3.7.2)$ of \cite{Fukaya2009} (c.f. Section $5.1$ of \cite{rabah2024fukaya}).

\begin{defn} [$A_\infty$-bimodule]
    $(D[1], \n)$ is said to be an $A_\infty$-bimodule over $(C_0[1], \m^0)$ and $(C_1[1], \m^1)$ if $\hat{d}\circ\hat{d}\equiv 0$.
\end{defn}

For our geometric purposes, we will also need the notion of $\mathfrak{G}$-gapped $A_\infty$-bimodules. Indeed, let $\mathfrak{G}\subseteq \R\times\Z$ be a submonoid and suppose that $(D[1], \n)$ is a filtered $A_\infty$-bimodule over $(C_i, \m^i)$ for $i=0,1$.

\begin{defn}[$\mathfrak{G}$-gapped $A_\infty$-bimodule]\label{AinftyBimodule}
    $(D[1], \n)$ is said to be a $\mathfrak{G}$-gapped $A_\infty$-bimodule if for any $k_0, k_1\geq 0$ integers, there exists a $\Z$-linear map $n_{k_1, k_0, \beta}:B_{k_1}(\overline{C}_1[1])\otimes_\Z\overline{D}[1]\otimes_\Z\overline{B}_{k_0}(\overline{C}_0[1])\rightarrow \overline{D}[1]$, on the $\Z$-reduction of $D[1]$ and $C_i[1]$ for $i=0,1$ such that:
    \begin{equation}
        n_{k_1, k_0}=\sum_{\beta\in\mathfrak{G}}T^{\omega(\beta)}e^{\frac{\mu(\beta)}{2}}n_{k_1, k_0, \beta}.
    \end{equation}
\end{defn}

Now suppose that, for $i=0,1$, $(C_i[1], \m^i)$ is a filtered, $\mathfrak{G}_i$-gapped $A_{n, K}$-algebra, as in Definition $7.2.67$ of \cite{fukaya2010lagrangian} or Section $5.1$ of \cite{rabah2024fukaya}. Let $G$ be a submonoid of $\R\times \Z$ and define
\[
\mathfrak{G}:=\{\beta_1+\alpha+\beta_0:\beta_i\in\mathfrak{G}_i, \alpha\in G\}.
\]
Following Section $7.2.6$ of \cite{fukaya2010lagrangian}, namely Definitions $7.2.61-7.2.63$, we introduced in Section $5.1$ of \cite{rabah2024fukaya} a partial order $\precsim$ on the set of equivalences classes $\mathfrak{G}_i\times\N/\sim$. Using the definition of $\mathfrak{G}$ as given above, we extend $\precsim$ to a partial order on the set of equivalences classes on $\mathfrak{G}\times \N$.

\begin{defn}[Filtered, Gapped $A_{n, K}$-bimodule]\label{AnKbimoduleDefn}
    For $i=0,1$, let $(C_i[1], \m^i)$ be filtered, $\mathfrak G_i$, $A_{n, K}$-bimodules and let $\mathfrak{G}$ be as above. Let $D$ be a graded, filtered $\Lambda_0$-module. An $A_{n, K}$-bimodule structure on $D$ is a sequence of filtration-preserving $\Lambda_0$-linear maps
    \begin{equation}
        n_{k_1, k_0, \beta}:B_{k_1}(C_1[1])\otimes D[1]\otimes B_{k_0}(C_0[1])\rightarrow D[1],
    \end{equation}
    indexed by $\mathfrak{G}\times \N$ such that if we consider its formal adjoint $\hat{d}_{k_1, k_0, \beta}$ as given above, we require that $\hat{d}_{k_1, k_0, \beta}\circ\hat{d}_{k_1, k_0, \beta}=0$ whenever $(\beta, k_0+k_1+1)\precsim (n, K)$. 
\end{defn}

As done in \cite{rabah2024fukaya}, \cite{rabah2025fukaya}, we will follow the $\textit{general recipe}$ following \cite{Fukaya2009}, \cite{fukaya2010lagrangian} by geometrically constructing the operations defining a structure of an $A_{n,K}$-bimodule structure and then use $\textit{homotopy limit arguments}$ as in Lemma $7.2.69$ and Theorem $7.2.72$ of \cite{fukaya2010lagrangian}. Nevertheless, since in this paper we are interested in defining just a Lagrangian Floer modules for a given transversaly intersecting pairs of Lagrangians $L_0,L_1$ we prove the following:

Set, $D:=\widehat{\bigoplus}_{p\in L_0\cap L_1}\Lambda_0\langle p\rangle$ and for a real number $E>0$, we denote by $D_E:=D/F^E D$. We will assume that $D$ is complete and seperated with respect to the $T$-adic topology induced by $F$. Such assumption is a consequence of the $\textit{refined}$ filtration that we use in the Spectral Sequences Section \ref{spectralSeqSec}.

\begin{prop}
    Let $\{\lambda_j\}_{j\in\N}$ be an increasing sequence of non-negative real numbers such that $\lambda_j\rightarrow\infty$ as $i\rightarrow\infty$. Suppose that:
    \begin{enumerate}
        \item For every $j\in\N$, there exists a degree-$1$ map $d_j:D_{\lambda_j}\rightarrow D_{\lambda_j}$ satisfying $d_j\circ d_j\equiv 0$.
        \item For every $j\in\N$, there exists a filtered chain map $f_j:(D_{\lambda_j}, d_j)\rightarrow (D_{\lambda_j},\bar d_{j+1})$, that is $\bar d_{j+1}\circ f_j=f_j\circ d_j$ where $\bar d_{j+1}:= d_{j+1}\text{ mod }F^{\lambda_j}\Lambda_0$; such that $f_j=\Id+N_j$, where $N_j$ has entries of positive-valuation.
    \end{enumerate}
    Then, there exists a unique differential $\delta:D\rightarrow D$ such that for every $j\in\N$, $\delta_j:=\delta\text{ mod }F^{\lambda_j}\Lambda_0$ is congruent to $d_j$.
\end{prop}

\begin{remark}
    We say $\delta_j$ is congruent to $d_j$ if there exists an invertible filtered map $P_j:D_{\lambda_j}\rightarrow D_{\lambda_j}$ such that $\delta_j=P^{-1}_j\circ d_j\circ P_j$.
\end{remark}

\begin{proof}
    Noting that $L_0$ and $L_1$ intersect transversely, it follows that $D$ is finitely-generated over $\Lambda_0$. On the other hand and for a fixed $j\in\N$, and since the entries of $N_j$ have positive-valuation, it follows that there exists a $\gamma_j>0$ such that $\gamma_j$ is less than or equal to the valuations of all entries of $N_j$. In particular, $N_j$ is idempotent in $D_{\lambda_j}$. That is, there exists an $E_j>0$ such that for $k\gamma_j\geq \lambda_j$, we have 
    \[
    N_j^k|_{D_{\lambda_j}}=N_j\dots N_j|_{D_{\lambda_j}}=0\in D_{\lambda_j}.
    \]
    Therefore, $f_j$ is invertible having $\sum_{k\geq0}(-1)^kN_j^k$ as an inverse. With this in mind, we can inductively construct invertible filtered maps $P_j:D_{\lambda_j}\rightarrow D_{\lambda_j}$ having $P_1:=\Id$.

    Now set, $\delta_j:=P^{-1}_j\circ d_j\circ P_j$ and note that, 
    \[
    \delta_{j+1}\text{ mod }F^{\lambda_j}\Lambda_0=(f_j\circ P_j)^{-1}\circ\bar d_{j+1}\circ (f_j\circ P_j)=\delta_j.
    \]
    Thus, using completeness and seperatedness of $D$ with respect to the $T$-adic topology given by the $F$-filtration, it follows that setting $\delta:=\varprojlim_j\delta_j$ completes the proof.
\end{proof}

\subsection{Geometric Realization}

Let $L_0, L_1$ be a pair of transversely intersecting closed, connected, relatively-spin Lagrangians. Denote by 
\[
\Omega(L_0, L_1):=\{\ell:[0,1]\rightarrow X\text{ smooth map}:\ell(0)\in L_0, \ell(1)\in L_1\}
\]
and fix a $\textit{base path}$ $\ell_0\in\Omega(L_0, L_1)$. Denote by $\Omega(L_0, L_1; \ell_0)$ the connected component of $\Omega(L_0, L_1)$ containing $\ell_0$. To describe the $\textit{Novikov}$ covering, that is a covering space with an abelian group of deck of transformations, of $\Omega(L_0, L_1; \ell_0)$, we consider pairs $(\ell, w)$ such that $\ell\in \Omega(L_0, L_1; \ell_0)$ and $w$ is a homotopy of paths in $\Omega(L_0, L_1; \ell_0)$ between $\ell$ and $\ell_0$. In equations, we require that $w:[0,1]_s\times[0,1]_t\rightarrow X$ is a smooth map such that:
\begin{enumerate}
    \item $w(0,.)\equiv\ell_0$ and $w(1, .)\equiv\ell$,
    \item $w(s, 0)\in L_0$ and $w(s,1)\in L_1$ for all $s\in[0,1]$.
\end{enumerate}

Now suppose that we are given two such pairs of the the form $(\ell, w)$ and $(\ell, w')$. We consider the concatenation $\Bar{w}\# w':[0,1]\times [0,1]\rightarrow X$ and notice that, as both $w, w'$ have the same base path, we can think of $\Bar{w}\# w'$ as a map $C:S^1_s\times[0,1]_t\rightarrow X$ such that, $C(s, 0)\in L_0$ and $C(s, 1)\in L_1$ for all $s\in S^1$.

\begin{defn}
    We define the energy functional $I_\omega:\pi_1(\Omega(L_0, L_1; \ell_0))\rightarrow\R$ by $C\mapsto \int _C\omega$.
\end{defn}

Now for $i=0,1$, we denote by $c_i\equiv C(.,i)$ and consider the bundle pair $(C^*TX, c_0^*TL_0\sqcup c_1^*TL_1)\rightarrow S^1\times [0,1]$.

\begin{defn}
    We define the index functional $I_\mu:\pi_1(\Omega(L_0, L_1; \ell_0))\rightarrow \Z$ by $C\mapsto \mu(C^*TX, c_0^*TL_0\sqcup c_1^*TL_1)$, where the right-hand side is as in Definition $2.1.15$ of \cite{Fukaya2009}.
\end{defn}

\begin{defn}
\begin{enumerate}
    \item We say $(\ell, w)$ and $(\ell, w')$ are equivalent and write $(\ell, w)\sim(\ell, w')$ if $I_\omega(\Bar{w}\# w')=0= I_\mu(\Bar{w}\# w')$.
    \item We denote by $\widetilde\Omega(L_0, L_1; \ell_0):=\Omega(L_0, L_1; \ell_0)/\sim$, which defines a covering space $\widetilde\Omega(L_0, L_1; \ell_0)\rightarrow\Omega(L_0, L_1; \ell_0)$ of deck of transformations denoted by $\Pi(L_0, L_1; \ell_0)$.
    \item We define the action functional $\mathcal A:\widetilde\Omega(L_0, L_1; \ell_0)\rightarrow \R$ by $[\ell, w]\mapsto \int w^*\omega$. We denote its set of critical points by $\operatorname{Crit}(\mathcal A)$.
\end{enumerate}
\end{defn}
We note that (c.f. Section $2.2$ of \cite{Fukaya2009}), $\operatorname{Crit}(\mathcal A)$ is of the form $[p, w]$, where $p$ denotes the constant ($[0,1]$-parametrized) path at $p\in L_0\cap L_1$ and $w:[0,1]\times [0,1]\rightarrow X$ is a homotopy of paths in $\Omega(L_0, L_1; \ell_0)$ joining $\ell_0$ to the constant path $p$.

Denote by $E, \mu$ the induced group morphisms from $I_\omega, I_\mu$ respectively, and consider $E\times\mu:\Pi(L_0, L_1; \ell_0)\rightarrow \R\times \Z$. By construction of $\widetilde\Omega(L_0, L_1; \ell_0)$, it follows that $E\times\mu$ is an injective group morphism and hence $\Pi(L_0, L_1; \ell_0)$ is an abelian group. Let $p,q\in L_0\cap L_1$ be two distinct intersection points and denote by $\Pi(L_0, L_1; p,q)$ the set of equivalence classes of elements of $\bigoplus_{\ell_0\in\pi_0(\Omega(L_0, L_1))}\Pi(L_0, L_1;\ell_0)$ such that $\alpha\sim\beta$ if $E(\alpha)\equiv E(\beta)$ and $\mu(\alpha)\equiv\mu(\beta)$.

\begin{defn}[Floer's Complex over the Novikov Field]
    We set:
    \begin{enumerate}
        \item $CF(L_0, L_1; \ell_0):=\widehat{\bigoplus}_{\operatorname{Crit}(\mathcal A)}\Z[p,w]$.
        \item $CF(L_0, L_1; \Z):=\widehat{\bigoplus}_{\ell_0\in\pi_0(\Omega(L_0, L_1))}CF(L_0, L_1; \ell_0)$.
        \item $CF(L_0, L_1; \Lambda):=CF(L_0, L_1)\widehat{\otimes}_\Z\Lambda$.
    \end{enumerate}
\end{defn}

Before we define our Floer's complex over the Novikov ring, we introduce the following equivalence relation on $CF(L_0, L_1; \Lambda)$:\\
We say $T^\lambda e^\mu[p, w]$ is equivalent to $T^{\lambda'}e^{\mu'}[q,v]$ and write $T^\lambda e^\mu[p, w]\sim T^{\lambda'}e^{\mu'}[q, v]$ if
\begin{enumerate}
    \item $p\equiv p'$,
    \item $\lambda+\int w^*\omega\equiv \lambda'+\int v^*\omega$ and
    \item $2\mu+\mu([p, w])\equiv 2\mu'+\mu([q, v])$.
\end{enumerate}

\begin{defn}($F$-filtration)
    We define the $F$-filtration on $CF(L_0, L_1; \Lambda)$ as follows:
    
    For $[p,w]\in \operatorname{Crit}(\mathcal A)$,  $T^\lambda e^\mu[p, w]\in F^{\lambda'}CF(L_0, L_1; \Lambda)$ if $\lambda+\int w^*\omega\geq \lambda'$.
\end{defn}

\begin{defn}[Floer's Complex over the Novikov Ring]\label{FloerComplexoverRing}
    We define the Floer's complex over the Novikov ring $\Lambda_0$ by $CF(L_0, L_1;\Lambda_0):=F^0CF(L_0,L_1;\Lambda)/\sim$.
\end{defn}

Now for $i=0,1$ let $\Pi(L_i):=\pi_2(X, L_i)/\sim$ where $\alpha\sim\beta$ if $\omega(\alpha)\equiv\omega(\beta)$ and $\mu(\alpha)\equiv\mu(\beta)$. Let $\mathfrak{G}(L_i)$ be the submonoid of $\R\times 2\Z$ generated by pairs $(\omega(\beta), \mu(\beta))$ where $\beta\in \Pi(L_i)$ and can be represented by $J$-holomorphic disk bounding $L_i$. Let $\mathfrak{G}_0$ be the monoid generated by $(\omega(\alpha), \mu(\alpha))$ such that $\alpha\in \Pi(L_0, L_1; p,q)$ and $p,q\in L_0\cap L_1$.

\begin{defn}\label{AinftyBimoduleGapped}
    We define $\mathfrak{G}:=\{\beta_1+\alpha+\beta_0:\beta_i\in\mathfrak G(L_i)\text{ for }i=0,1\text{ and }\alpha\in\mathfrak{G}_0\}$.
\end{defn}

\begin{theorem}[$A_\infty$-bimodule Structure]\label{AinftyBimoduleThm}
    For a pair of transversaly intersecting, closed, connected, relatively-spin Lagrangians $L_0,L_1$, the complex $CF(L_0,L_1;\Lambda_0)$ as given in Definition \ref{FloerComplexoverRing}, admits the structure of a curved, filtered, $\mathfrak{G}$-gapped $A_\infty$-bimodule over $(\mathrm{CM}^*(f_0;\Lambda_0),\{\m_k^0\}_{k\geq1}), (\mathrm{CM}^*(f_1;\Lambda_0),\{\m_k^1\}_{k\geq1})$.
\end{theorem}

\subsubsection{Preliminaries}

Our geometric realization of an $A_{n,K}$-bimodule structure on $CF(L_0,L_1;\Lambda_0)$ will be given using the Morse model as in \cite{charest2015floer}, \cite{venugopalan2020fukaya} which will be modeled over $\textit{two-paged ribbon trees with}$ $\textit{metric on its branches}$ as in \cite{ma2010geometric}, \cite{ma2015quilted} and Section $4$ of \cite{charest2015floer}. We start by recalling the necessary notions needed in this paper.

\begin{defn}[Stable Two-paged Ribbon Tree, Definition $3.1$ of \cite{ma2015quilted}]\label{two-pagedtreeDefn}
    Let $\Gamma$ be a ribbon tree such that:
    \begin{enumerate}
        \item The set of vertices of $T$ decomposes into $V(\Gamma)=V_{sp}(\Gamma)\sqcup V_0(\Gamma)\sqcup V_1(\Gamma)$.
        \item The set of internal edges of $\Gamma$ decomposes into $E_{int}(\Gamma)=E_{sp}(\Gamma)\sqcup E_0(\Gamma)\sqcup E_1(\Gamma)$.
        \item The exterior edges of $\Gamma$ are of the form:
        \begin{enumerate}
            \item $e_-,e_+$ where $e_+$ is the root of $\Gamma$,
            \item $a_1,\dots, a_{k_0}$ and 
            \item $b_1,\dots, b_{k_1}$.
        \end{enumerate}
    \end{enumerate}
    We require the following conditions:
        \begin{enumerate}
            \item $\textbf{The Spine}$ $(V_{sp},E_{sp}\cup\{e_-,e_+\})$ is a path.
            \item After ordering the elements of $V_{sp}$ according to the spine. That is, according to path joining $e_-$ to $e_+$ in $\Gamma$, the $i^{th}$-branch given by the union of the spine and all the edges connected to the $i^{th}$-vertex of $V_{sp}$, is a ribbon tree.
            \item We have no edges connected elements of $V_0$ and $V_1$.
            \item The valency of each vertex is at least $3$.
        \end{enumerate}
        In this case, we call $\Gamma$ a two-paged ribbon tree of type $(k_0,k_1)$.
\end{defn}

\begin{defn}\label{isoTwo-page}
    Let $\Gamma, \Gamma'$ be two two-paged ribbon trees as above. We say $\Gamma$ and $\Gamma'$ are isomorphic if there exists a map $F:\Gamma\rightarrow \Gamma'$ defining a bijection on the vertices, internal edges and external edges (leaves) between $\Gamma$ and $\Gamma'$, such that:
    \begin{enumerate}
        \item If an edge of $\Gamma$ joins $v$ to $w$ then, its image is an edge in $\Gamma'$ joining $F(v)$ to $F(w)$.
        \item If $i\in\{sp, 0,1\}$ then $F(V_i(\Gamma))=V_i(\Gamma')$ and $F(E_i(\Gamma))=E_i(\Gamma')$.
        \item $F(e_\pm)=e_\pm'$, $F(a_i)=a_i'$ for all $i=1,\dots, k_0$ and $F(b_j)=b_j'$ for all $j=1,\dots, k_1$.
        \item $F$ preserves the ribbon structure on each branch.
    \end{enumerate}
\end{defn}

We note that, by Theorem $1.1$ and Proposition $3.3$ of \cite{ma2015quilted}, it follows that $T_{k_0,k_1}$ the set of all isomorphism classes of two-paged ribbon trees has a topology making it homeomorphic to the Stasheff's associahedron as in \cite{stasheff1963homotopy}.

Following Section $4$ of \cite{charest2015floer}, we introduce the following Definitions:

\begin{defn}[Metric Two-paged Ribbon Tree]
    A metric two-paged ribbon tree of type $(k_0,k_1)$ is a pair $(\Gamma,l)$ such that:
    \begin{enumerate}
        \item $\Gamma$ is a two-paged ribbon tree of type $(k_0,k_1)$ as in Definition \ref{two-pagedtreeDefn}, together with an extra decoration $b(e)\geq0$ integer for each $e\in E_0(\Gamma)\sqcup E_1(\Gamma)$ corresponding to the number of $\textit{breakings}$ at edge $e$.
        \item $l:E_0(\Gamma)\sqcup E_1(\Gamma)\rightarrow[0,\infty]$ is a function such that, if $e$ is an internal edge and $l(e)=\infty$ then, $b(e)\geq 1$.
    \end{enumerate}
\end{defn}

\begin{defn}
    Let $(\Gamma,l), (\Gamma',l')$ be two metric two-paged ribbon trees. We say $(\Gamma,l)$ and $(\Gamma',l')$ are isomorphic, if there exists a map $F:(\Gamma,l)\rightarrow(\Gamma',l')$ such that:
    \begin{enumerate}
        \item $F:\Gamma\rightarrow \Gamma'$ is an isomorphism of two-paged ribbon trees as above, with the extra condition that $b(e)\equiv b(F(e))$.
        \item For any $e\in E_0(\Gamma)\sqcup E_1(\Gamma)$, $l'(F(e))\equiv l(e)$.
    \end{enumerate}
\end{defn}

Now let $(\Gamma,l)$ be a metric two-paged ribbon tree and suppose that there exists an edge $e$ of $\Gamma$ with $l(e)\equiv0$. We identify $(\Gamma,l)\sim(\Gamma/e,l_e)$, where $\Gamma/e$ is obtained from $\Gamma$ after $\textit{contracting}$ edge $e$ and setting $l_e\equiv l$ on all edges other than $e$. Continuing in this fashion for all edges $e$ of $l(e)\equiv0$, defines a $\textit{reduced broken}$ metric two-paged ribbon tree.

On the other hand, if $e$ is a non-spine interior edge of $T$ such that $l(e)=\infty$ and hence $b(e)\geq1$, then we $\textit{cut}$ $e$ into $[0,\infty)\sqcup\R\sqcup\dots\sqcup\R\sqcup(-\infty,0]$ where the number of the $\R$-factors is $b(e)-1$. This procedure will provide us with a metric two-paged ribbon tree with a spine as that of $\Gamma$ together with $b(e)$ metric rooted ribbon trees labeled by the corresponding page of $e$. In particular, we will have a surjective map 
\[
q_\Gamma:\prod_{v\in V_{sp}}T_{k_v,k_v'}\times\prod_{v\notin V_{sp}}\mathrm{Gr}_{k_v+1}\times[0,\infty]^\#\rightarrow\overline{T}_{k_0,k_1}
\]
where $\#$ denotes the number of edges having $b(e)\geq1$, for $v\in V_{sp}$ $k_v, k_v'$ denotes the number of incidences at $v$ from the $0$ and $1$ page respectively, and for $v\notin V_{sp}$, $k_v$ denotes the valency of $v$.

\begin{defn}[Moduli Space of Metric Two-paged Ribbon Trees]
    Let $\overline{T}_{k_0,k_1}$ be the set of all equivalence classes of all $\textit{reduced broken}$ metric two-paged ribbon trees. We give $\overline{T}_{k_0,k_1}$ the biggest topology so that $q_\Gamma$ is a continuous map as $\Gamma$ varies in the moduli space in $T_{k_0,k_1}$.
\end{defn}

\begin{prop}[\cite{ma2010geometric}, \cite{ma2015quilted}, \cite{charest2015floer}]\label{SpineProp}
    The moduli space $\overline{T}_{k_0,k_1}$ admits the structure of a compact smooth manifold with corners such that:
    \begin{enumerate}
        \item $\overline{T}_{k_0,k_1}$ is homeomorphic to the Stasheff's associahedra. 
        \item $\dim_\R\overline{T}_{k_0, k_1}=k_0+k_1-1$.
        \item The codimension-$1$ corner of $\overline{T}_{k_0,k_1}$ is of the union of the following factors:
        \begin{enumerate}
            \item $\overline{T}_{k_0-k+1, k_1}\times\overline{\mathrm{Gr}}_{k+1}$ for $2\leq k\leq k_0$,
            \item $\overline{T}_{k_0,k_1-k+1}\times\overline{\mathrm{Gr}}_{k+1}$ for $2\leq k\leq k_1$,
            \item $\overline{T}_{r_0,r_1}\times \overline{T}_{s_0,s_1}$ where $r_0+s_0=k_0$ and $r_1+s_1=k_1$.
        \end{enumerate}
        \item $\overline{T}_{k_0,k_1}$ has a stratification by smooth manifolds indexed by combinatorial types. That is, let $[\Gamma]$ be an isomorphism class of an element in $T_{k_0,k_1}$, then $\overline{T}_{k_0,k_1}=\bigsqcup_\Gamma T_\Gamma$ where $T_\Gamma\cong \prod_{v\in V_{sp}(\Gamma)}T_{k_v,k_v'}\times\prod_{v\notin V_{sp}(\Gamma)}\mathrm{Gr}_{k_v+1}$.
    \end{enumerate}
\end{prop}

Let $[\Gamma,l]\in\overline{T}_{k_0,k_1}$ be an isomorphism class of metric two-paged ribbon tree. The metric $l$ on $\Gamma$ will associate coordinates on each non-spine edge of $\Gamma$ as follows:

\begin{enumerate}
    \item If $e$ is an internal edge such that $0<l(e)<\infty$ then, we identify $e\cong I_e\equiv[0,l(e)]_s$. 
    \item If $e$ is an internal edge with $l(e)=\infty$ and hence $b(e)\geq1$, we identify $e\cong I_e\equiv[0,\infty)_s\sqcup\R_s\sqcup\dots\sqcup\R_s\sqcup(-\infty,0]_s$ where the number of $\R$-factors is $b(e)-1$.
    \item If $e$ is an external edge of $T$ other then $e_-, e_+$, that is a non-spine leaf, then we identify $e\cong I_e\equiv \R_s\sqcup\dots\sqcup\R_s\sqcup(-\infty,0]_s$ where the number of $\R$-factors is $b(e)$.
\end{enumerate}

Now let $L_0, L_1$ be a pair of transversely intersecting, closed, connected, relatively-spin as a pair Lagrangians in $X$. For $i=0,1$, let $(f_i,g_i)$ be a Morse-Smale pair on $L_i$. Let $\beta_v\in\pi_2(X,\mathbf{L})$ be a finite sequence of relative spherical classes indexed by the internal vertices of $\Gamma$, $x_1,\dots,x_{k_0}\in\mathrm{Crit}(f_0)$, $y_1,\dots,y_{k_1}\in\mathrm{Crit}(f_1)$ and $p,q\in L_0\cap L_1$. For the edges $e_\pm$ of $\Gamma$, we associate $q,p$ respectively. While for the rest of spine-edges and using the ordering given by the path joining $e_-$ to $e_+$ we associate $p_1,\dots, p_r\in L_0\cap L_1$. On the other hand, for each external vertex of the $0$-page of $T$, we associate $x_i$ respecting the ribbon structure. Moreover, for each external vertex of the $1$-page, we associate $y_i$ respecting the ribbon structure on $T$. With this in mind, we give the following Defintion:

\begin{defn}[Floer-Morse Trajectories of Fixed Type]\label{FloerMorseStripFixedCombDefn}
    Let $\mathcal{M}(\Gamma;\underline{x},p,q,\underline{y};(\beta_v)):=\{((\Gamma, l),u):(\Gamma,l)\in T_\Gamma, u:\Gamma\rightarrow X\}$ be the moduli space of all continuous maps $u$ given by:
    \begin{enumerate}
        \item For each $v$ interior vertex of the spine, $u_v\in\M_{|v|+2}(L_0,L_1, p_i,p_{i+1}; J,\beta_v)$ with the understanding that $|v|$ denotes the valency of $v$, and $p_0\equiv p$, $p_{r+1}\equiv q$.
        \item For each $v$ interior vertex of the $i$-page of $T$, $u_v\in\M_{|v|+1}(X,L_i;J,\beta_v)$.
        \item For each $e$ edge of the $i$-page of $\Gamma$, we associate a $\textit{domain-dependent}$ Morse trajectory of $(f_i,g_i)$. Namely, $u_e:I_e\rightarrow L_i$ smooth map such that:
        \begin{equation*}
            \frac{d}{ds}u_e(s)+\nabla^{g_i}\tilde f_i(s,u_e(s))=0,
        \end{equation*}
        where $\tilde f_i:\widetilde{\overline{\mathrm{Gr}}}_{k+1}\times L_i\rightarrow\R$ where $\widetilde{\overline{\mathrm{Gr}}}_{k+1}\rightarrow\overline{\mathrm{Gr}}_{k+1}$ is the universal curve associated to the compactified metric rooted ribbon trees, as in Section $4$ of \cite{rabah2024fukaya}, of $(k+1)$-leaves where $k$ is given by the number of external leaves in the $i$-page of the associated branch, as in item $(2)$ of Definition \ref{two-pagedtreeDefn}, other than the spine ones. Moreover, if $e$ is an external leaf on the $i$-page and under the identification $e\cong I_e\equiv(-\infty,0]\sqcup\R\sqcup\dots\sqcup\R$ given by $l$, we require $u_e(-\infty)=x_j, y_j$ according to which page $e$ lies on.
    \end{enumerate}
    We give $\mathcal{M}(\Gamma;\underline{x},p,q,\underline{y};(\beta_v))$ the $C^\infty_{loc}$-topology.
\end{defn}

\begin{lemma}\label{FloerMorseStripFixLem}
    For a generic choices of domain-dependent Morse functions $\tilde f_i:\widetilde{\overline{\mathrm{Gr}}}_{k_i+1}\times L_i\rightarrow\R$, the moduli space $\mathcal{M}(\Gamma;\underline{x},p,q,\underline{y};(\beta_v))$ admits a Global Kuranishi chart representation.
\end{lemma}

\begin{proof}
    We follow the same argument as in the proof of Proposition $4.37$ of \cite{rabah2024fukaya}. Using the cyclic ordering on each branch of $\Gamma$ and after taking the minimal path joining $e_+$ to any interior vertex of $\Gamma$, we can order all interior vertices of $\Gamma$. Say the interior vertices on the $0$-page are labeled by $a_1,\dots,a_r$, the interior vertices on the $1$-page are labeled by $b_1,\dots, b_s$ and the vertices on the spine are labeled by $v_1,\dots, v_{r-1}$. Such ordering will also induce an ordering on the relative spherical classes $\beta_v$. Now on each of the non-spine vertices, we apply Proposition $4.37$ of \cite{rabah2024fukaya} and get that for a generic choice of domain-dependent Morse function $\tilde f_i$ associated to $f_i:L_i\rightarrow\R$, we have a Global Kuranishi chart representation $(G_{a_j}, \T_{a_j},\E_{a_j},s_{a_j})$ for $j=1,\dots, r$ and $(G_{b_j}, \T_{b_j},\E_{b_j},s_{b_j})$ for $j=1,\dots,s$ of the corresponding Floer-Morse trajectories on each branch without the spine of $\Gamma$. On other hand, for each interior vertex of the spine of $\Gamma$, we adhere to Theorem \ref{GKCStrips} to get a Global Kuranishi chart $(G_{v_j}, \T_{v_j}, \E_{v_j}, s_{v_j})$ for $j=1,\dots, r-1$. Now using the submersivity property and same argument as in the proof of Proposition $4.37$ of \cite{rabah2024fukaya}, it follows that $\T:=\prod(\T_{a_j}\times\T_{b_j}\times\T_{v_j})\times_{ev}(L_0^{k_0}\times L_1^{k_1})$ is a Thickening of $\mathcal{M}(\Gamma;\underline{x},p,q,\underline{y};(\beta_v))$ having a group of symmetry $G:=\prod (G_{a_j}\times G_{b_j}\times G_{v_j})$, an obstruction bundle given by $\E:=\iota^*(\bigoplus\E_{a_j}\boxplus\E_{b_j}\boxplus\E_{v_j})$ where $\iota:\T\hookrightarrow\prod(\T_{a_j}\times\T_{b_j}\times\T_{v_j})$, and Kuranishi section $s:=\iota^*(\bigoplus s_{a_j}\boxplus s_{b_j}\boxplus s_{v_j})$.
\end{proof}

\begin{cor}\label{FloerMorseStripFixedCor}
    For $\beta\in\pi_2(X,\mathbf{L})$ and for generic choice of domain-dependent Morse functions $\tilde f_i$ associated to $f_i$, the moduli space $\mathcal M(\Gamma;\underline{x},p,q,\underline{y};\beta):=\bigsqcup_{\sum\beta_v=\beta}\mathcal{M}(\Gamma;\underline{x},p,q,\underline{y};(\beta_v))$ admits a Global Kuranishi chart representation. 
\end{cor}

Now let $\Gamma'$ be obtained from $\Gamma$ after collapsing an interior edge. That is, $\Gamma'\cong\Gamma/e$ isomorphic as two-paged ribbon trees as in Definition \ref{isoTwo-page}. By Proposition \ref{SpineProp}, we choose and fix a metric on $\overline{T}_{k_0,k_1}$ and for $\epsilon>0$ small we denote by $\mathbb D_\epsilon$ an $\epsilon$-neighborhood of $T_{\Gamma'}\subseteq \overline{T}_\Gamma$.

\begin{defn}
    For $\beta\in\pi_2(X, \mathbf{L})$ and $\tilde f_i$ domain-dependent Morse functions associated to $f_i$, we denote by $\mathcal M^\epsilon(\Gamma;\underline{x},p,q,\underline{y};\beta)=\{((\Gamma,l),u)\}\subseteq\mathcal{M}(\Gamma;\underline{x},p,q,\underline{y};\beta)$ where $(\Gamma,l)\in\mathbb D_\epsilon$.
\end{defn}

\begin{lemma}\label{FloerMorseStripNbhd}
    If $\epsilon>0$ is small enough then, for the same choice of domain-dependent Morse functions $\tilde f_i$ associated to $f_i$, as in Corollary \ref{FloerMorseStripFixedCor} or Lemma \ref{FloerMorseStripFixLem}, the moduli space $\mathcal M^\epsilon(\Gamma;\underline{x},p,q,\underline{y};\beta)$ admits a Global Kuranishi chart representation.
\end{lemma}

\begin{proof}
    We first apply Corollary \ref{FloerMorseStripFixedCor} and fix a choice of domain-dependent Morse functions $\tilde f_i$ associated to $f_i$ to get a Global Kuranishi chart representation of the moduli space $\mathcal M(\Gamma';\underline{x},p,q,\underline{y};\beta)$. Now we consider the following two cases:\\
    $\textbf{Case 1:}$ The collapsed edge is an interior edge on the Spine of $\Gamma$.

    In this case, Items $(4)$ and $(5)$ of the compatibility conditions of Theorem \ref{GKCStrips} allow us to apply the same argument of the proof of Lemma \ref{FloerMorseStripFixLem} such choice of $\tilde f_i$, to get a Global Kuranishi chart representation for $\mathcal{M}^\epsilon(\Gamma;\underline{x},p,q,\underline{y};\beta)$, for any $\epsilon>0$.\\
    $\textbf{Case 2:}$ The collapsed edge is a non-spine interior edge of $\Gamma$.

    In this case, Items $(4)$, $(5)$ and $(6)$ of the compatibility conditions of Theorem \ref{GKCStrips}, allow us to use Proposition $4.42$ of \cite{rabah2024fukaya} for $\epsilon>0$ small enough. Noting that, by definition of domain-dependent Morse functions, we have $\tilde f_i\equiv f$ on a small enough neighborhood of $\mathrm{Crit}(f_i)$, it follows, possibly after taking yet a smaller $\epsilon>0$, that we can use the same construction as given in the proof of Lemma \ref{FloerMorseStripFixLem} to give a Global Kuranishi chart representation of $\mathcal M^\epsilon(\Gamma;\underline{x},p,q,\underline{y};\beta)$ using the fixed choice $\tilde f_i$.
\end{proof}

\begin{defn}[Floer-Morse Trajectories]
    For $\beta\in\pi_2(X,\mathbf{L})$, $x_,\dots,x_{k_0}\in\mathrm{Crit}(f_0)$, $y_1,\dots,y_{k_1}\in\mathrm{Crit}(f_1)$, $p,q\in L_0\cap L_1$, and a choice of domain-dependent Morse functions $\tilde f_i$ associated to $f_i$, we define $\M(\underline{x}, p,q,\underline{y};\beta):=\{((\Gamma,l),u)\}$ such that:
    \begin{enumerate}
        \item $(\Gamma,l)\in \overline{T}_{k_0,k_1}$ and
        \item $u\in\mathcal{M}(\Gamma;\underline{x},p,q,\underline{y};\beta)$.
    \end{enumerate}
\end{defn}

By Gromov's Compactness, Morse trajectory compactness and the gluing results of \cite{charest2015floer}, it follows that the codimension-$1$ (normalized) corner of $\M(\underline{x};p,q;\underline{y};\beta)$ is of the form
\[
\partial\M(\underline{x};p,q;\underline{y};\beta)=\M_{sp}\sqcup\M_0\sqcup\M_1,
\]
such that:
\[
\M_{sp}=\bigsqcup\M(\underline{x}^-;p,r;\underline{y}^-;\beta_1)\times\M(\underline{x}^+;r,q;\underline{y}^+;\beta_2)
\]
where the union is taken over $r\in L_0\cap L_1$, $\beta_1+\beta_2=\beta$, and all decompositions $\underline{x}=\underline{x}^-\cup\underline{x}^+$ and $\underline{y}=\underline{y}^-\cup\underline{y}^+$.
\[
\M_0=\bigsqcup\M(z,x_{t+1},\dots,x_{t+s};\beta_1)\times\M(x_1,\dots,x_t,z,x_{t+s+1},\dots,x_{k_0};p,q;\underline{y};\beta_2)
\]
where the union is taken over $z\in\mathrm{Crit}(f_0)$, $\beta_1+\beta_2=\beta$, $0\leq s\leq k_0$ and $0\leq t\leq k_0-s$.
\[
\M_1=\bigsqcup\M(\underline{x};p,q;y_1,\dots,y_{t},z,y_{t+s+1},\dots,y_{k_1};\beta_1)\times \M(z,y_{t+1},\dots,y_{t+s};\beta_2)
\]
where the union is taken over $z\in\mathrm{Crit}(f_1)$, $\beta_1+\beta_2=\beta$, $0\leq s\leq k_1$ and $0\leq t\leq k_1-s$.

\begin{theorem}\label{GKCStripMorseFloerThm}
    For generic choices of domain-dependent Morse functions $\tilde f_i:\widetilde{\overline{\mathrm{Gr}}}_{k_i+1}\times L_i\rightarrow\R$ associated to $f_i$, we have the following:

    \begin{enumerate}
        \item The moduli space $\M(\underline{x}, p,q,\underline{y};\beta)$ admits a normally complex Global Kuranishi chart representation $(G_\beta,\T(\underline{x},p,q,\underline{y},\beta),\E(\underline{x},p,q,\underline{y},\beta),s_{\underline{x},p,q,\underline{y},\beta})$.
        \item Each of factors of the codimension-$1$ corner of $\M(\underline{x}, p,q,\underline{y};\beta)$, admits a normally complex Global Kuranishi chart representation $(G_j, \T_j,\E_j,s_j)$ such that:
        \begin{enumerate}
            \item $G_1\times G_2\hookrightarrow G$.
            \item The following diagram commutes:
        \[\begin{tikzcd}
	{\mathcal E_1\boxplus\mathcal E_2} && {\mathcal E(\underline x;p,q;\underline y;\beta)} \\
	\\
	{\mathcal T_1\times \mathcal T_2} && {\mathcal T(\underline x;p,q;\underline y;\beta)} \\
	\\
	{\mathcal B_1\times \mathcal B_2} && {\overline T_{k_0,k_1}}
	\arrow[hook, from=1-1, to=1-3]
	\arrow[from=1-1, to=3-1]
	\arrow[from=1-3, to=3-3]
	\arrow["{s_1\boxplus s_2}", bend left=30pt, , from=3-1, to=1-1]
	\arrow[hook, from=3-1, to=3-3]
	\arrow[from=3-1, to=5-1]
	\arrow["{s_{\underline x,p,q,\underline y,\beta}}"', bend right=30pt, , from=3-3, to=1-3]
	\arrow[from=3-3, to=5-3]
	\arrow["{\text{glue}}", from=5-1, to=5-3]
\end{tikzcd}\]       
  where $\mathcal B_1\times\mathcal B_2$ is either:
  \[
  \overline{T}_{r_0,r_1}\times \overline{T}_{s_0,s_1}\text{ where } r_0+r_1=k_0, s_0+s_1=k_1
  \]
            \[
            \overline{\mathrm{Gr}}_{s+1}\times \overline{T}_{k_0-s+1,k_1}\text{ where }0\leq s\leq k_0
            \]
            or
            \[
            \overline{T}_{k_0, k_1-s+1}\times\overline{\mathrm{Gr}}_{s+1}\text{ where }0\leq s\leq k_1.
            \]
            Moreover, the above diagram is compatible with the normal complex structure and orientation.
            \item If $(G_j,\T_j,\E_j,s_j)$ is a Global Kuranishi chart representation of a factor of the codimension-$1$ corner of $\M(\underline{x};p,q;\underline{y};\beta)$ representing a disk component, that is either the first factor of $\M_0$ or the second factor of $\M_1$ above, then it is given by Theorem $4.50$ of \cite{rabah2024fukaya}.
\end{enumerate}
    \end{enumerate}
\end{theorem}

\begin{proof}
    The proof follows after using the same arguments in Section $4$ of \cite{rabah2024fukaya}.

    $\textbf{Step 1:}$ For a given $\beta\in\pi_2(X,\mathbf{L})$, $x_1,\dots,x_{k_0}\in\mathrm{Crit}(f_0)$, $p,q\in L_0\cap L_1$, $y_1,\dots, y_{k_1}$ and for a generic choice of domain-dependent Morse functions $\tilde f_i$ associated to $f_i$. There exists a Global Kuranishi chart representation $(G,\T,\E,s)$ of the moduli space $\mathcal{M}:=\{((\Gamma, l),u)\}$ where $(\Gamma, l)$ is a reduced, unbroken two-paged metric ribbon tree and $u\in\mathcal{M}((\Gamma,l),\underline{x},p,q,\underline{y};\beta)$, together with a smooth submersion $\T\rightarrow T_{k_0,k_1}$.
    \begin{proof}[Proof of Step 1]
        Indeed, using Lemma \ref{FloerMorseStripFixLem} or Corollary \ref{FloerMorseStripFixedCor} and for a fixed $(\Gamma,l)$, it follows that, for a generic choice of domain-dependent Morse functions $\tilde f_i$, there exists a Thickening $\T(\Gamma,l)$ of the moduli space $\mathcal{M}((\Gamma,l),\underline{x},p,q,\underline{y};\beta)$. Now note that, for fixed $k_0$ and $k_1$, it follows by Euler characteristic considerations, that the number of combinatorial types allowed, that is the number of interior vertices, is finite. Thus, applying Corollary \ref{FloerMorseStripFixedCor} and Lemma \ref{FloerMorseStripNbhd} implies that if we set $\T:=\{((\Gamma,l),u):(\Gamma,l)\in T_{k_0,k_1}. u\in\T(\Gamma,l)\}$ then by the parametric implicit function theorem, the collection $\{U\times V_\epsilon\}$ where $u\in U\subseteq\T(\Gamma,l)$ is a coordinate chart around $u\in\T(\Gamma,l)$ and $V_\epsilon$ is an $\epsilon$-neighborhood of the zero section of the normal bundle of $T_\Gamma\subseteq T_{k_0,k_1}$, form an atlas on $\T$. Moreover, from Lemma \ref{FloerMorseStripNbhd}, it follows that locally, $\T\rightarrow T_{k_0,k_1}$ is a projection map and hence a smooth submersion.

        Now let $d:=\Omega(\beta)$ where $\Omega$ is as in Lemma \ref{Existinceof2-form} and denote by $G:=O(d)$. Fix $(\Gamma,l)\in T_{k_0,k_1}$ and $u\in\T(\Gamma,l)$. The $G$-action on $((\Gamma,l),u)$ is given by the action of $\prod_{v\in V_{int}(\Gamma)}G_v$ and the trivial action on $(\Gamma, l)$. Thus, possibly after shrinking $U$, it follows that the collection $\{U\times V_\epsilon\}$ forms an atlas on $\T$ making it a smooth $G$-manifold. Notice that by the proof of Lemma \ref{FloerMorseStripFixLem} and Lemma \ref{compatibleObsBundle}, it follows that, if we define $\E\rightarrow\T$ to be the pullback of the trivial bundle $T_{k_0,k_1}\times\E_\beta\rightarrow T_{k_0,k_1}$, where $\E_\beta$ is as in Theorem \ref{GKCStrips} then for $s:\T\rightarrow\E$ given as in Definition $4.46$ of \cite{rabah2024fukaya}, the tuple $(G,\T,\E,s)$ forms a Global Kuranishi chart representation of $\mathcal{M}$, as desired.
    \end{proof}

    $\textbf{Step 2:}$ The compactification $\M$ of $\mathcal{M}$ as given above, where $\mathcal M$ is as in $\textbf{Step 1}$ admits a Global Kuranishi chart representation $(G,\T,\E,s)$ together with a smooth submersion $\T\rightarrow \overline{T}_{k_0,k_1}$.

    \begin{proof}[Proof of Step 2]
        We argue by induction on the number of breakings. Set $\T:=\{((\Gamma,l),u)\}$, where $(\Gamma,l)$ is a reduced broken two-paged metric ribbon tree and if we write $\Gamma\cong \Gamma_1\sqcup\dots\sqcup\Gamma_r/\sim$ as a $\textit{gluing}$ of unbroken reduced two-paged metric ribbon trees $\Gamma_i$ then, $u=u_1\vee\dots\vee u_r$. In particular, $u_i\in \T(\Gamma, l_i)$ for all $i=1,\dots, r$. We give $\T$ the Gromov-Hausdorff topology. For simplicity, assume that $\Gamma\cong\Gamma_1\sqcup\Gamma_2/\sim$ is once broken and write $u=u_1\vee u_2$. For $i=1,2$, let $u_i\in U_i\subseteq\T(\Gamma_i,l_i)$ be a coordinate chart around $u_i\in\T(\Gamma_i,l_i)$ after fixing domain-dependent Morse functions $\tilde f_0, \tilde f_1$ as in $\textbf{Step 1}$. Furthermore, assume that the breaking occurs at a non-spine edge, say on the $0$-page. Let $z\in\mathrm{Crit}(f_0)$ at which the breaking occurs and take small enough neighborhoods $O_i\subseteq L_0$ such that on their closure $\tilde f_0\equiv f_0$. On such neighborhoods, let $\gamma_1:[s_1,\infty)\rightarrow L_0$ and $\gamma_2:(-\infty,s_2]\rightarrow L_0$ be Morse trajectories of $f_0$ asymptotic to $z$. By the Gluing results as found in \cite{schwarz1993morse} and for $R>>1$, we get a smooth embedding
        \[
        \mathcal{M}(\gamma_1(s_1),z)\times\mathcal M(z,\gamma_2(s_2))\times [R,\infty)\hookrightarrow \mathcal{M}(\gamma_1(s_1),\gamma_2(s_2))
        \]
        on the corresponding Morse trajectories of $u_1$ and $u_2$ on such broken edge. Set $U:=U_1\times U_2\times[R,\infty)$. Now consider the case when the breaking occurs at an internal spine edge. Let $\epsilon>0$ be as in Lemma \ref{FloerMorseStripNbhd} and consider $V_\epsilon$ an $\epsilon$-neighborhood of the $T_{\Gamma_1}\times T_{\Gamma_2}\hookrightarrow T_\Gamma$. For $i=1,2$, let $u_i\in U_i\subseteq \T(\Gamma_i,l_i)$ be a smooth coordinate chart and set $U:=U_1\times U_2\times V_\epsilon$.
        
        Note that, as $L_0, L_1$ intersect transversely, $\mathrm{Crit}(f_i)$ is a finite set and the number of $\beta_i\in\pi_2(X,\mathbf{L})$ such that $\omega(\beta_i)\leq \omega(\beta)$ is also finite, we can repeat the above construction finitely many times to get for each $((\Gamma,l),u)\in\T$ a coordinate chart $U$ around it. By the proof Lemma \ref{FloerMorseStripFixLem}, Lemma \ref{FloerMorseStripNbhd}, the submersive property of the evaluation map and items $(2), (4)$ and $(5)$ of the $\textbf{compatibility conditions}$ of Theorem \ref{GKCStrips}, it follows by the parametric implicit function theorem, that such collection of charts forms an atlas on $\T$. Now the claim follows from the last paragraph of the proof of $\textbf{Step 1}$.
    \end{proof}

    To finish the proof of the Theorem, we argue by upward induction of the energy $\omega(\beta)$ and the number of leaves on each page, that is $k_0$ and $k_1$.

    The base case follows after generically choosing $\tilde f_i$ using Theorem $4.31$ of \cite{rabah2024fukaya}. As for the inductive step, suppose that we have $(G_j,\T_j,\E_j,s_j)$ Global Kuranishi charts representing a connected component of the codimension-$1$ corner of $\M$ with the understanding that $(G_j,\T_j,\E_j,s_j)$ is constructed inductively over the number of breakings usings $\textbf{Step 1}$ and $\textbf{Step 2}$. Let $(G,\T,\E,s)$ be the Global Kuranishi chart representing $\M$ as constructed in $\textbf{Step 2}$. By construction of $\T$ and by our inductive hypothesis, it follows that we have an embedding $\T_1\times\T_2\hookrightarrow\T$ compatible with the group embedding $G_1\times G_2\hookrightarrow G$ making the following diagram commutes:
    \[\begin{tikzcd}
	{\mathcal T_1\times\mathcal T_2} && {\mathcal T} \\
	\\
	{\mathcal B_1\times\mathcal B_2} && {\overline T_{k_0,k_1}}
	\arrow[hook, from=1-1, to=1-3]
	\arrow[from=1-1, to=3-1]
	\arrow[from=1-3, to=3-3]
	\arrow["{\text{glue}}", from=3-1, to=3-3]
\end{tikzcd}\]
with the understanding that the $G_j, G$ action is trivial on $\mathcal B_j, \mathcal{B}$ respectively. On the other hand, the commutativity of 
\[\begin{tikzcd}
	{\mathcal E_1\boxplus\mathcal E_2} && {\mathcal E} \\
	\\
	{\mathcal T_1\times\mathcal T_2} && {\mathcal T}
	\arrow[hook, from=1-1, to=1-3]
	\arrow[from=1-1, to=3-1]
	\arrow[from=1-3, to=3-3]
	\arrow["{s_1\boxplus s_2}", bend left=30, from=3-1, to=1-1]
	\arrow[hook, from=3-1, to=3-3]
	\arrow["s"', bend right=30, from=3-3, to=1-3]
\end{tikzcd}\]
follows from our inductive hypothesis and the construction of $\E$ as given in $\textbf{Step 1}$. Lastly, to show compatibility with the normal complex structure and orientations, we first give the Global Kuranishi chart construction in $\textbf{Step 1}$ a stable normal complex structure as in the proof of Theorem $4.48$ of \cite{rabah2024fukaya}. Then still under the same induction argument as above, namely upward multi-layered induction over $\omega(\beta)$, $k_0,k_1$, we use the same argument as given in the second paragraph of the proof of Theorem $4.50$ of \cite{rabah2024fukaya} to deduce the claim.
\end{proof}

\begin{proof}[Proof of Theorem \ref{AinftyBimoduleThm}]
    Having Lemma $7.2.69$ and Theorem $7.2.72$ of \cite{fukaya2010lagrangian} in mind, it suffices to construct a filtered $A_{n,K}$-bimodule structure on $CF(L_0,L_1;\Lambda_0)$, provided that such construction can be extended to an $A_{n',K'}$-bimodule structure on $CF(L_0,L_1;\Lambda_0)$ for $(n,K)\precsim (n',K')$. 

    To be more precise, let $\frak{G}$ be as given in Definition \ref{AinftyBimoduleGapped} and the partial order $\precsim$ is as given in the paragraph under Definition \ref{AnKbimoduleDefn} above. For $i=0,1$, let $(f_i,g_i)$ be a Morse-Smale pair on $L_i$ and for $k_0,k_1\geq0$ integers, let $x_1,\dots,x_{k_0}\in\mathrm{Crit}(f_0)$ and $y_1,\dots, y_{k_1}\in\mathrm{Crit}(f_1)$. For $p\in L_0\cap L_1$ and $\beta\in\pi_2(X,\mathbf{L})$, we define
    \[
    \n_{k_0,k_1,\beta}(p):=\sum_{q\in L_0\cap L_1}(-1)^*|s^{-1}_\epsilon(0)\cap \mathcal U_{free}|q,
    \]
    where the sum is taken over all $q\in L_0\cap L_1$ such that the expected dimension of $\M(\underline{x};p,q;\underline{y};\beta)$ is zero and $*$ is as in Definition $8.7.1$ of \cite{fukaya2010lagrangian} with the understanding that the unstable manifold of a given Morse critical point is a defines a singular chain after taking an apropriate tringulation in $L_i$. Moreover, here $(\mathcal U, E,\mathcal S)$ is the associated derived orbifold chart from a choice of Global Kuranishi chart $(G,\T,\E,s)$ representing $\M(\underline{x};p,q;\underline{y};\beta)$. Now set
    \[
    \n_{k_0,k_1}:=\sum_{\frak{G}}\n_{k_0,k_1,\beta}T^{\omega(\beta)}e^\frac{\mu_L(\beta)}{2}
    \]
    and extend $\Lambda_0$-linearily. To show that such operations $\{\n_{k_0,k_1}\}_{k_0,k_1\geq0}$ defines the structure of a filtered, $\frak{G}$-gapped, curved $A_{n,K}$-bimodule structure on $CF(L_0,L_1;\Lambda_0)$ over the $A_\infty$-algebras $(\mathrm{CM}^*(f_i;\Lambda_0),\{\m^i_{k}\}_{k\geq0})$ given by Corollary $5.4$ of \cite{rabah2024fukaya}, we argue by induction over $(\frak{G},\precsim)$. For the base case, we use Theorem \ref{GKCStripMorseFloerThm} to give a Global Kuranishi chart representation of $\M(\underline{x};p,q;\underline{y};\beta)$ and then perturb $s$ to an FOP-perturbation $s_\epsilon$ using Proposition $2.10$ of \cite{bai2022arnold} (c.f. Appendix B of \cite{rabah2024fukaya}). As for the inductive step, let $\M(\underline{x};p,q;\underline{y};\beta)$ be such that its expected dimension is $1$ and suppose that any component of its codimension-$1$ corners is equipped with a Global Kuranishi chart representation $(G_j,\T_j,\E_j,s_j)$ together with an FOP-perturbation $s_{j, \epsilon}$. Moreover, if $(G_j,\T_j,\E_j,s_j)$ is a Global Kuranishi chart representation of some disk components then it is equivalent to the one given by Theorem $4.50$ of \cite{rabah2024fukaya}. From the proof of Theorem \ref{GKCStripMorseFloerThm}, it follows that we can construct a Global Kuranishi chart representation $(G,\T,\E,s)$ of $\M(\underline{x};p,q;\underline{y};\beta)$ such that it satisfies the conditions of the statement of Theorem \ref{GKCStripMorseFloerThm}. In particular, we have the following commutative diagram of smooth manifolds:
    \[\begin{tikzcd}
	{\mathcal E_1\boxplus\mathcal E_2} && {\mathcal E} \\
	\\
	{\mathcal T_1\times\mathcal T_2} && {\mathcal T}
	\arrow[hook, from=1-1, to=1-3]
	\arrow[from=1-1, to=3-1]
	\arrow[from=1-3, to=3-3]
	\arrow["{s_{1}\boxplus s_2}", bend left=30, from=3-1, to=1-1]
	\arrow[hook, from=3-1, to=3-3]
	\arrow["s"', bend right=30, from=3-3, to=1-3]
\end{tikzcd}\]
which is compatible with the group embedding $G_1\times G_2\hookrightarrow G$.
By Corollary $2.14$ of \cite{bai2022arnold}, it follows that $s_{1,\epsilon}\boxplus s_{2,\epsilon}$ is a strongly transversal FOP section as in Definition $2.8$ of \cite{bai2022arnold} (c.f. Appendix B of \cite{rabah2024fukaya} for a recollection of all the necessary ingredients and definitions used in this proof). After fixing a $G$-equivariant metric on $\T$ we can find an invariant neighborhood of $\T_1\times \T_2$ in $\T$. Then applying the relative extension property of FOP sections, as in Proposition $2.14$ of \cite{bai2022arnold}, it follows that we can find a strongly transversal FOP perturbation $s_\epsilon$ of $s$.

With such construction at hand, we can us the arguments of the proof of Theorem $3.7.1$ of \cite{Fukaya2009} to complete the proof.
\end{proof}

\section{Lagrangian Floer Theory over $\Z$}
     
\subsection{Single Lagrangian}\label{SelfFloerSect}

Suppose that $L\subset X$ is a closed, connected, relatively-spin Lagrangian. We recall that the main Theorem of \cite{rabah2024fukaya}, was the construction of a filtered, curved, gapped, $A_\infty$-algebra over the Novikov ring with integers on the $\textit{Morse Model}$ of the cohomology of a $L$. To be more precise, for any Morse-Smale pair $(f,g)$, the Morse complex $\mathrm{CM}^*(f; \Lambda_0)\equiv\mathrm{CM}^*(f;\Z)\widehat{\otimes}_\Z\Lambda_0$ admits the structure of a filtered, curved, gapped $A_\infty$-algebra, which we will denote by $\{\m_k\}_{k\geq0}$. We recall the following Definitions:

\begin{defn}[Bounding Cochain]\label{boundingCochainDefn}
  An element $b\in\mathrm{CM}^1(f;\Z)\widehat{\otimes}_\Z\Lambda_0^+\equiv\mathrm{CM}^1(f;\Lambda_0^+)$ is said to be a bounding cochain for $(\mathrm{CM}^*(f;\Lambda_0),\{\m_k\}_{k\geq0})$, if $b_i$ solves the $\textit{Maurer-Cartan}$ Equation
  \[
  \sum_{k\geq0}\m_k(b,\dots,b)=0,
  \]
  where $\Lambda_0^+$ denotes the elements of $\Lambda_0$ with positive-valuations.
\end{defn}

\begin{defn}[Deformed $A_\infty$-algebra]\label{deformedAinftyDefn}
    Suppose that $b$ is a bounding cochain for $(\mathrm{CM}^*(f;\Lambda_0),\{\m_k\}_{k\geq0})$. We define the $b$-deformed $A_\infty$-algebra on $\mathrm{CM}^*(f;\Lambda_0)$ by
    \[
    \m_k^{b}(x_1, \dots,x_k):=\sum_{s_0\geq0}\dots\sum_{s_k\geq0}\m_{k+\sum s_i}(b^{\otimes s_0}, x_1, b^{\otimes s_1}, \dots, x_k, b^{\otimes s_k}),
    \]
    where $x_1,\dots, x_k\in \mathrm{CM}^*(f; \Lambda_0)$.
\end{defn}

\begin{defn}[Unobstructed over $\Lambda_0$]\label{UnobstructedDefn}
    We say $L$ is unobstructed over $\Lambda_0$ if there exists a bounding cochain $b$, as in Definition \ref{boundingCochainDefn}, for the $A_\infty$-algebra $(\mathrm{CM}^*(f;\Lambda_0),\{\m_k\}_{k\geq0})$. In this case, we will say $(L, b)$ is unobstructed over $\Lambda_0$.
\end{defn}

We note that, by the definition of $b$ being a bounding cochain and using Proposition $3.6.10$ of \cite{Fukaya2009}, it follows that $(\mathrm{CM}^*(f;\Lambda_0),\{\m_k^{b}\}_{k\geq0})$ is filtered, uncurved, gapped $A_\infty$-algebra. In particular, $\m_0^{b}(1)=0$ and hence $\m_1^{b}\circ \m_1^{b}=0$. With this in mind, we define the following homology modules:

\begin{defn}[Lagrangian Floer Theory over $\Z$]\label{LagFloerHomSingle L defn}
Suppose that $L$ is unobstructed over $\Lambda_0$, as in Definition \ref{UnobstructedDefn}, and let $\m_1^{b}:\operatorname{CM}^*(f;\Lambda_0)\rightarrow\mathrm{CM}^*(f; \Lambda_0)$ be as above. We define a Lagrangian Floer homology group
\[
HF^*((L, b);\Lambda_0):=\ker(\m_1^{b})/\mathrm{Im}(\m_1^{b}).
\]    
\end{defn}

\begin{remark}
    Even though it does not show in our notation and as for now, the definition of $HF^*((L, b); \Lambda_0)$ depends on various choices, such as the choice of:
    \begin{enumerate}
        \item Morse function on $L$.
        \item Almost complex structure on $X$.
        \item Auxiliary data used in the construction of Global Kuranishi charts.
        \item Choice of perturbations.
        \item Choice of bounding cochain.
    \end{enumerate}
    We will deal with such $\textit{independece question}$ in \cite{RabahIntegralFukaya}, while use such constructions to give geometric applications as in Section \ref{AppSection}.
\end{remark}

    In principle and since the authors argue over any arbitrary ring, we can use Sections $3.6.2$ of \cite{Fukaya2009} to give conditions for the existence of bounding cochains for a given Lagrangian as above. Indeed, we have the following result whose proof is a matter of notational change from the one given in \cite{Fukaya2009}.

    \begin{theorem}[Theorem 15.2 \cite{fukaya2013lagrangian}]
        For $k\geq1$ integer, there exists:
        \begin{enumerate}
            \item a sequence of non-negative integers $\{n_k\}_{k\geq1}$ such that $n_k<\dim_\R L/2$ and
            \item a sequence of classes $[o^{2n_k}_k(L)]\in H^{2n_k}(L;\Z)$
        \end{enumerate}
        such that, if $o^{2n_k}_k\equiv 0$ for all $k\geq 1$ then, $L$ is unobstructed over $\Lambda_0$ as in Definition \ref{UnobstructedDefn}.
    \end{theorem}
    
    Nevertheless, the arguments and results of Section $3.8$ of \cite{Fukaya2009} do not apply immediately to our setting. For this reason, we prove the following:

\begin{prop}
    Suppose that, $L\subset X$ is a connected, closed, relatively-spin Lagrangian such that:
    \begin{enumerate}
        \item $\mu_L:\pi_2(X, L)\rightarrow \Z$ is trivial and
        \item $H^2(L;\Z)=\{0\}$.
    \end{enumerate}
    Then, $L$ is unobstructed over $\Lambda_0$, as in Definition \ref{UnobstructedDefn}.
\end{prop}

\begin{proof}
    In the view of the first assumption and after some abuse of notation, we consider the $e$-reduced Novikov ring $\Lambda_0:=\{\sum_{\lambda\geq0}a_\lambda T^\lambda:a_\lambda\in \Z, \#\{\lambda<E:a_\lambda\neq0\}<\infty, \forall E\in\R\}$, let $b\in\mathrm{CM}^1(f; \Lambda_0)$ with positive-valuation and denote by $\kappa(b):=\sum_{k\geq0}\m_k(b,\dots,b)$. Let $\mathfrak{G}:=\{\omega(\beta)\in\R:\beta\in \Pi(L), \M(\beta)\neq\emptyset\}$. By Gromov's compactness, we can find a monoid $\mathfrak{G}\subseteq G$ such that $G=\{\lambda_i\in\R:\lambda_{i}\leq\lambda_{i+1}\}$. We will argue by induction over such monoid $G$.

    First, we observe that using Lemma \ref{Morse=FloerenergyZero} and the fact that $b$ has positive-valuation, it follows that if $\kappa(b)\equiv 0\text{ mod }T^\lambda$ and after writing $\kappa(b)=\sum_{\lambda\in G}c_\lambda T^\lambda$ then, $c_\lambda\in \mathrm{CM}^1(f;\Z)$ and $d^\text{Morse}(c_\lambda)=0$. Moreover, by linearity of the $\m_k$ operations, it follows that, for every $\lambda_0>0$ and $a\in\mathrm{CM}^1(f; \Z)$, we can write $\kappa(b+T^{\lambda_0}a)=\kappa(b)+T^{\lambda_0}d^\text{Morse}a+\sum_{\lambda>\lambda_0}a_\lambda T^\lambda$.

    Second, we note that, as $\Z$ is a PID and since $\mathrm{CM}^2(f;\Z)$ is a free $\Z$-module, we have that $\ker(d^\text{Morse}:\mathrm{CM}^2(f;\Z)\rightarrow \mathrm{CM}^3(f; \Z))$ is a free abelian group. On the other hand, under the second assumption, it follows that there exists a splitting $h:\ker(d^\text{Morse}:\mathrm{CM}^2(f;\Z)\rightarrow \mathrm{CM}^3(f; \Z))\rightarrow\mathrm{CM}^1(f;\Z)$.

    Now we construct a bounding cochain $b\in\mathrm{CM}^1(f; \Lambda_0)$ inductively over $G$. For the base case, we set $b\equiv 0$ as the Morse differential is unobstructed and after using Lemma \ref{Morse=FloerenergyZero}. Using the ordering of $G$, we set $b^0=0$ and $b^j:=\sum_{i\leq j}a_{\lambda_i}T^{\lambda_i}$ such that $b^j=b^{j-1}+T^{\lambda_j}a_j$ where $a_j:=-d^\text{Morse}c_{\lambda_j}$. Then, $\kappa(b)\equiv 0\text{ mod }T^\lambda$ for every $\lambda\geq0$. By seperatedness property of the $T$-filtration on $\Lambda_0$, it follows that $\kappa(b)=0$. Therefore, $L$ is unobstructed over $\Lambda_0$, as desired.
\end{proof}

\begin{cor}\label{unobstructedSpheres}
    Any Lagrangian which is an integral homology sphere in a Calabi-Yau manifold of real dimension at least $6$, is unobstructed over $\Lambda_0$.
\end{cor}

\subsubsection{Spectral Sequence}\label{spectralSeqSec}
Now suppose that $(L,b)$ is unobsructed over $\Lambda_0$ and to avoid notational clutter, we denote by $C:=\mathrm{CM}^*(f; \Lambda_0)$ and $\overline{C}=\mathrm{CM}^*(f; \Z)$ it $\Z$-reduction. Moreover, denote by $d^\text{Morse}:C\rightarrow C$ the Morse differential on $C$.

\begin{lemma}\label{Morse=FloerenergyZero}
    $\m^b_{1, \beta=0}\equiv d^\text{Morse}$.
\end{lemma}

\begin{proof}
    This follows by construction of $\m_{1, \beta=0}$ as in \cite{rabah2024fukaya}, after noting that $b$ has positive-valuation. (c.f. Corollary $4.32$ of \cite{rabah2024fukaya}).
\end{proof}

We will use the discrete energy $F$-filtration on $C$, as in Equation $(6.3.19)$ of \cite{fukaya2010lagrangian}. Namely, since the energy of pseudo-holomorphic curves are $\textit{quantized}$, it follows that there exists a real number $\eta>0$ such that if $p\geq0$ integer then $F^pC= T^{p\eta}C$ and $F^pC=C$ if $p<0$. Moreover, we will use the indexing convention as in \cite{Fukaya2009}, \cite{fukaya2010lagrangian} and denote by $gr_p\Lambda_0:=F^p\Lambda_0/F^{p+1}\Lambda_0$ and 
\[
F^pHF^k((L,b);\Lambda_0):=\mathrm{im}(H^k(F^pC, \m_1^b)\rightarrow HF^k((L,b);\Lambda_0)).
\]

\begin{theorem}[Spectral Sequence]\label{SpecSeqThmI}
    There exists a spectral sequence $\{E^{k,p}_r\}$ of differentials $\delta_r:E^{k,p}_r\rightarrow E^{k+1, p+r-1}_r$, whose $E_2$-page is of the form $E^{k,p}_2\cong \bigoplus_j HM^*(L;\Z)\otimes_\Z (gr_p\Lambda_0)^{k-j}$, where $HM^*(L;\Z)$ is the Morse cohomology of $L$ given by the Morse-Smale pair $(f,g)$. Moreover, this spectral sequence stabilizes. That is, there exists an $r_0\geq0$, independent of $k,p$, such that for all $r\geq r_0$, we have $\delta_r\equiv 0$ and $E^{k,p}_r\cong E^{k,p}_\infty\cong gr_pHF^k((L,b);\Lambda_0)$. Furthermore, the $F$-filtration on $HF^*((L,b); \Lambda_0)$ is exhaustive, complete and separated.
\end{theorem}

The proof that we will present here is based on our understanding of Chapter $6$ of \cite{fukaya2010lagrangian}, yet we try to not argue in the complete generality as they do. We start by fixing some notations following Chapter $1$ of \cite{Fukaya2009}:\\
Set $K:=\Z[e^{-1}, e]$ and denote by $\Lambda$ the universal Novikov ring over the integers. Namely,
\begin{equation} \label{NovikovFieldUsher}
    \Lambda:=\{\sum_{\lambda\in\R}a_\lambda T^\lambda:a_\lambda\in K, \#\{\lambda<E:a_\lambda\neq0\}<\infty, \forall E\in\R\}.
\end{equation}

For $a\in\Lambda$ for the form $a=\sum_{\lambda\in\R}a_\lambda T^\lambda$, we will denote its valuation by $v(a):=\min\{\lambda: a_\lambda\neq 0\}$. Denote by $d^\text{Floer}:=\m_1^b$ and $\tilde C:=\overline{C}\widehat{\otimes}_\Z\Lambda$. Following Equations $(6.3.18), (6.3.19)$ of \cite{fukaya2010lagrangian}, we consider
\[
Z_r^{k,p}:=(F^pC^k\cap (d^\text{Floer})^{-1}(F^{p+r-1}C^{k+1}))+F^{p+1}C^k,
\]
\[
B_r^{k,p}:=(d^\text{Floer}(F^{p-r+2}C^{k-1})\cap F^pC^k)+F^{p+1}C^k
\]
and set
\[
E^{k,p}_r:=Z_r^{k,p}/B_r^{k,p}.
\]
We note that, as $d^\text{Floer}\circ d^\text{Floer}=0$, it follows that $B_r^{k,p}\subseteq Z_{r}^{k,p}$. Now for $r\geq 1$ integer, we define \[
\delta_r:E^{k,p}_r\rightarrow E_r^{k+1, p+r-1}\]by $\delta_r[x]:=[d^\text{Floer}x]$.

\begin{lemma}[Lemma $6.3.20$ \cite{fukaya2010lagrangian}]\label{LemmaSpecSeq}
   $\delta_r$ is well-defined and $\delta_r\circ\delta_r\equiv0$. Moreover, we have an identification $E^{k,p}_{r+1}\cong \ker(\delta_r)/\mathrm{im}(\delta_r)$.
\end{lemma}

\begin{proof}
    Let $[x]\in E^{k, p}_r$ and suppose that such class is also represented by $y\in Z^{k,p}_r$. Then, $x-y=d^\text{Floer}z+a$ such that $z\in F^{p-r+2}C^{k-1}$ and $a\in F^{p+1}C^k$. In particular, $d^\text{Floer}a\in F^{p+r-1}C^{k+1}$ and as $(p+r-1)-r+2=p+1$, it follows that $\delta_r$ is well-defined. On the other hand, $\delta_r\circ\delta_r=0$ as $d^\text{Floer}\circ d^\text{Floer}=0$. Similar reasoning and using the Isomorphism Theorem of modules, shows the last claim. 
\end{proof}

\begin{lemma}\label{LemmaSpecSeqI}
    The first two pages of the spectral sequence $\{E^{k,p}_r\}$ and of the form:
    \begin{enumerate}
        \item $E^{k,p}_1\cong (\overline{C}\otimes_\Z gr_p\Lambda_0)^k$ and $\delta_1\equiv d^\text{Morse}$.
        \item $E^{k,p}_2\cong \bigoplus_j HM^*(L;\Z)\otimes_\Z(gr_p\Lambda_0)^{k-j}$.
    \end{enumerate}
\end{lemma}

\begin{proof}
    This follows from Lemma \ref{LemmaSpecSeq} and Lemma \ref{Morse=FloerenergyZero}.
\end{proof}

It is at this step where we defer from the more general discussion as in \cite{fukaya2010lagrangian} and make use of the fact that $\Tilde{C}$ is a finite-dimensional $\Lambda$-vectorspace by applying Theorem $2.5$ of \cite{usher2008spectral}.

\begin{lemma}\label{LemmaSpecSeqII}
    There exists a real number $\lambda_0\geq0$ such that for every $x\in\mathrm{im}(d^\text{Floer}:\tilde C\rightarrow\tilde C)$ there exists a $y\in \Tilde{C}$ where, $d^\text{Floer}y=x$ and $v(y)\geq v(x)-\lambda_0$. 
\end{lemma}

\begin{proof}
    As noted above, $\tilde C$ is a finite-dimensional $\Lambda$-vectorspace, say $\dim_\Lambda\tilde C\equiv N$. After fixing a basis of $\tilde C$ and using the same notation as in \cite{usher2008spectral} or Chapter $1$ of \cite{Fukaya2009}, we get an identification $\tilde C\cong \Lambda_K(\R)^N\equiv\Lambda^N$. Noting that $K\equiv \Z[e^{-1}, e]$ is a Noetherian ring, we apply Theorem $2.5$ of \cite{usher2008spectral} to the matrix representation of $d^\text{Floer}$ given by the above identification.
\end{proof}

Now set $Z^k:=\ker(d^\text{Floer}:C^k\rightarrow C^{k+1})$ and $B^k:=\mathrm{im}(d^\text{Floer}:C^{k-1}\rightarrow C^k)$. Moreover, choose an integer $r_0$ such that $(r_0-2)\eta>\lambda_0$.

\begin{lemma}\label{LemmaSepSeqIII}
    Given integers $k,p$ and $r$ such that $p\geq 0$ and $r\geq r_0$ then, $Z^{k,p}_r=(Z^k\cap F^pC^k)+F^{p+1}C^k$, $B^{k,p}_r=(B^k\cap F^p C^k)=F^{p+1}C^k$ and $\delta_r=0$.
\end{lemma}

\begin{proof}
    We only show that $Z^{k,p}_r\subseteq (Z^k\cap F^pC^k)+F^{p+1}C^k$ as the opposite inclusion follows from their definitions. To do so, let $x\in Z^{k,p}_r$ such that $x\in F^pC^k$ and $d^\text{Floer}x\in F^{p+r-1}C^{k+1}$. By Lemma \ref{LemmaSepSeqIII}, it follows that there exists $y\in\tilde C$ such that $d^\text{Floer}y\equiv d^\text{Floer}x$. Moreover, $v(y)\geq(p+r-1)\eta-\lambda_0>(p+1)\eta$. Therefore, $y\in F^{p+1}C^{k}$ and $x-y\in Z^k\cap F^pC^k$. Similar reasoning shows that $B^{k,p}_r=(B^k\cap F^p C^k)=F^{p+1}C^k$. As for the last claim, it follows from the definition of $\delta_r$ and $Z^k$.
\end{proof}

\begin{proof}[Proof of Theorem \ref{SpecSeqThmI}]
    After using Lemmas \ref{LemmaSpecSeqI}, \ref{LemmaSpecSeqII} and \ref{LemmaSepSeqIII}, the only thing that remains to prove is that the $F$-filtration on $HF^*((L, b); \Lambda_0)$ is exhaustive, complete and separated. Indeed, exhaustive follows from the fact that $F^0C\equiv C$. While separatedness follows from the fact that $B^k\subseteq C^k$ is closed in the $T$-adic topology. Furthermore, completeness follows from the fact that $Z^k\subseteq C^k$ is closed and complete in the $T$-adic topology.
\end{proof}

In what follows, we will give sufficient conditions for the degeneration of the spectral sequence in Theorem \ref{SpecSeqThmI} at its $E_2$-page. Namely, we will prove the following:

\begin{theorem}\label{degenerationSpecSeq}
    Let $L$ be a connected, closed, relatively-spin, unobstructed over $\Lambda_0$ Lagrangian and let $b$ be a bounding cochain of $L$. If the canonical inclusion $\iota:L\hookrightarrow X$ induces an injective group morphism $\iota_*:H_*(L;\Z)\hookrightarrow H_*(X;\Z)$ then, the spectral sequence $\{E^{k,p}_r\}$ of Theorem \ref{SpecSeqThmI} degenerates at its $E_2$-page.
\end{theorem}

To give a proof of the above Theorem and following Section $6$ of \cite{fukaya2010lagrangian} and Section $15$ of \cite{fukaya2013lagrangian}, we will consider the Open-Closed map of Theorem \ref{OCmapTheorem}. Let $g_X$ be a metric on $X$ extending that on $L$. That is, $\iota^*g=g_X$ and let $F:X\rightarrow\R$ be a Morse function on $X$ such that the pair $(F, g_X)$ is a Morse-Smale pair and $F$ has the form
\[
F_X(x,v)\equiv f(x)+|v|^2
\]
on a tubular neighborhood of $L$.

\begin{lemma}\label{GysinMapLemma}
    The inclusion map $\mathrm{Crit}(f)\hookrightarrow\mathrm{Crit}(F)$ is a chain map inducing the Gysin morphism $\iota_!:H^{*+n}(X;\Z)\rightarrow H^*(L;\Z)$.
\end{lemma}

To avoid notational clutter, let $B:=\mathrm{CM}^*(F;\Lambda_0)[n]$ and $d_B:=(-1)^nd^\text{Morse}_X$. Note that, $d_B$ is $F$-filtered differential and hence by standard spectral sequence constructions (c.f. Theorem $5.4.1$ of \cite{weibel1994introduction}), it follows that there exists a spectral sequence $\{E^{k,p}_r(B)\}$ collapsing at its $E_2$-page, which is of the form
\[
E^{k,p}_2(B)\equiv\bigoplus_jHM^{j+n}(F;\Z)\otimes_\Z(gr_p\Lambda_0)^{k-j},
\]
with differential $\delta_r^B\equiv0$ for all $r\geq 2$.
\begin{lemma}\label{OCfilterdChainMapSpec}
    There exists a filtered chain map $P_r:E_r(C)\rightarrow E_r(B)$ such that $P_2\equiv\bigoplus_j \iota_!\otimes\Id_{(gr_p\Lambda_0)^{k-j}}$.
\end{lemma}

\begin{proof}
    For $x\in\mathrm{Crit}(f)$, we set $P_r(x):=\sum_{j\geq0}\p_{j+1}(x;b,\dots,b)$ and define $P_r$ after extending $\Lambda_0$-linearily. The fact that $P_r$ is a filtered chain map, follows from Theorem \ref{OCmapTheorem}, and completeness and seperatedness property of the $F$-filtration on $\Lambda_0$. While the second claim follows from Lemma \ref{GysinMapLemma} and the gappedness condition on $\Lambda_0$.
\end{proof}

Observe that, as $\delta_r^B\equiv 0$ for all $r\geq 2$, it follows that $P_r\circ\delta_r=\delta_r^B\circ P_r=0$ for $r\geq2$. In particular, $\mathrm{im}(\delta_r)\subset\ker(P_r)$ for $r\geq2$. Now for $r\geq2$, define $K_r$ inductively by
\[
K_{r+1}:=\frac{K_r\cap\ker(\delta_r)}{K_r\cap\mathrm{im}(\delta_r)}
\]
with $K_2:=\ker(P_2)$.

\begin{lemma}\label{degenerationLemma}
    For $r\geq2$, we have $K_r=\ker(P_r)$.
\end{lemma}

\begin{proof}
    We argue by induction on $r$. Indeed, assume that $K_r=\ker(P_r)$ for some $r\geq 2$ and consider a $\delta_r$-class $[a]\in E_r(C)$. Write such class as $a+\delta_rb$ and observe that $P_r(a+\delta_rb)=P_r(a)$. Therefore, $\ker(P_{r+1})=\frac{\ker(P_r)\cap\ker(\delta_r)}{\mathrm{im}(\delta_r)}$ which in return is identitically equal to $\frac{K_r\cap\ker(\delta_r)}{K_r\cap\mathrm{im}(\delta_r)}$.
\end{proof}

\begin{proof}[Proof of Theorem \ref{degenerationSpecSeq}]
After noting that $\iota_*:H_*(L;\Z)\hookrightarrow H_*(X;\Z)$ being an injective morphism implies that $\iota_!:H^{*+n}(X;\Z)\rightarrow H^*(L;\Z)$ is an injective morphism and noting that $(gr_p\Lambda_0)^{k-j}$ is a free $\Z$-module and hence by exactness of tensor products, we have $K_2\equiv\{0\}$; the result follows from Lemma \ref{degenerationLemma}. 
\end{proof}

\subsection{Pair of Transversely Intersecting Lagrangians}

Let $L_0,L_1$ be a pair of transversaly intersecting, connected, closed, relatively-spin as a pair Lagrangians of $X$. For $i=0,1$, let $(f_i,g_i)$ be a Morse-Smale pair on $L_i$ and consider the $A_\infty$-algebras $(\mathrm{CM}^*(f_i;\Lambda_0),\{\m_k^i\}_{k\geq0})$ as given by the main theorem of \cite{rabah2024fukaya}, \cite{rabah2025fukaya}. 

\begin{defn}[Unobstructed over $\Lambda_0$ as a pair]\label{UnobstructedPair}
    The pair $L_0,L_1$ is said to be unobstructed over $\Lambda_0$ if each of $L_0$ and $L_1$ is unobstructed over $\Lambda_0$ as in Definition \ref{UnobstructedDefn}.
\end{defn}

Suppose that $L_0,L_1$ are unobstructed over $\Lambda_0$ as a pair and let $b_i$ be a choice of bounding cochain on the $A_\infty$-algebra $(\mathrm{CM}^*(f_i;\Lambda_0),\{\m_k^i\}_{k\geq0})$ for $i=0,1$. Consider, \[
CF((L_0,b_0),(L_1,b_1);\Lambda_0):=CF(L_0,L_1;\Lambda_0)
\] where $CF(L_0,L_1;\Lambda_0)$ is as in Definition \ref{FloerComplexoverRing}. Using the identification $CF(L_0,L_1;\Lambda_0)\cong \widehat{\bigoplus}_{p\in L_0\cap L_1}\Lambda_0\langle p\rangle$ we consider the $\{\n_{k_0,k_1}\}_{k_0,k_1\geq0}$ the $A_\infty$-bimodule structure on $CF(L_0,L_1;\Lambda_0)$ as given by Theorem \ref{AinftyBimoduleThm}.

\begin{defn}[Lagrangian Floer Homology of Transversely Intersecting Lagrangians]\label{HF2Lags Defn}

We define $\delta_{b_0,b_1}:CF((L_0,b_0), (L_1,b_1);\Lambda_0)\rightarrow CF((L_0,b_0), (L_1,b_1);\Lambda_0)$ by 
\[
\delta_{b_0,b_1}(p):=\sum_{k_0,k_1\geq0}\n_{k_0,k_1}(b_0^{\otimes k_0},p,b_1^{\otimes k_1}).
\]
Moreover, we define 
\[
HF((L_0,b_0),(L_1, b_1);\Lambda_0):=\ker(\delta_{b_0,b_1})/\mathrm{Im}(\delta_{b_0,b_1}),
\]
which we will call Lagrangian Floer homology modules of the pair $(L_0,b_0)$ and $(L_1,b_1)$ over $\Lambda_0$.
\end{defn}

\begin{prop}\label{UnobstructedHamiltonianImage}
    If $L$ is unobstructed over $\Lambda_0$, $b$ is a bounding cochain for $L$ and $\phi:X\rightarrow X$ is a Hamiltonian diffeomorphism then, $\phi(L)$ is unobstructed over $\Lambda_0$ and $\phi_*(b)$ is a bounding cochain of $\phi(L)$.
\end{prop}

\begin{proof}
We first note that, as $\phi:X\rightarrow X$ is a Hamiltonian diffeomorphism and in particular, $\phi$ is isotopic to the identity, it follows that $\phi(L)$ is a relatively-spin Lagrangian with respect to the same background class $st\in H^2(X;\Z/2)$ as for $L$. It is important to recall that our definition of unobstructedness over $\Lambda_0$, as in Definition \ref{UnobstructedDefn}, depends on the choice of Morse-Smale pair we put on the Lagrangian (in \cite{RabahIntegralFukaya}, we show that different choices give pseudo-isotopic $A_\infty$-algebras). With this in mind, let $(f,g)$ be any Morse-Smale pair on $L$ and denote by:
\begin{enumerate}
    \item $f_\phi:=f\circ\phi^{-1}$,
    \item $g_\phi:=(\phi^{-1})^*g$ and 
    \item $J_\phi$ an almost complex structure on $X$ given by $(J_{\phi})_{\phi(x)}:=d_x\phi\circ J_x\circ(d_x\phi)^{-1}$ for any $x\in X$.
\end{enumerate}

As $\phi:X\rightarrow X$ is an orientation-preserving diffeomorphism, where $X$ is given the orientation defined by $\omega^n$, it follows that the pair $(f_\phi,g_\phi)$ is a Morse-Smale pair on $\phi(L)$. By Chain Rule, it follows that $\nabla^{g_\phi}f_\phi=d\phi(\nabla^gf)$ and by compactness of $X$, we have that $||d\phi||_\infty\leq C<\infty$ is uniformly bounded. Therefore, the assignment
\[
I_\phi:\gamma\mapsto\phi\circ\gamma,
\]
defines an isomorphism between the compactified moduli space of trajectories of $(f,g)$ and $(f_\phi, g_\phi)$. Similar reasoning implies that if $\tilde f:\widetilde{\overline{\mathrm{Gr}}}_{k+1}\times L\rightarrow \R$ is a domain-dependent Morse function associated to $f$ and $\tilde f_\phi:\widetilde{\overline{\mathrm{Gr}}}_{k+1}\times \phi(L)\rightarrow \R$ is given by $\tilde f_\phi((T,p),x):=\tilde f((T,p),\phi^{-1}(x))$ then, $I_\phi$ defines an isomorphism between the compactified moduli space of domain-dependent Morse trajectories of $(f,g)$ and $(f_\phi, g_\phi)$.

Now to deal with our construction of Global Kuranishi charts, we denote by:
\begin{enumerate}
    \item $\Omega_\phi:=(\phi^{-1})^*\Omega$ where $\Omega$ is as in Lemma \ref{Existinceof2-form} and 
    \item $( H_\phi,\nabla_\phi,\tau_\phi):=(\phi^{-1})^*(H,\nabla,\tau)$ where $(H,\nabla,\tau)$ is as in Lemma \ref{LineBundleLemma}.
\end{enumerate}
If $u:(\Sigma,\partial\Sigma)\rightarrow(X,L)$ is a smooth map from a bordered nodal Riemann surface of genus-zero and we write $u_\phi:=u\circ\phi$, we get that 
\[
u^*_\phi(H_\phi,\nabla_\phi,\tau_\phi)\cong u^*(H,\nabla,\tau),
\]
\[
u^*_\phi\Omega_\phi\equiv u^*\Omega.
\]
Therefore, the resulting space of framing as in Definitions \ref{spaceFramings} and \ref{framedStripsConstraint} is the same and hence the consistent choice of domain metric as in Definition \ref{consistentMetricDefn} can be taken the same. For such choice and as $u^*_\phi\Omega\equiv u^*\Omega$, it follows that for a given choice of framing $F$, we have the matrix $A_F^\phi\equiv A_F$ as given in item $(3)$ of Definition \ref{FramedStrips}. With this in mind, let $\mathcal C^*\rightarrow \F$ be the universal framed curve away from its nodes and marked points and suppose that $(V_m,\lambda_m)$ is a finite-dimensional approximation scheme for $\bigwedge^{0,1}_{\mathcal{C}^*/\F}\otimes_\C TX$ with respect to the almost complex structure $J$ on $X$. The map given by $(c,x;\eta)\mapsto(c,\phi(x);d_x\phi\circ\eta)$ defines a bundle isomorphism between $\bigwedge^{0,1}_{\mathcal{C}^*/\F}\otimes_\C TX$ for $(J,J_\phi)$ respectively. Now set
\[
\lambda^\phi_m(\eta)(c,\phi(x)):=d_x\phi\circ\lambda_m(\eta)(c,x).
\]
Then, $(V_m,\lambda_m^\phi)$ is a finite-dimensional approximation scheme for the moduli space of Floer-Morse trajectories, as in Definition $4.29$ of \cite{rabah2024fukaya}. 

Note that, as $\phi$ is a Hamiltonian diffeomorphism, it follows that for a given $\beta\in\pi_2(X,L)$, we have $\omega(\beta)=\omega(\phi_*(\beta))$ and $\mu_L(\beta)=\mu_{\phi(L)}(\phi_*(\beta))$. Now by definition of $J_\phi$ and using Chain Rule, it follows that $\bar\partial_{J_\phi}(\phi\circ u)=d\phi\circ\bar\partial_J u$. Thus,
\[
\bar\partial_{J_\phi}(\phi\circ u)+\lambda^\phi_m(\eta)\circ\Gamma_{\phi\circ u}=d\phi\circ(\bar\partial_J u+\lambda_m(\eta)\circ\Gamma_u).
\]
Therefore, for a given $\beta\in\pi_2(X,L)$ and if we denote by $(G,\T,\E,s)$ the Global Kuranishi chart representation of $\M(X,L;J,\beta)$ as given by Theorem $3.51$ of \cite{rabah2024fukaya}, the assignment $I_\phi:\T\rightarrow\T^\phi$ given by 
\[
((\Sigma,\partial\Sigma,\underline{z}),u,F,\eta;r)\mapsto((\Sigma,\partial\Sigma,\underline{z}),u,F,\eta;r)
\]
define a Thickening $\T^\phi$ on $\M(X,\phi(L);J_\phi,\phi_*(\beta))$ such that if $\T^\phi$ is given the induced smooth $G$-equivariant atlas from the one $\T$ then, $I_\phi$ is a smooth $G$-equivariant diffeomorphism. Define $\E^\phi$ to be the trivial vectorbundle over $\T^\phi$ of fibre over $p\in\T^\phi$ given by $\E_{I_{\phi^{-1}}(p)}$. Similarly, define the Kuranishi section $s^\phi:\T^\phi\rightarrow \E^\phi$ by $s^\phi(p):=s(I_{\phi^{-1}}(p))$. Such construction and using Theorem $3.51$ of \cite{rabah2024fukaya}, provides us with $(G,\T^\phi,\E^\phi,s^\phi)$ a Global Kuranishi chart representation of $\M(X,\phi(L);J_\phi,\phi_*(\beta))$ together with a commutative diagram as smooth $G$-manifolds:
\[\begin{tikzcd}
	{\mathcal E} && {\mathcal E^\phi} \\
	\\
	{\mathcal T} && {\mathcal T^\phi}
	\arrow["{\hat I_\phi}", from=1-1, to=1-3]
	\arrow[from=1-1, to=3-1]
	\arrow[from=1-3, to=3-3]
	\arrow["s", bend left=30,  from=3-1, to=1-1]
	\arrow["{I_\phi}"', from=3-1, to=3-3]
	\arrow["{s^\phi}"', bend right=30, from=3-3, to=1-3]
\end{tikzcd}\]
From the commutativity of the above diagram as smooth $G$-manifolds and noting that $\hat I_\phi$ is an equivariant bundle isomorphism covering $I_\phi$ and since the universal zero locus as given in Section $2.1.2$ of \cite{bai2022arnold} has a diffeomorphism invariant Whitney stratification (c.f. Appendix $B$ of \cite{rabah2024fukaya}), it follows that for a given FOP-perturbation $s_\epsilon$ of $s$, we can find an FOP-perturbation $s^\phi_\epsilon$ of $s^\phi$ such that $s^\phi_\epsilon\circ I_\phi\equiv \hat I_\phi\circ s_\epsilon$.

Having the above constructions at hand and using the arguments of the proofs of Theorem $4.50$, Theorem $5.3$ and Corollary $5.4$ of \cite{rabah2024fukaya}, it follows that we can assign an $A_\infty$-algebra structure $\{\m_k^\phi\}_{k\geq0}$ on the Morse complex $\mathrm{CM}^*(f_\phi;\Lambda_0)$ together with a strict, filtered, curved $A_\infty$-isomorphism \[
\Phi:(\mathrm{CM}^*(f;\Lambda_0),\{\m_k\}_{k\geq0})\rightarrow (\mathrm{CM}^*(f_\phi;\Lambda_0),\{\m_k^\phi\}_{k\geq0})
\]
in the sense of \cite{FukayaUnobstructed}, given by $\Phi(x):=\phi(x)$ for $x\in\mathrm{Crit}(f)$ such that for any integer $k\geq0$, we have 
\[
\m_k^\phi(\phi(x_1),\dots,\phi(x_k))=\Phi(\m_k(x_1,\dots,x_k)).
\]
\end{proof}

\begin{defn}[Full Novikov Ring over the Integers]
    We denote by 
    \[
    \Lambda:=\Lambda_0[T^{-1}],
    \]
    the $\textit{full}$ Novikov ring over the integers.
\end{defn}

\begin{remark}
    It is important to note that $\Lambda$ is not a field as $2\in\Lambda$ is not invertible. For this reason, we call it the full Novikov ring over the integers, while call $\Lambda_0$ the Novikov ring over the integers (c.f. \cite{fukaya2013lagrangian}).
\end{remark}

\begin{lemma}\label{generatorBoundOverField}
    Suppose that $L_0,L_1$ are transversely intersecting unobstructed Lagrangians over $\Lambda_0$ as in Definition \ref{UnobstructedPair} and let $b_0, b_1$ be a choice of bounding cochains on $L_0,L_1$, respectively. Then, the minimum number of generators of $HF((L_0,b_0),(L_1,b_1);\Lambda_0)$ is greater than or equal to $\#(L_0\cap L_1)$, the number of intersection of points between $L_0$ and $L_1$.
\end{lemma}

\begin{proof}
    This is an algebraic lemma and follows from the fact that the Novikov ring
    \[
    R:=\{\sum_i a_iT^{\lambda_i}:a_i\in\Z,\#\{\lambda\leq E\}<\infty,\forall E\in\R\},
    \]
    is a PID by \cite{hofer1995floer}. After taking $\text{mod } 2$ grading on $CF(L_0,L_1;\Lambda)$, we can decompose 
    \[
    CF(L_0,L_1;\Lambda)\cong CF^\text{even}(L_0,L_1;\Lambda)\oplus CF^\text{odd}(L_0,L_1;\Lambda).
    \]
    As $\delta_{b_0,b_1}$ is of degree-$1$ and noting that $CF(L_0,L_1;\Lambda)\cong CF(L_0,L_1;R)\otimes_R\Lambda$, it follows that we can decompose $\delta_{b_0,b_1}$ under the above decomposition into
    \[
    A:=\bar\delta_{b_0,b_1}|_{CF^\text{even}(L_0,L_1;R)}:CF^\text{even}(L_0,L_1;R)\rightarrow CF^\text{odd}(L_0,L_1;R),
    \]
    \[
    B:=\bar\delta_{b_0,b_1}|_{CF^\text{odd}(L_0,L_1;R)}:CF^\text{odd}(L_0,L_1;R)\rightarrow CF^\text{even}(L_0,L_1;R),
    \]
    with the understanding that $\delta_{b_0,b_1}\equiv\bar\delta_{b_0,b_1}\otimes_R\Lambda$ and using the fact that $\bullet\otimes_R\Lambda$ is exact as $R$ is a PID. Under the above decomposition of $CF(L_0,L_1;R)$, we can write 
    \[
\bar\delta_{b_0,b_1}=
\begin{pmatrix}
0 & eB\\
A & 0
\end{pmatrix}.
\]
Once again, by exactness of $\bullet\otimes _R\Lambda$, it follows that 
\[
\ker(A\otimes_R\Lambda)\cong\ker(A)\otimes_R\Lambda,
\]
\[
\ker(B\otimes_R\Lambda)\cong\ker(B)\otimes_R\Lambda.
\]
As $e$ is invertible in $\Lambda$, it follows that 
\[
\ker(\delta_{b_0,b_1})\cong(\ker(A)\oplus\ker(B))\otimes_R\Lambda.
\]
Now the result follows from the fact that $R$ is PID and hence $\ker(A)$ and $\ker(B)$ are free $R$-modules.
\end{proof}

\begin{theorem}\label{HamiltonianInvariance}
    Suppose that $L$ is unobstructed over $\Lambda_0$, $b$ is a bounding cochain for $L$ and $\phi:X\rightarrow X$ is a Hamiltonian diffeomorphism such that $L$ and $\phi(L)$ intersect transversely. Then, we have the following identifications $HF((L,b);\Lambda)\cong HF((L,b),(\phi(L), \phi_*(b));\Lambda)\cong HF((\phi(L), \phi_*(b));\Lambda)$ as $\Lambda$-modules.
\end{theorem}

Lemma \ref{generatorBoundOverField} together with Theorem \ref{HamiltonianInvariance} will allow us to use the $\textit{Spectral Sequence}$ results above to obtain stronger bounds on the number of intersection points of a Lagrangian with its Hamiltonian image (c.f. Section \ref{AppSection}). 

\begin{proof}[Proof of Theorem \ref{HamiltonianInvariance}]
    We will use the moduli space of pseudo-holomorphic curves with $\textit{moving Lagrangian }$ $\textit{boundary}$ conditions, as in Sections $5.3.2$ of \cite{Fukaya2009}, to construct a filtered, curved $A_\infty$-bimodule morphism between the diagonal bimodule structure given by the $A_\infty$-algebra $(\mathrm{CM}^*(f;\Lambda),\{\m_k\}_{k\geq0})$ and the $A_\infty$-bimodule structure on $CF(L,\phi(L);\Lambda)$ as given above. 

Let $H:[0,1]\times X\rightarrow\R$ be a time-dependent generating the given Hamiltonian diffeomorphism $\phi:X\rightarrow X$. That is, we use the sign convention 
\[
dH=\omega(V,.),
\]
for the Hamiltonian vectorfield $V$ given by $H$ and note that $\phi\equiv \phi^1_V$ the time-$1$ flow of $V$. To avoid notational clutter and for $r\in[0,1]$, we denote by $\phi_r:=\phi^r_V$ the time-$r$ flow of $V$.

Let $\chi:\R_s\rightarrow[0,1]_r$ be a smooth bump function such that:
\[
\chi(s)=
\begin{cases}
    &1\text{ for }s\geq K\\
    & 0\text{ for }s\leq -K
\end{cases}
\]
where $K>0$ is a real number. Let $\{J_{s,t}\}$ be a $\R\times[0,1]$-smooth family of $\omega$-tame almost complex structures on $X$ such that:
\begin{enumerate}
    \item $J_{0,t}\equiv J$,
    \item $J_{s,0}\equiv J$ and 
    \item $J_{s,1}\equiv (\phi_s)_*J=d\phi_s\circ J\circ (d\phi_s)^{-1}$.
\end{enumerate}
We will be interested in the moduli space of Floer-Morse trajectories, modeled over two-paged metric ribbon trees, having a distinguished spine vertex $v_*$ in the moduli space of maps of the form $u:\R_s\times[0,1]_t\rightarrow X$ such that:
\begin{enumerate}
    \item $\partial_s u+J_{\chi(s),t}(u)\partial_tu=0$,
    \item $u(s, 0)\in L$,
    \item $u(s,1)\in L_{\chi(\tau)}:=\phi_{\chi(s)}(L)$ and
    \item $\int |\partial_s u|^2\leq C<\infty$.
\end{enumerate}

Now let $(f,g)$ be any Morse-Smale pair on $L$ and for $r\in[0,1]$ we set:
\begin{enumerate}
\item $f_r:=f\circ\phi_r^{-1}$,
    \item $g_r:=(\phi^{-1}_r)^*g$ and
\end{enumerate}

Suppose that $\tilde f$ is a generic domain-dependent Morse function associated to $f$, as in Theorem $4.50$ of \cite{rabah2024fukaya} and let $\Gamma$ be a two-paged, possibly broken, metric ribbon tree. All the edge on the $0$-page of $\Gamma$ will correspond to domain-dependent Morse trajectories of $\tilde f$. We use the ordering of the spine edge given by the minimal path joining $e_-$ to $e_+$. Suppose that $v$ is a spine vertex between $e_-$ and $v_*$ and suppose that $B$ is a branch of $\Gamma$ in the $1$-page attached to $v$. Then all the edges of $B$ will correspond to domain-dependent Morse trajectories of the restriction of $\tilde f$ onto such branch. On the hand, if $v$ is between $v_*$ and $e_+$, then all the edges of $B$ will correspond to domain-dependent Morse trajectories of $\phi^*(\tilde f|_{B,e})(B,s;x):=(B,s;\tilde f|_{B,e}(B,s;\phi^{-1}(x)))$. While, if $B$ is a branch of the $1$-page of $\Gamma$ that is attached to $v_*$, then its edges will correspond to domain-dependent Morse trajectories of $f^\text{move}(B,s;x):=\tilde f|_{B,e}(B,s;\phi_{\chi(r)}^{-1}(x))$, where $r\in[0,1]$ such that $(r,1)\in\R\times[0,1]$ is the point of the domain of the corresponding pseudo-holomorphic curve associated to $v_*$ such that it is attached to the root of $B$. Now observe that, using the same arguments of the proof of Proposition \ref{UnobstructedHamiltonianImage}, for any $r\in[0,1]$, the assignment $\gamma\mapsto \phi_r\circ\gamma$ gives an orientation-preserving diffeomorphism between the compactified moduli space of Morse trajectories of $(f,g)$ and $(f_r,g_r)$. Then using parametric implicit function theorem and parametric Sard's Theorem, it follows that for a generic $g$, the compactified moduli space of Morse trajectories of $(f_r,g_r)$ forms a smooth $[0,1]_r$-family. On the other hand, from the proof of Lemma $4.22$ of \cite{rabah2024fukaya}, $\tilde f$ is of the form of $\tilde f|_{B,e}(B,s;x)=f(x)+\beta_e(s)h(x)$ where $\beta_e$ is a smooth bump function. Hence for any $r\in[0,1]$, $\phi_r^*(\tilde f)|_{B,e}=f(\phi^{-1}_r(x))+\beta_e(s)h(\phi^{-1}_r(x))$. Therefore, using the arguments of Lemma $4.28$ and Theorem $4.31$ of \cite{rabah2024fukaya}, we find that for a generic choice of domain-dependent Morse function $\tilde f$ associated to $f$ the moduli space of domain-dependent Morse trajectories as described above for a fixed $\Gamma$ two-paged metric ribbon tree, is a smooth $\overline{T}_{k_0,k_1}$-family, where $k_i\geq0$ is an integer associated to the number of leaves on the $i^{th}$-page.

Now we deal with the auxiliary data needed to give the moduli space of pseudo-holomorphic curves associated to vertices, a Global Kuranishi chart representation. We first note that, for any $r\in[0,1]$ $\phi_r$ is a Hamiltonian diffeomorphism and hence for any given $\beta\in\pi_2(X,L)$, we have $\omega(\beta)=\omega((\phi_r)_*(\beta))$ and $\mu_L(\beta)=\mu_{\phi_r(L)}((\phi_r)_*(\beta))$. Let $\Gamma$ be a two-paged metric ribbon tree and recall that we have a distinguished spine vertex $v_*$ such that we associate to it a pseudo-holomorphic curve as defined above. If $v$ is a spine vertex between $e_-$ and $v_*$, we associate to it a pseudo-holomorphic curve $u_v:(\Sigma,\partial\Sigma)\rightarrow(X,L)$ for the almost complex structure $J$. While, if $v$ is a spine vertex between $v_*$ and $e_+$, we associate to it a  $J_{1,t}$-strip bounding $L$ and $\phi(L)$, as in Definition \ref{J-stripsDefn}, where $J_{1,t}$ is as given above. On the other hand, all internal vertices of the $0$-page of $\Gamma$ correspond to $J$-disks $u:(\Sigma,\partial\Sigma)\rightarrow (X,L)$. While for internal vertices of the $1$-page of $\Gamma$, we associate to them $(\phi_r)_*J$-disks where $r\in[0,1]$ is determined as above. Namely, if $v$ is a vertex on a branch attached to a spine vertex between $e_-$ and $v_*$ then we set $r=0$. If $v$ is a vertex on a branch attached to a spine vertex between $v_*$ and $e_+$ then we set $r=1$. Lastly, if $v$ is a vertex on branch attached to $v_*$ then we set $r$ to the same as given above.

Let $(H,\nabla,\tau)\rightarrow X$ be the Hermitian line bundle as in the proof of Proposition \ref{UnobstructedHamiltonianImage} or Lemma $3.6$ of \cite{rabah2024fukaya}. To be able to construct a framing for the moduli space as described in this proof, we need to interporlate our connection $\nabla$ along the edges and vertices of a given two-paged metric ribbon tree in a coherent way. To be more precise, let $\Gamma$ be a two-paged metric ribbon tree and denote by $u:\R_s\times[0,1]_t\rightarrow X$ the pseudo-holomorphic map associated to the distinguished spine vertex $v_*$. Set $v(s,t):=\phi^{-1}_{\chi(s)}\circ u(s,t)$. Using the $\nabla$-parallel transport along the curve $s\mapsto \phi_{\chi(s)}(x)$, we get a unitary isomorphism $U_{\chi(s)}:H_x\rightarrow H_{\phi_{\chi(s)}(x)}$ such that $U_0=\Id_{H_x}$. This induces a bundle isomorphism 
\[
I_u:v^*H\rightarrow u^* H,
\]
covering the identity map on $\R\times [0,1]$. Now let $\beta:[0,1]_t\rightarrow [0,1]_r$ be a smooth function such that $\beta\equiv 0$ near $0$ and $\beta\equiv 1$ near $1$ and denote by $\hat{\beta}:\R_s\times[0,1]_t\rightarrow [0,1]_r$ its pullback on the domain of $u$. We set 
\[
H_u:=u^*H,
\]
with the connection
\[
\nabla_u:=(1-\hat{\beta}(s,t))u^*\nabla+\hat{\beta}(s,t)I_u(v^*\nabla)I_u^{-1},
\]
and boundary trivializations given by 
\[
\tau_u(s,0):=\tau(u(s,0)),
\]
\[
\tau_u(s,1):=U_{\chi(s)}(\tau(v(s,1))).
\]
Having such construction in mind, we use $\phi_r^*(H,\nabla,\tau)\equiv(H_r,\nabla_r,\tau_r)$ to frame the associated psuedo-holomorphic curves at each vertex. Now using the arguments of Proposition $4.37$, Proposition $4.42$, Theorem $4.44$ and Theorem $4.48$ of \cite{rabah2024fukaya}, and using Theorem \ref{GKCStripMorseFloerThm} allows us to give the moduli space as described in this proof, a normally complex Global Kuranishi chart representation $(G,\T,\E,s)$ together with Global Kuranishi chart representation $(G_j,\T_j,\E_j,s_j)$ of the components of its codimension-$1$ such that, the following diagram commutes:

\[\begin{tikzcd}
	{\mathcal E_1\boxplus\mathcal E_2} && {\mathcal E} \\
	\\
	{\mathcal T_1\times \mathcal T_2} && {\mathcal T} \\
	\\
	{\mathcal B_1\times \mathcal B_2} && {\overline T_{k_0,k_1}}
	\arrow[hook, from=1-1, to=1-3]
	\arrow[from=1-1, to=3-1]
	\arrow[from=1-3, to=3-3]
	\arrow["{s_1\boxplus s_2}", bend left=30pt, , from=3-1, to=1-1]
	\arrow[hook, from=3-1, to=3-3]
	\arrow[from=3-1, to=5-1]
	\arrow["{s_{\underline x,p,q,\underline y,\beta}}"', bend right=30pt, , from=3-3, to=1-3]
	\arrow[from=3-3, to=5-3]
	\arrow["{\text{glue}}", from=5-1, to=5-3]
\end{tikzcd}\]
where $\mathcal B_j$ is as in the statement of Theorem \ref{GKCStripMorseFloerThm}. Lastly, using the same arguments of the proof of Proposition \ref{UnobstructedHamiltonianImage} to get the desired FOP-perturbation, we can apply the arguments of the proof of Theorem $6.1.25$ of \cite{fukaya2010lagrangian} to finish the proof.

\end{proof}

\section{Applications}\label{AppSection}
Having a Lagrangian Floer theory with values in the Novikov ring with integer coefficients, enable us to lift such modules, after some additional geometric conditions such as $\textit{polarization}$ of the ambient symplectic manifold relative to the Lagrangian for instance, to take values in ring spectrums (c.f. \cite{abouzaid2024foundation}). A lot of $\textit{intersecting}$ work along these lines have been done in the following non-exhaustive list \cite{porcelli2024spectral}, \cite{porcelli2024bordism}, \cite{porcelli2025open}, \cite{blakey2025parameterized} and \cite{rezchikov2024cyclotomic}.

\subsection{Integral Lagrangian Intersection Bounds}

In this subsection, we give lower bounds on the number of intersection points of an unobstructed Lagrangian $L$, as in Definition \ref{UnobstructedDefn}, with its Hamiltonian image $\phi(L)$. Using our Lagrangian Floer modules, as in Definitions \ref{LagFloerHomSingle L defn} and \ref{HF2Lags Defn}, we will be able to give stronger bounds opposed to the cases as in \cite{Fukaya2009} by incorporating torsion elements from the singular homology of the underlying Lagrangians (c.f. \cite{bai2022arnold} and \cite{rezchikov2022integral}).

Let $L\subset X$ be a closed, connected, relatively-spin Lagrangian such that $L$ is unobstructed over $\Lambda_0$, as in Definition \ref{UnobstructedDefn}. Fix $b$ a bounding cochain of $L$ as in Definition \ref{boundingCochainDefn}. Let $\phi:X\rightarrow X$ be a Hamiltonian diffeomorphism such that $L$ and $\phi(L)$ intersect transversely.

The goal of this subsection is to prove the following Theorems: 

\begin{theorem}[Integral Bound, Ungraded Case]\label{IntergralArnoldConjI}
    Suppose that the canonical inclusion $\iota:L\hookrightarrow X$ induces a split-injective group morphism $\iota_*:H_j(L;\Z)\hookrightarrow H_j(X;\Z)$ for every $j$. Then, the number of intersection points of $L$ and $\phi(L)$ is bounded from below by $\#(L\cap\phi(L))\geq \sum \mathrm{rk}_\Z H_*(L;\Z)+2t_L$, where $t_L$ is the minimum number of generators of the torsion group of $\bigoplus H_*(L;\Z)$.
\end{theorem}

\begin{theorem}[Integral Bound, Graded Case]\label{IntegralBoundGraded}
    Suppose that $2c_1(X)=0\in H^2(X;\Z)$ and $L\subseteq X$ is a graded Lagrangian in the sense of \cite{seidel2000graded}. Furthermore, suppose that the canonical inclusion $\iota:L\hookrightarrow X$ induces a split-injective group morphism $\iota_*:H_j(L;\Z)\hookrightarrow H_j(X;\Z)$ for every $j$. Then, the number of intersection points of $L$ and $\phi(L)$ is bounded from below by $\#(L\cap\phi(L))\geq \sum_j\mathrm{rk}_\Z H_j(L;\Z)+2t_j$, where $t_j$ is the minimum number of generators of the torsion group of $H_j(L;\Z)$.
\end{theorem}

\begin{remark}
    The key reason why we assume split-injective rather than just injective as in Theorem \ref{degenerationSpecSeq}, is that we cannot deduce as isomorphism at the level modules from Theorem \ref{degenerationSpecSeq} (c.f. \cite{fukaya2013lagrangian}).
\end{remark}

\subsubsection{Algebraic Preliminaries}
Before we give a proof of Theorems \ref{IntergralArnoldConjI} and \ref{IntegralBoundGraded}, we will need some algebraic results to get convergence of the spectral sequence as in Theorem \ref{SpecSeqThmI}. It is important to note that Theorem \ref{degenerationSpecSeq} only shows degeneration at the $E_2$-page, yet it does not provide us with an isomorphism as $\Lambda_0$-modules. To be precise, in this subsection we show the following:

\begin{prop}\label{PropConvergence}
    Under the setting of Theorem \ref{degenerationSpecSeq} and if we further assume that $\iota_*:H_*(L;\Z)\rightarrow H_*(X;\Z)$ admits a group morphism left inverse for every $*$. Then, for any integer $k$, we have $HF^k((L,b);\Lambda_0)\cong\bigoplus H_*(L;\Z)\otimes_\Z\Lambda_0^{k-*}$ as $\Lambda_0$-modules.
\end{prop}

\begin{lemma}\label{LemmaIArnold}
    Suppose that $\bar C, \bar B$ are freely finitely-generated, bounded chain complexes over $\Z$ and $H(f):H^*(\bar B)\rightarrow H^*(\bar C)$ is a degree-preserving group morphism. Then, there exists a chain map $f:\bar B\rightarrow \bar C$ inducing such morphism.
\end{lemma}

\begin{proof}
    Following the same notations as in the $\textit{Spectral Sequence}$ arguments above, we set $Z^n_{\bar B}:=\ker d^n_{\bar {B}}$ and $B^n_{\bar B}:=\mathrm{im}(d^{n-1}_{\bar {B}})$, and similarly for $\bar C$. As $\bar B^n/Z^n_{\bar B}$ is a free group, it follows that we have a (non-canonical) splitting $\bar B^n=Z^n_{\bar B}\oplus F^n_{\bar B}$ such that $F^n_{\bar B}$ is a free group. Denote by $q^n_{\bar B}:Z^n_{\bar B}\rightarrow H^n(\bar B)$ the quotient map and similarly for $\bar C$. As $\bar B^n$ is a free group, it follows that $Z^n_{\bar B}$ is a free group and hence there exists a lift of $H^n(f)$, $a_n:Z^n_{\bar B}\rightarrow Z_{\bar C}^n$ such that $q_{\bar C}^n\circ a_n\equiv H^n(f)\circ q_{\bar B}^n$. In particular, $a_n(B_{\bar B}^n)\subseteq B^n_{\bar C}$. Therefore, there exists a group morphism $a_{n+1}\circ d^n_{\bar B}|_{F^n_{\bar B}}:F^n_{\bar B}\rightarrow B^{n+1}_{\bar C}$. As $F^n_{\bar B}$ is a free group, it follows that there exists a group morphism $c_n:F^n_{\bar B}\rightarrow \bar C^n$ such that $d^n_{\bar C}\circ c_n\equiv a_{n+1}\circ d^n_{\bar B}|_{F^n_{\bar B}}$. Now define $f^n:\bar B^n\rightarrow \bar C^n$ using the above splitting by $f^n(x+y):=a_n(x)+c_n(y)$.
\end{proof}

\begin{lemma}\label{LemmaIIArnold}
    Let $f:B\rightarrow C$ be a filtered chain map such that $gr_p(f):F^pB/F^{p+1}B\rightarrow F^pC/F^{p+1}C$ is a quasi-isomorphism for every $p$. Then, $f$ is a quasi-isomorphism. 
\end{lemma}

\begin{proof}
    This follows from separatedness and completeness of $C, B$ and $\mathrm{Cone}(f)$ with respect to the topology induced from the $F$-filtration. Indeed, we give $\mathrm{Cone}(f)$ the induced direct-sum filtration and note that our assumption implies that $gr_p\mathrm{Cone}(f):=F^p\mathrm{Cone}(f)/F^{p+1}\mathrm{Cone}(f)$ is a acyclic complex for any $p$. It suffices to show that $\mathrm{Cone}(f)$ is acyclic. For instance, given a cocycle $z\in\mathrm{Cone}(f)$ of degree $k$, we will argue by induction using the $F$-filtration to show that there exists $y\in\mathrm{Cone}(f)$ of degree $k-1$ such that $dy=z$. The base case follows from the fact that $gr_p\mathrm{Cone}(f)$ is acyclic for any $p$. Now assume that we have constructed $y_1,\dots,y_p$ such that $y_i\in\mathrm{Cone}(f)/F^i\mathrm{Cone}(f)$ and $dy_i=z\text{ mod }F^i\mathrm{Cone}(f)$. As $gr_{p+1}\mathrm{Cone}(f)$ is acyclic, it follows that there exists $a\in gr_{p+1}\mathrm{Cone}(f)$ such that $da=z\text{ mod }F^{p+1}\mathrm{Cone}(f)$. Denote by $\pi:\mathrm{Cone}(f)/F^{p+1}\mathrm{Cone}(f)\rightarrow\mathrm{Cone}(f)/F^p\mathrm{Cone}(f)$ the projection map and note that $\pi(a)-y_p$ is closed in $\mathrm{Cone}(f)$. Once again, using acyclicity of $gr_{p+1} \mathrm{Cone}(f)$, it follows that there exists $v\in\mathrm{Cone}(f)/F^p\mathrm{Cone}(f)$ of degree $k-2$ such that $dv=\pi(a)-y_p$. Choose a lift $\tilde v\in\mathrm{Cone}(f)/F^{p+1}\mathrm{Cone}(f)$ and set $y_{p+1}:=a-d\tilde v$. Now by completeness and separatedness, it follows that $y:=\varprojlim_py_p\in\mathrm{Cone}(f)$ of degree $k-1$ such that $dy=z$. 
\end{proof}

\begin{proof}[Proof of Proposition \ref{PropConvergence}]
    Set $C:=\mathrm{CM}^*(f;\Lambda_0)$ and $B:=\mathrm{CM}^*(F;\Lambda_0)$. Our assumption provides us with left-inverses of the Gysin map $H^*(P_0):H^*(L;\Z)\rightarrow H^{*+n}(X;\Z)$. Lemma \ref{LemmaIArnold} realizes such inverse from a chain map $r_0:(\bar B, d_B)\rightarrow (\bar C, d^\text{Morse})$ where we are following the same notation as in Lemma \ref{OCfilterdChainMapSpec}. Set $\Phi:=(r_0\otimes \mathbf{1})\circ P$ where $P$ is as in Lemma \ref{OCfilterdChainMapSpec}. Noting that, $gr_p\Lambda_0$ is a free abelian group, using the above filtration, it follows that $gr_p\Phi=(r_0\circ P_0)\otimes\mathbf 1$ for all $p$. Now the result follows from Lemma \ref{LemmaIIArnold}. 
\end{proof}

Denote by 
\[
R:=\{\sum_{\lambda\in\R}a_\lambda T^\lambda:a_\lambda\in\Z, \#\{a_\lambda\neq0:\lambda<E\}<\infty, \forall E\in\R\}.
\]

\begin{lemma}[c.f. \cite{bai2022arnold}]\label{LemmaArnoldIII}
    Suppose that $A$ is a finitely-generated abelian group. Then, $\mathrm{rk}_R(A\otimes_\Z R)=\mathrm{rk}_\Z(A)$ and $t(A\otimes_\Z R)=t(A)$.
\end{lemma}

\begin{proof}
    This follows from the fact that $R$ is PID (c.f. \cite{hofer1995floer}) and the Fundamental Theorem of Modules over a PID.
\end{proof}

\begin{lemma}\label{LemmaArnoldIV}
    Let $C_0\xrightarrow{\partial_0}C_1\xrightarrow{\partial_1}C_0$ be a $2$-periodic complex of finite free $R$-modules. For $\epsilon\in\{0,1\}$ denote by $N_\epsilon:=\mathrm{rk}_R (C_\epsilon)$, $H_\epsilon:=\ker (\partial_\epsilon)/\mathrm{im}(\partial_{1-\epsilon})$, $h_\epsilon:=\mathrm{rk}_R(H_\epsilon)$ and $t_\epsilon:=t(H_\epsilon)$. Then, $N_\epsilon\geq h_\epsilon+t_\epsilon+t_{1-\epsilon}$.
\end{lemma}

\begin{proof}
    Set $Z_\epsilon:=\ker(\partial_\epsilon)$ and $B_\epsilon:=\mathrm{im}(\partial_{1-\epsilon})$. Then, we get the following short exact sequences
    \[
    0\rightarrow Z_\epsilon\rightarrow C_\epsilon\rightarrow B_{1-\epsilon}\rightarrow0,
    \]
    \[
    0\rightarrow B_\epsilon\rightarrow Z_\epsilon\rightarrow H_\epsilon\rightarrow 0.
    \]
    By exactness, it follows that $N_\epsilon=h_\epsilon+\mathrm{rk}_R(B_\epsilon)+\mathrm{rk}_R(B_{1-\epsilon})$. Now the result follows after noting that $t_\epsilon\leq\mathrm{rk}_R(B_\epsilon)$ for any $\epsilon\in\{0,1\}$.
\end{proof}

\begin{proof}[Proof of Theorem \ref{IntergralArnoldConjI}]
    Having Lemma \ref{UnobstructedHamiltonianImage} and Theorem \ref{HamiltonianInvariance} in mind, we argue as follows:

    As multiplication by $e$ on $\mathrm{CF}((L,b), (\phi(L), \phi_*(b));\Lambda)$ is invertible and since $\deg(e)=2$, it follows that we have the following $2$-periodic complex
    \[
    C_0\xrightarrow{d^\text{Floer}} C_1\xrightarrow{e^{-1}d^\text{Floer}}C_0
    \]
    where $C_\epsilon$ is generated over $R$ by elements of $\mathrm{CF}((L, b), (\phi(L), \phi_*(b));\Lambda)$ having index equal to $\epsilon$ after modding out by $2$. Now apply Lemma \ref{LemmaArnoldIII} and Lemma \ref{LemmaArnoldIV} on Proposition \ref{PropConvergence} to get the desired result. 
\end{proof}

\subsection{Lagrangian Intersection Bounds over Finite Fields}

In this subsection we will interested in giving a lower bound on the number of intersection points between two transversely intersecting Lagrangians without assuming that $\iota:L\hookrightarrow X$ the canonical inclusion of one of these Lagrangians, induces a split-injective map on homology. Our approach is the same as in \cite{Fukaya2009}, \cite{fukaya2010lagrangian} and \cite{fukaya2013lagrangian}, namely we will use a Lagrangian Floer theory over finite field $\Z/p$, where $p$ is a prime number, to obtain such bounds. 

Fix a prime $p$ and denote by $\mathbb F:=\Z/p$. For $L\subset X$ a closed, connected relatively-spin Lagrangian, as in Definition \ref{unobstructedoverFdefn}. We denote by $A_p:=\sum\dim_\mathbb FH_*(L;\mathbb F)$ and $B_p:=\sum\dim_\mathbb F\ker(\iota_*:H_*(L;\mathbb F)\rightarrow H_*(X;\mathbb F))$. The main goal of this subsection is to prove the following theorem:

\begin{theorem}[Theorem I \cite{Fukaya2009}, Theorem $15.6$ \cite{fukaya2010lagrangian}]\label{LagIntFThm}
    Suppose that $\phi:X\rightarrow X$ is a Hamiltonian diffeomorphism such that $L$ and $\phi(L)$ intersect transversely. Then, the number of intersection points $\#(L\cap\phi(L))$ is bounded from below by $\#(L\cap\phi(L))\geq A_p-2B_p$.
\end{theorem}

We denote by $\Lambda_0^\mathbb F:=\{\sum_ia_iT^{\lambda_i}e^{n_i}:a_i\in\mathbb F, n_i\in\Z,0\leq\lambda_1\leq\lambda_2\leq\dots,\lim_i\lambda_i=\infty\}$, the Novikov ring over $\mathbb F$. We show the following easy Lemma:

\begin{lemma}\label{LemmaReductionModp}
    Reducing the coefficients of $\Lambda_0$ $\text{mod }p$, defines a continuous graded surjective ring morphism $\Lambda_0\rightarrow \Lambda_0^\mathbb F$. Moreover, we have a ring isomorphism $\Lambda_0^\mathbb F\cong\Lambda_0/p\Lambda_0$.
\end{lemma}

\begin{proof}
    Consider $\rho:\Lambda_0\rightarrow\Lambda_0^\mathbb F$ given by $\sum_ia_iT^{\lambda_i}e^{n_i}\mapsto\sum_i\bar a_iT^{\lambda_i}e^{n_i}$ where $a_i\equiv\bar a_i\text{ mod }p$. Observe that by finiteness of $\{a_i\neq0:\lambda_i<E, \forall E\in\R\}$, it follows that $\rho$ is continuous, with respect to the $T$-adic topology, and grading-preserving. Thus by completeness, it follows that $\rho$ is surjective and by separatedness that $\ker\rho=p\Lambda_0$.
\end{proof}

Let $(f,g)$ be a Morse-Smale pair as above (c.f. \cite{rabah2024fukaya}), and denote by $C:=\mathrm{CM}^*(f;\Lambda_0)$ and by $C_p:=\mathrm{CM}^*(f;\Lambda_0^\mathbb F)=\mathrm{CM}^*(f;\mathbb F)\widehat{\otimes}_\mathbb F\Lambda_0^\mathbb F$. Moreover, we denote by $\rho_C:C\rightarrow C_p$ the induced map from Lemma \ref{LemmaReductionModp} and mod $p$ reduction of elements $\mathrm{Crit}(f)$.

\begin{theorem}\label{A_inftyoverF}
    The Morse complex $C_p$ admits the structure of a curved filtered $A_\infty$-algebra over $\Lambda_0^\mathbb F$ given by $\{\m^{(p)}_k\}_{k\geq0}$, where for any $x_1,\dots,x_k\in\mathrm{Crit}(f)$ and if we denote by $\bar x_1,\dots,\bar x_k$ their mod $p$ reduction then, $\rho_C\circ\m_k(x_1,\dots, x_k)=\m_k^{(p)}(\bar x_1,\dots, \bar x_k)$.
\end{theorem}

\begin{proof}
    The key point is to use the same $\textit{perturbation}$ data as in \cite{rabah2024fukaya}, which is used in the construction of the $A_\infty$-algebra structure on $C$ and then apply reduction $\text{mod }p$ on the pseudo-cycle obtained from the zero-locus of the single-valued FOP perturbations used in \cite{rabah2024fukaya} and adhere to the results of Appendix $A$ and Section $2.2.1$ of \cite{bai2024cohomological}.

    $\textbf{Step 1: }\textit{Construction of an }A_{n,K}$-$\textit{algebra}$.

Let $x_0,x_1,\dots, x_k\in\mathrm{Crit}(f)$, $\beta\in\pi_2(X,L)$ and $\tilde f\in C^\infty(\overline{\mathrm{Gr}}_{k+1}\times L)$ be a generic domain-dependent Morse function associated to $f$, as in Theorem $4.50$ of \cite{rabah2024fukaya}. Denote by $(G_\beta,\T(\underline{x},\beta),\E(\underline{x},\beta), s_{\underline{x},\beta})$ the Global Kuranishi chart of $\M_{k+1}(\underline{x};\beta,\tilde f)$ given by Theorem $4.50$ of \cite{rabah2024fukaya}. After taking quotient by $G$ and using the fact that $G$ acts linearily and with finite-isotropy on both $\T(\underline{x}, \beta)$ and $\E(\underline{x},\beta)$, we consider its associated derived orbifold chart $(\mathcal{U},E,\mathcal S)$. Let $s_\epsilon$ be the same FOP-perturbation of $\mathcal S:\mathcal U\rightarrow E$ as in the first paragraph of the proof of Theorem $5.3$ of \cite{rabah2024fukaya} and denote by $\overline{|s_\epsilon^{-1}(0)\cap \mathcal{U}_{free}|}$ the $\text{mod }p$ reduction of the pseudo-cycle used in the proof of Theorem $5.3$ of \cite{rabah2024fukaya}. Set 

\[
\m_{k,\beta}^{(p)}(\bar x_1,\dots, \bar x_k):=\sum_{x_0\in\mathrm{Crit}(f)}(-1)^*\overline{|s_\epsilon^{-1}(0)\cap \mathcal{U}_{free}|}\bar x_0,
\]
where the expected dimension $\dim_\R\M(\underline{x};\beta,\tilde f)=0$ and $*$ is as in Theorem $5.3$ of \cite{rabah2024fukaya}.
Moreover, we set
\[
\m_k^{(p)}:=\sum_{\beta\in\Pi(L)}\m_{k,\beta}^{(p)}T^{\omega(\beta)}e^\frac{\mu(\beta)}{2}.
\]
Now after using Lemma \ref{LemmaReductionModp} to deduce that $\rho_C:C\rightarrow C_p$ is continuous with respect to the $T$-adic topology, the claim that $\m_k^{(p)}$ defines a structure of an $A_{n, K}$-algebra on $C_p$ follows from our construction and Theorem $5.3$ of \cite{rabah2024fukaya}.

$\textbf{Step 2: }\textit{Extending to an }A_\infty$-$\textit{algebra}$.

In principle, as we are working over Novikov rings with coefficients in $\mathbb F$, one can adhere to the $\textit{homotopy limit}$ results found in Section $13.3$ of \cite{fukaya2013lagrangian} or Theorem $7.2.72$ of \cite{fukaya2010lagrangian}. Nevertheless, in order to prove the desired equality and following the construction in $\textbf{Step 1}$, we use the same arguments in the proof of Corollary $5.4$ of \cite{rabah2024fukaya}, to deduce our claim.

$\textbf{Step 3: }\rho_C\circ\m_k(x_1,\dots, x_k)=\m_k^{(p)}(x_1,\dots,x_k)$.

From $\textbf{Step 1}$ and $\textbf{Step 2}$, it follows that for any $E\in\R$, we have 
\[
\rho_C\circ\m_k(x_1,\dots, x_k)\equiv\m_k^{(p)}(\bar x_1,\dots, \bar x_k)\text{ mod }T^E\Lambda_0^\mathbb F.
\]
Thus by completeness and separatedness, it follows that such equality holds at the chain level and hence the desired result.
\end{proof}

\begin{defn}[Unobstructed over $\Lambda_0^\mathbb F$]\label{unobstructedoverFdefn}
    $L$ is said to be unobstructed over $\Lambda_0^\mathbb F$ if there exists a $b\in C_p$ of degree-$1$ and of positive-valuation such that 
    \[
    \sum_{k\geq0}\m_k^{(p)}(b,\dots,b)=0.
    \]
\end{defn}

Similarly as above, we will call such $b$ a boundig cochain of $(C_p,\{\m_k^{(p)}\}_{k\geq0})$ and as in Definition \ref{deformedAinftyDefn}, we can deform our $A_\infty$-operations over $\Lambda_0^\mathbb F$ using such $b$.

\begin{defn}[Lagrangian Floer Theory over $\Z/p$]
    For a closed connected relatively-spin Lagrangian that is unobstructed over $\Lambda_0^\mathbb F$, as in Definition \ref{unobstructedoverFdefn}, and a choice of bounding cochain $b\in C_p$, we define 
    \[
    HF((L,b);\Lambda_0^\mathbb F):=\ker(\m_1^{(p), b})/\mathrm{im}(\m_1^{(p), b}),
    \]
    where $\m_k^{(p), b}$ is as given in Definition \ref{deformedAinftyDefn}.
\end{defn}

\begin{prop}\label{UnobstroverZimpliesF}
    If $L$ is unobstructed over $\Lambda_0$, as in Definition \ref{UnobstructedDefn}, then $L$ is unobstructed over $\Lambda_0^\mathbb F$, as in Definition \ref{unobstructedoverFdefn}.
\end{prop}

\begin{proof}
    Suppose that $L$ is unobstructed over $\Lambda_0$ and let $b\in C$ be a bounding cochain of $(C, \{\m_k\}_{k\geq0})$. Then, by completeness and separatedness of $C$ and $C_p$ with respect to the $T$-adic topology, and by continuity of $\rho_C:C\rightarrow C_p$, it follows that $\rho_C(b)$ satisfies the $\textit{Maurer-Cartan}$ equation for $(C_p, \{\m_k^{(p)}\}_{k\geq0})$ and hence the desired result.
\end{proof}

We can give a partial converse of the above Proposition.

\begin{prop}\label{UnobstructedoverFimpliesZ}
    Given $L$ a closed connected relatively-spin monotone Lagrangian of $X$. That is, $\omega(\beta)=\tau\mu_L(\beta)$ for any $\beta\in\pi_2(X,L)$ where $\tau>0$. Suppose that, $L$ is unobstructed over $\Lambda_0^{\Z/p}$, as in Definition \ref{unobstructedoverFdefn}, for all prime $p$. Then, $L$ is unobstructed over $\Lambda_0$, as in Definition \ref{UnobstructedDefn}.
\end{prop}

\begin{proof}
    Let $(f,g)$ be a Morse-Smale pair on $L$ with the condition that $f$ has a unique minimum $x_0\in\mathrm{Crit}(f)$. That is, $\mathrm{CM}^0(f;\Z)=\Z\langle x_0\rangle$. Note that, for $\beta\in\pi_2(X,L)$, $\m_{0,\beta}(1)\in\mathrm{CM}^{2-\mu_L(\beta)}(f;\Z)$ and hence we only need to consider Maslov-index $2$ $J$-disks contributions. In particular and using our monotonicity assumption, we can write $\m_0=wT^Ae$ such that $A=2\tau$. Let $b=\sum_ib_iT^{\lambda_i}e^{n_i}$ be a degree-$1$ element of positive-valuation. Thus, $\deg b_i=1+2n_i$, where $n_i\geq0$ and $\lambda_i>0$. Thus, the coefficient of $e$ in the expansion of $\m_k(b,\dots,b)$ will be from $J$-disks of degree $\beta$ such that $\mu_L(\beta)/2-\sum n_i=1$. On the other hand, the energy of such configurations will be $\omega(\beta)+\sum\lambda_i$ and using monotonicity, it follows that $\omega(\beta)+\sum\lambda_i>2\tau=A$. Therefore, the expansion of $\m_k(b,\dots, b)$ has no terms having $T^Ae$ as a coefficient. Now under the assumption that $L$ is unobstructed over $\Lambda_0^{\Z/p}$ for all prime $p$, it follows that $w\equiv 0\text{ mod }p$ for all prime $p$ and hence $L$ is unobstructed over $\Lambda_0$ having $b=0$ as a bounding cochain. 
\end{proof}

\begin{remark}
    We do not know how to prove Proposition \ref{UnobstructedoverFimpliesZ} in full generality. $\textit{In fact, we expect it to be false}$. Namely, the argument we presented above, applied on the non-monotone case, boils down to saying that any polynomial of integer coefficients that has a solution over $\Z/p$ for all prime $p$, has an integer solution. This is false in general. For instance, we can give the following $\textit{adhoc}$ algebraic example: 
    
    Consider the curved dg-algebra $C_p:=\Lambda_0^\mathbb F\mathbf{1}\oplus\Lambda_0^\mathbb Fa\oplus\Lambda_0^\mathbb Fc$ where $a,c$ are degrees $1,2$ respectively. The differential $d$ on $C_p$ is given by $da=4c$ and $dc=0$, while its curvature is given by $h=2T^Ec$ such that $E>0$. We solve the Maurer-Cartan Equation for $b=xa$ then, we will have $h+db+b^2=(2T^E+4x)c=0$. Thus, the coefficient of $T^E$ should satisfy $2x=-1$, which has no solution over $\Z$. On the other hand, for $p=2$, $x=0$ is a solution, while for $p$ odd prime $x=-\frac{p-1}{2}\text{ mod }p$ is a solution. Nevertheless, we do not know of any geometric counterexample of the converse of Propostion \ref{UnobstroverZimpliesF}.
\end{remark}

\begin{prop}\label{HamiltonianInvarianceF}
    Suppose that $L$ is unobstructed over $\Lambda_0^\mathbb F$, $b$ is a bounding cochain of $L$ and $\phi:X\rightarrow X$ is a Hamiltonian diffeomorphism such that $L$ and $\phi(L)$ intersect transversely. Then, we have $HF((L,b);\Lambda^\mathbb F)\cong HF((L,b), (\phi(L), \phi_*(b));\Lambda^\mathbb F)\cong HF((\phi(L), \phi_*(b));\Lambda^\mathbb F)$ as $\Lambda^\mathbb F$-vectorspaces. 
\end{prop}

\begin{proof}
    This follows from the arguments of the proof of Theorem \ref{HamiltonianInvariance} together with Lemma \ref{LemmaReductionModp}.
\end{proof}

\begin{proof}[Proof of Theorem \ref{LagIntFThm}]
    Given $L$ a closed connected relatively-spin Lagrangian which is unobstructed over $\Lambda_0^\mathbb F$ and let $b\in C_p$ be a bounding cochain of $(C_p, \{\m_k^{(p)}\}_{k\geq0})$. Let $(F, g_X)$ be a Morse-Smale pair on $X$, extending $(f,g)$ as above. Following the same notation as above, we set $B:=\mathrm{CM}^{*+n}(F;\mathbb F)\widehat{\otimes}_\mathbb F\Lambda_0^\mathbb F$ and $d_B:=(-1)^nd^\text{Morse}_X\otimes\mathbf{1}$. We follow the same Spectral Sequence constructions as above. Namely, for $k,q,r$ non-negative integers we set
    \[
    Z^{k,q}_r:=(F^qC^k_p\cap d^{-1}F^{q+r-1}C^{k+1}_p)+F^{q+1}C^k_p,
    \]
    \[
    B^{k,q}_r:=(dF^{q-r+2}C_p^{k-1}\cap F^qC^k_p)+F^{q+1}C^k_p,
    \]
    \[
    E^{k,q}_r(C_p):=Z^{k,q}_r/B_r^{k,q},
    \]
    with the understanding that $d\equiv\m_1^{(p),b}$. Once agian, as $\m_1^{(p), b}\circ \m_1^{(p), b}=0$, we have $B^{k,q}_r\subset Z^{k,q}_r$. The differential $\delta_r$ on $E_r(C_p)$ is given on $x\in F^qC_p^k$ such that $dx\in F^{q+r-1}C^{k+1}_p$ by $\delta_r[x]:=[dx]\in E^{k+1, q+r-1}_r(C_p)$. Using the same choice of $\eta$ to define our filtration as above, we have
    \[
    E^{k,q}_1(C_p)=\bigoplus\mathrm{CM}^*(f;\mathbb F)\otimes_\mathbb F (gr_q\Lambda_0^\mathbb F)^{k-*},
    \]
    \[
    E_2^{k,q}(C_p)\cong\bigoplus H^*(L;\mathbb F)\otimes_\mathbb F (gr_q\Lambda_0^\mathbb F)^{k-*}
    \]
    and $\delta_1\equiv d^\text{Morse}$.
    We do the same thing for the filtered complex $(B, d_B)$ and note that its $E_2$-page is $H^{*+n}(X;\mathbb F)$ and its differential $\delta^B_r\equiv 0$ for all $r\geq2$.

    Using the same train of thought as in the proof of Theorem \ref{A_inftyoverF}, applied on the proof of Theorem \ref{TruncatedOCmap} and Theorem \ref{OCmapTheorem}, we get a filtered chain map
\[
P_b:(C_p,\m_1^{(p), b})\rightarrow (B, d_B).
\] 
Moreover, using the same arguments as above, such $P_b$ induces a chain map on the energy-filtration associated spectral sequences 
\[
P_r:E_r(C_p)\rightarrow E_r(B),
\]
such that $P_2=\iota_!\otimes \Id$ where $\iota_!:H^{*}(L;\mathbb F)\rightarrow H^{*+n}(X;\mathbb F)$ is the Gysin map. Moreover, as $P_r:E_r(C_p)\rightarrow E_r(B)$ is a chain map, it sends $Z_r(C_p)$ to $Z_r(B)$ and similary for the $B_r$ submodules. In addition, we have $P_r\circ\delta_r=\delta_r^B\circ P_r$ and $P_{r+1}=H(P_r)=\ker(P_r)/\mathrm{im}(P_{r})$. Now denote by $K_r$ defined inductively by setting
\[
K_{r+1}:=(K_r\cap\ker(\delta_r))/\mathrm{im}(\delta_r)
\]
and $K_2:=\ker(P_2)$. Note that by Poincare duality, it follows that $\mathrm{rk}_\mathbb F(\iota_!)=\mathrm{rk}_\mathbb F(\iota_*)=A_p-B_p$. With this in mind and using the same algebraic arguments as in Theorem $15.3$ of \cite{fukaya2010lagrangian}, it follows that $\dim_\mathbb F HF((L,b);\Lambda^\mathbb F)\geq A_p-2B_p$. Now the result follows from Proposition \ref{HamiltonianInvarianceF}.
\end{proof}

\begin{cor}\label{corLagIntoverF}
    If $L$ is unobstructed over $\Lambda_0$ and $\phi:X\rightarrow X$ is a Hamiltonian diffeomorphism such that $L$ and $\phi(L)$ intersect transversely then, $\#(L\cap\phi(L))\geq\max_{p \text{ prime}}(A_p-2B_p)$.
\end{cor}

\begin{proof}
    This follows from Proposition \ref{UnobstroverZimpliesF} and Theorem \ref{LagIntFThm}.
\end{proof}

\subsection{Examples}
In this section we give some examples and computations.

\subsubsection{Integral Lagrangian spheres in Symplectic Calabi-Yau manifolds}

For $0<\epsilon<<1$ consider the smooth quintic $3$-fold
\[
X_\epsilon:=\{z_0^5+\dots+z_4^5-(5+\epsilon)z_0\dots z_4=0\}\subset\CP^4.
\]
Using Theorem $2.7$ or Equation $(3.7)$ in Section $3$ of \cite{smith2002symplectic}, we can view $X_\epsilon$ as a smoothing of the $\textit{symplectic conifold}$ $X_0$. More importantly, it follows from Theorem $2.7$ of \cite{smith2002symplectic} that there exists an embedded Lagrangian sphere $L\cong S^3\subset X_\epsilon$. By adjunction, it follows that $c_1(X_\epsilon)=0$. Moreover, as $\pi_1(L)=0$, it follows that the map by the long exact sequence on homotopy groups of the pair $(X_\epsilon,L)$, $\pi_2(X_\epsilon)\rightarrow\pi_2(X_\epsilon,L)$ is surjective. Therefore, $\mu_L(\beta)=2c_1(X_\epsilon)[A]=0$ where $A\mapsto\beta$ under the above map. Thus, Corollary \ref{unobstructedSpheres} tells us that $L$ is unobstructed over $\Lambda_0$. Now let $f:L\rightarrow\R$ be a perfect Morse function and $\textit{generically}$ choose a metric $g$ on $L$ so that $(f,g)$ is a Morse-Smale pair. Write $\mathrm{CM}^*(f;\Z)=\Z\langle x_0\rangle\oplus\Z\langle x_3\rangle$ where $x_i$ has Morse index $i$. As $\mu_L\equiv0$ on $\pi_2(X_\epsilon, L)$, it follows that for any $\beta\in\pi_2(X_\epsilon,L)$ and integer $k\geq0$, we have $\deg(\m_{k,\beta})=2-k$. In particular, $\m_0\equiv 0$ as $f$ has no critical points of Morse index $2$. Thus, setting $b\equiv 0$ is a bounding cochain for the $A_\infty$-algebra $(\mathrm{CM}^*(f;\Lambda_0),\{\m_k\}_{k\geq0})$. Moreover, as $\deg(\m_1)=1$, it follows that $\m_1\equiv0$. Therefore, $HF^*(L;\Lambda_0)\cong H^*(S^3;\Z)\widehat{\otimes}_\Z\Lambda_0$ and in particular, $L$ is not Hamiltonianly displaceable.

The above discussion can be generalized to any integral Lagrangian sphere $L$ in a symplectic Calabi-Yau manifold $X$ of real dimension $2n$ at least $6$. Namely, we will have that $L$ is unobstructed over $\Lambda_0$ and if we choose the perfect Morse function on $L$ to generate our $A_\infty$-algebra, we can choose $b\equiv 0$ as a bounding cochain. Moreover, by the same $\textit{degree}$ argument, we can show that $\m_1\equiv0$ and hence, $HF^*(L;\Lambda_0)\cong H^*(S^n;\Z)\widehat{\otimes}_\Z\Lambda_0$. In another words, $L$ is an $\textit{essential}$ object in the $\textit{integral Fukaya category}$ of $X$ (c.f. \cite{RabahIntegralFukaya}). \\

\subsubsection{Rational Lagrangian Homology Sphere in a Calabi-Yau 3-fold}\label{Example}

We follow Example $2.17$ of \cite{OguisoSakurai2001}.

Identify $\C\cong\R\oplus i\R$ and consider the elliptic curve $E:=\C/(\Z\oplus i\Z)$. Let $X:=E^3/G$ where $G\cong\Z/2\oplus\Z/2$, acting on $E^3$ by:
\[
g_1.(z_1,z_2,z_3):=(z_1+\frac{1}{2}, -z_2,-z_3),
\]
\[
g_2.(z_1,z_2,z_3):=(-z_1,z_2+\frac{1}{2},-z_3+\frac{1}{2}),
\]
where $g_1,g_2$ are the generators of $\Z/2\oplus\Z/2$. Note that, $G$ acts freely and holomorphically on $E^3$. Moreover, if we denote by $\omega:=\sum_{j=1}^3dx_j\wedge dy_j$ then $G$ preserves $\omega$. Now set $L:=(\R/\Z)^3/G$ and note that $L$ is a closed, orientable (and hence spin as $\dim_\R L=3$) Lagrangian submanifold of $X$. Moreover, following Section $4$ of \cite{DekimpePetrosyan2009} we have the following:
\[
H_0(L;\Z)\cong H_3(L;\Z)\cong\Z,
\]
\[
H_1(L;\Z)\cong (\Z/4)^2,
\]
\[
H_2(L;\Z)=\{0\}.
\]
To show that $L$ is unobstructed over $\Lambda_0$, we note that $\C^3\rightarrow X$ is the universal cover and hence, $\pi_2(X)=\{0\}$. On the other hand, as $\pi_1(L)\rightarrow\pi_1(X)$ is injective, it follows that $\pi_2(X,L)=\{0\}$. In particular, for any choice of Morse-Smale pair on $L$, we have $\m_0\equiv 0$ and thus $b\equiv0$ serves as a bounding cochain. Now observe that, as $G$ acts by affine transformations with translations in the real-direction, it follows that the map $r:X\rightarrow L$ given by $[z_1,z_2,z_3]\mapsto[\mathrm{Re}z_1,\mathrm{Re}z_2,\mathrm{Re}z_3]$ satisfies $r\circ\iota=\Id_L$ where $\iota:L\hookrightarrow X$ is the canonical inclusion. Therefore, $\iota$ induces a split-injective morphism on homology and hence
\[
HF^*((L,0);\Lambda_0)\cong H^*(L;\Z)\otimes_\Z\Lambda_0.
\]
Same thing is true over $\Lambda$, namely $HF^*((L,0);\Lambda)\cong H^*(L;\Z)\otimes_\Z\Lambda$. In particular, given a Hamiltonian diffeomorphism $\phi:X\rightarrow X$ such that $L$ and $\phi(L)$ intersect transversaly, it follows that $\#(L\cap\phi(L))\geq 6$, compared to having $2$ as a lower bound when working over Novikov rings over fields with characteristic zero. 

\subsubsection{A Lagrangian in a non-spherically positive Symplectic Manifold with Torsion Components}

We will follow the above Example \ref{Example} closely yet with a minor change of notations.

Let $E$ and $G$ be as above, and set $X_0:=E^3/G$ and $L_0:=(\R/\Z)^3/G$. Denote by $q:E^3\rightarrow X_0$ the quotient map. To get an example that does not fall into the context of \cite{fukaya2013lagrangian}, we will blow-up $X_0$ away from $L_0$ (c.f. \cite{venugopalan2020fukaya}). To be more precise, let $c:=i/4\in E$ and set
\[
C_1:=q(E\times\{0\}\times\{c\}),
\]
\[
C_2:=q(\{0\}\times E\times\{c\}).
\]
Note that, each of $C_i\subset X_0$ is an embedded elliptic curve disjoint from $L_0$. Moreover, $C_1$ and $C_2$ intersect transversely at $p:=q(0,0,c)\in X_0$. Let,
\[
\pi_1:X_1\equiv\mathrm{Bl}_{C_1}(X_0)\rightarrow X_0,
\]
\[
\pi_2:X\equiv\mathrm{Bl}_{C_2'}(X_1)\rightarrow X_1,
\]
where $C_2'$ is the strict transform of $C_2$ by $\pi_1$, and set $\pi:=\pi_2\circ\pi_1$. As $L_0$ is disjoint from $C_1, C_2$, we can choose the Kahler blow-up parameters, that is the area of the exceptional divisors, so that the Kahler form $\tilde\omega$ on $X$ is identically equal to $\omega$, the Kahler form on $X_0$ as given in Example \ref{Example}, on a small enough tubular neighborhood of $L$, the proper transform of $L_0$ in $X$ (c.f. \cite{venugopalan2020fukaya}). Moreover, for such choices, we have $L\cong L_0$ are diffeomorphic. Such construction is the standard Atiyah flop (c.f. \cite{charest2015floer}) and in particular, if we denote by $F_0:=\pi_1^{-1}(p)\cong\CP^1$ and $F$ its strict transform then, 
\[
\nu_X(F)\cong\mathcal O_{\CP^1}(-1)\oplus\mathcal O_{\CP^1}(-1),
\]
where $\nu_X(F)$ is the normal bundle of $F\subset X$. Noting that, $\mathrm{GW}_{0,[F]}^X\neq 0$ (c.f. Section $2$ of \cite{rabah2024fukaya}), it follows by the independence of the Gromov-Witten invariants on the choice of almost complex structure, that $X$ does not admit a spherically-positive almost complex structure in the sense of Definition $1.1$ of \cite{fukaya2013lagrangian}. 

We give $X$ its integrable blow-up complex structure $J$ and observe that every $J$-disk $u:(\Sigma,\partial\Sigma)\rightarrow (X,L)$ is a constant map. This follows from the same argument as in the proof of Lemma $2.11$ of \cite{rabah2024fukaya}. Indeed, for such choice of almost complex structure, we have that $\pi:X\rightarrow X_0$ is a holomorphic map and hence, by the arguments of Example \ref{Example}, it follows that $\pi\circ u:(\Sigma,\partial\Sigma)\rightarrow (X_0,L_0)$ is a constant map. Therefore, $L$ is unobstructed over $\Lambda_0$ and moreover, after fixing any Morse-Smale pair on $L$, we can take $b\equiv 0$ as a bounding cochain. Once again, $r:X\rightarrow L$ given by $[z_1,z_2,z_3]\mapsto [\mathrm{Re}z_1,\mathrm{Re}z_2,\mathrm{Re}z_3]$ satisfies $r\circ\iota=\Id_L$ and hence the canonical inclusion $\iota:L\hookrightarrow X$ induces a split-injective morphism $\iota_*:H_*(L;\Z)\hookrightarrow H_*(X;\Z)$. Therefore,
\[
HF^*((L,0);\Lambda_0)\cong H^*(L;\Z)\otimes_\Z\Lambda_0,
\]
\[
HF^*((L,0);\Lambda)\cong H^*(L;\Z)\otimes_\Z\Lambda.
\]
Thus, given any Hamiltonian diffeomorphism $\phi:X\rightarrow X$ such that $L$ and $\phi(L)$ intersect transversely, it follows that
\[
\#(L\cap\phi(L))\geq 6,
\]
serving as a lower bound strictly greater than the one given over any field of charactersitic zero.

    \bibliography{ref}
\bibliographystyle{alpha}
\nocite{rabah2025fukaya}
\bigskip
\noindent
Mohamad Rabah\\
Morningside Center of Mathematics,
Chinese Academy of Sciences.\\
\textit{E-mail address:} \texttt{mohamad@amss.ac.cn},
\texttt{mohammadbrabah@gmail.com}
		
\end{document}